\documentclass{nyjm}
 \newtheorem{thm}{Theorem}[section]
 
 \newtheorem{lem}[thm]{Lemma}
 
 \theoremstyle{definition}
 
 \newtheorem{defn}[thm]{Definition}
 \theoremstyle{remark}
   
    \newtheorem{remark}[thm]{Remark}

 \newtheorem{notation}[thm]{Notation}
 \newtheorem{rem}[thm]{Remark}
  \newtheorem*{ack}{Acknowledgements}
 
 \numberwithin{equation}{section}

\usepackage{rotating}

\newcommand\cA{{\mathcal A}}

\newcommand{\cL}{{{\mathcal{L}}}}
\newcommand{\cC}{{{\mathcal{C}}}}

\newcommand{\tropical}[1]{\left[\,#1\,\right]}

\newcommand{\MCG}{{\operatorname{MCG}}}

\newcommand{\sgn}{{\operatorname{sgn}}}
\usepackage[utf8]{inputenc}
\usepackage[T1]{fontenc}
\usepackage{pinlabel}
\usepackage[pdftex,dvipsnames]{xcolor}
\usepackage{enumitem}
\newlist{subquestion}{enumerate}{1}
\setlist[subquestion,1]{label=(\alph*)}
\usepackage{amssymb}
\usepackage[noend]{algpseudocode}
\errorcontextlines\maxdimen

\makeatletter
    \newcommand*{\algrule}[1][\algorithmicindent]{\makebox[#1][l]{\hspace*{.5em}\thealgruleextra\vrule height \thealgruleheight depth \thealgruledepth}}%
\newcommand*{\thealgruleextra}{}
\newcommand*{\thealgruleheight}{.75\baselineskip}
\newcommand*{\thealgruledepth}{.25\baselineskip}

\newcount\ALG@printindent@tempcnta
\def\ALG@printindent{%
    \ifnum \theALG@nested>0% is there anything to print
        \ifx\ALG@text\ALG@x@notext% is this an end group without any text?
        \else
            \unskip
            \addvspace{-5.5 pt}% FUDGE to make the rules line up
            \ALG@printindent@tempcnta=1
            \loop
                \algrule[\csname ALG@ind@\the\ALG@printindent@tempcnta\endcsname]%
                \advance \ALG@printindent@tempcnta 1
            \ifnum \ALG@printindent@tempcnta<\numexpr\theALG@nested+1\relax% can't do <=, so add one to RHS and use < instead
            \repeat
        \fi
    \fi
    }%
\usepackage{etoolbox}

\usepackage[utf8]{inputenc}
\usepackage[english]{babel}

\usepackage{geometry}
\usepackage{marginnote}
\patchcmd{\ALG@doentity}{\noindent\hskip\ALG@tlm}{\ALG@printindent}{}{\errmessage{failed to patch}}
\makeatother

\newbox\statebox
\newcommand{\myState}[1]{%
    \setbox\statebox=\vbox{#1}%
    \edef\thealgruleheight{\dimexpr \the\ht\statebox+1pt\relax}%
    \edef\thealgruledepth{\dimexpr \the\dp\statebox+1pt\relax}%
    \ifdim\thealgruleheight<.75\baselineskip
        \def\thealgruleheight{\dimexpr .75\baselineskip+1pt\relax}%
    \fi
    \ifdim\thealgruledepth<.25\baselineskip
        \def\thealgruledepth{\dimexpr .25\baselineskip+1pt\relax}%
    \fi
 
    \State #1%
      \def\thealgruleheight{\dimexpr .75\baselineskip+1pt\relax}%
    \def\thealgruledepth{\dimexpr .25\baselineskip+1pt\relax}%
}

\begin{document}

\title[Curves  on non--orientable  surfaces and crosscap transpositions   ]
 {Curves  on non--orientable  surfaces and crosscap transpositions    }

%his generalizes (and simplifies) the results of a previous paper [3].

%----------Author 1
\author[\"{O}yk\"{u} Yurtta\c{s}]{\"{O}yk\"{u} Yurtta\c{s}}

\address{Dicle University Science Faculty Mathematics Department\\21280, Diyarbak{\i}r, Turkey}



\begin{abstract}

Let  $N_{g,n}$ be an $n$--punctured non--orientable surface of genus $g$ with one boundary component. For $g\geq 2$ one of the generators of the mapping class group of $N_{g,n}$ is a crosscap transposition.  We give explicit formulae for the action of crosscap transpositions and their inverses  on the set  of multicurves in $N_{g,n}$ in terms of  generalized Dynnikov coordinates.

\end{abstract}
%%% ----------------------------------------------------------------------
\maketitle

\section{Introduction}
 
%(g,n)\neq (2,0)$

%($n\geq 1, (g,n)\neq (1,1)$) 

 Let $N_{g,n}$~ ($g\geq 2$) be  a non--orientable surface of genus $g$ with $n$ punctures and one boundary component.    In all figures of this paper  each disk with a cross represents a crosscap, a graphical representation of a M\"{o}bius band. This means that  the interior of each such disc is removed, and the antipodal points on the resulting boundary component are identified. Throughout, we take a standard model of  $N_{g,n}$ where the punctures and the crosscaps are arranged  along the horizontal diameter of $N_{g,n}$ as shown in Figure \ref{arcsproof}.  A simple closed curve in $N_{g,n}$ is  {\it{inessential}} if it bounds an unpunctured disk, once punctured disk or an unpunctured annulus.  It is {\it{essential}}, otherwise. If a regular neighborhood of an essential simple closed curve in  $N_{g,n}$ is an annulus it is called  {\it{2-sided}}, and if it is a M\"{o}bius band it is called {\it{1-sided}}.  We call  the core curves and  the double covers of  the core curves  M\"{o}bius curves.  A {\it{multicurve}} in $N_{g,n}$ is a disjoint union of finitely many essential simple closed curves in $N_{g,n}$ modulo isotopy.  We write $\mathfrak{L}_{g,n}$ to denote the set of multicurves in $N_{g,n}$.

Multicurves on orientable surfaces are usually described by techniques such as the Dehn–Thurston coordinate system \cite{penner}. An alternative way to describe multicurves on finitely punctured disks is to use the Dynnikov coordinate system \cite{thurstonw}.  In 2016, Papadopoulos and Penner \cite{papapenner}  provided analogues for non--orientable surfaces of several results from Thurston theory of surfaces including the Dehn-Thurston coordinate function. Inspired by their work, the generalized Dynnikov coordinate system  was introduced in \cite{pamukyurttas} for multicurves in $N_{g,n}$ which provides an explicit bijection between $\mathfrak{L}_{g,n}$ and a certain subset of $(\mathbb{Z}^{2(n+g-2)}\times \mathbb{Z}^{ g})\setminus\{0\}$.  Here, we give a modified version of the generalized Dynnikov coordinate system together with the formulae in Theorem \ref{lem:inverse} (a corrected version of Theorem 2.14 in \cite{pamukyurttas}) for the inverse of the Dynnikov coordinate function.  Furthermore, with a slight modification, we also describe generalized Dynnikov coordinates for multicurves in $N_{g,0}$, which wasn't covered in \cite{pamukyurttas}. Let $n >1$. The generalized Dynnikov coordinates can be described as follows:

  \begin{figure}[h!]
\begin{center}
\labellist
\small\hair 2pt
  \pinlabel {\begin{turn}{-90}$\scriptstyle{\alpha_{1}}$\end{turn}} [ ] at  105 250
          \pinlabel {\begin{turn}{-90}$\scriptstyle{\alpha_{2}}$\end{turn}} [ ] at  105 100
          \pinlabel {\begin{turn}{-90}$\scriptstyle{\alpha_{2i-3}}$\end{turn}} [ ] at  170 250
                    \pinlabel {\begin{turn}{-90}$\scriptstyle{\alpha_{2i-4}}$\end{turn}} [ ] at  170 90
                     \pinlabel {\begin{turn}{-90}$\scriptstyle{\alpha_{2i-1}}$\end{turn}} [ ] at  253 250
                    \pinlabel {\begin{turn}{-90}$\scriptstyle{\alpha_{2i-2}}$\end{turn}} [ ] at  253 90
                    
                     \pinlabel {\begin{turn}{-90}$\scriptstyle{\alpha_{2i+1}}$\end{turn}} [ ] at  330 250
                    \pinlabel {\begin{turn}{-90}$\scriptstyle{\alpha_{2i+2}}$\end{turn}} [ ] at  330 90
                    
                     \pinlabel {\begin{turn}{-90}$\scriptstyle{\alpha_{2n-3}}$\end{turn}} [ ] at  380 250
                    \pinlabel {\begin{turn}{-90}$\scriptstyle{\alpha_{2n-2}}$\end{turn}} [ ] at  380 90

     \pinlabel {$\scriptstyle{\beta_{1}}$} [ ] at  68 155
     \pinlabel {\begin{turn}{-90}$\scriptstyle{\beta_{i}}$\end{turn}} [ ] at  215 100
     \pinlabel {\begin{turn}{-90}$\scriptstyle{\beta_{i+1}}$\end{turn}} [ ] at  295 90
     \pinlabel {\begin{turn}{-90}$\scriptstyle{\beta_{n}}$\end{turn}} [ ] at  425 100
     \pinlabel {\begin{turn}{-90}$\scriptstyle{\beta_{n+1}}$\end{turn}} [ ] at  502 90
     \pinlabel {\begin{turn}{-90}$\scriptstyle{\beta_{n+i-1}}$\end{turn}} [ ] at  545 85
          \pinlabel {\begin{turn}{-90}$\scriptstyle{\beta_{n+i}}$\end{turn}} [ ] at  625 85

     \pinlabel {\begin{turn}{-90}$\scriptstyle{\beta_{n+g-2}}$\end{turn}} [ ] at  672 85
     \pinlabel {\begin{turn}{-90}$\scriptstyle{\beta_{n+g-1}}$\end{turn}} [ ] at  744 150
  \pinlabel {\begin{turn}{-90}$\scriptstyle{\gamma_{1}}$\end{turn}} [ ] at  465 270
          \pinlabel {\begin{turn}{-90}$\scriptstyle{\gamma_{2}}$\end{turn}} [ ] at  473 95
          
            \pinlabel {\begin{turn}{-90}$\scriptstyle{\gamma_{2i-1}}$\end{turn}} [ ] at  585 250
          \pinlabel {\begin{turn}{-90}$\scriptstyle{\gamma_{2i}}$\end{turn}} [ ] at  588 90
 \pinlabel {\begin{turn}{-90}$\scriptstyle{\gamma_{2g-3}}$\end{turn}} [ ] at  716 220
          \pinlabel {\begin{turn}{-90}$\scriptstyle{\gamma_{2g-2}}$\end{turn}} [ ] at  710 93
\endlabellist  
 \includegraphics[scale=0.48]{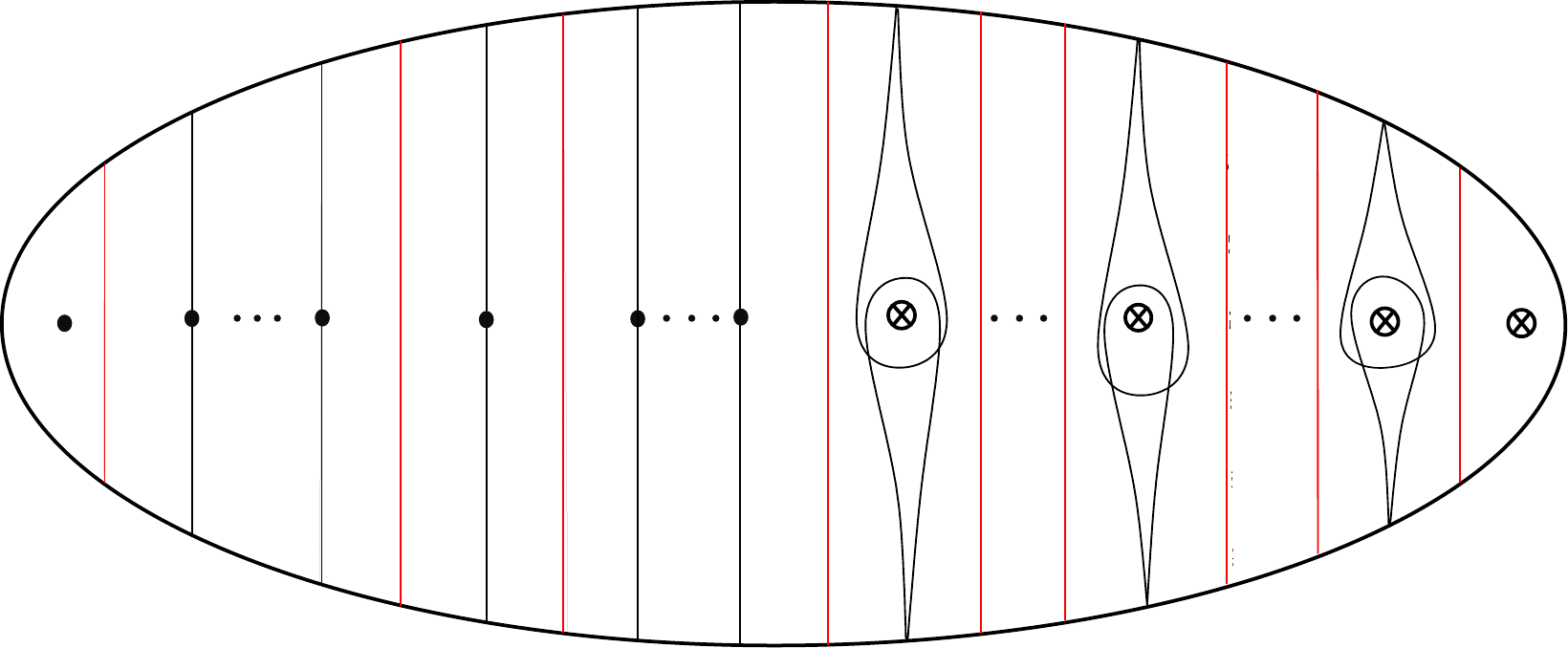}
\caption{The arcs $\alpha_i$, $\beta_i$, $\gamma_i$ and the core curves $\cC_i$ on $N_{g,n}$}\label{arcsproof}
\end{center}
\end{figure}

Let $\mathcal{A}_{g,n}$ be the set of arcs $\alpha_i$ ($1\leq i\leq 2n-2$), $\gamma_i$ ($1\leq i\leq 2g-2$) and $\beta_i$ ($1\leq i\leq n+g-1$) as depicted in Figure \ref{arcsproof}: the arcs $\alpha_{2i-3}$ and $\alpha_{2i-2}$  ($2\leq i\leq n$) join the $i$-th puncture to the boundary,  the {\it {teardrops}} $\gamma_{2i-1}$ and $\gamma_{2i}$ encircle the $i$--th crosscap and have  endpoints on the boundary, and the arc $\beta_i$ $(1\leq i\leq n-1)$ has endpoints on the boundary and passes between the $i$--th and $(i+1)$--th punctures, $\beta_n$  passes between the $n$--th puncture and the first crosscap, and  $\beta_{n+i}$  $(2\leq i\leq g-1)$ passes between the $i$--th and $(i+1)$--th crosscaps. Finally, $c_i$ ($1\leq i\leq g$) denotes the core curve of the $i$--th crosscap.

%lying between punctures and/or crosscaps;  and the core curves $\cC_i$ ($1\leq i\leq g$) of the crosscaps 

  Given $\cL\in \mathfrak{L}_{g,n}$  let $L$ be a minimal representative of $\cL$ (that is, $L$ intersects each of the arcs and curves minimally). For the sake of brevity, let $\alpha_i,~  \beta_i$,~ $\gamma_i$ also denote  the number of  intersections of $L$ with the  corresponding arcs.  We write $c_i=-1$ if $\cL$ contains the $i$--th core curve, $c_i=-2m$ if $\cL$ contains $m$ disjoint copies of the double cover of the $i$--th core curve and $c=-2m-1$ if $\cL$ contains $m$ disjoint copies of the double cover of the $i$--th core curve  plus the core curve itself.   Otherwise $c_i$ denotes  the number of  intersections of $L$ with the core curve of the $i$--th crosscap.  It will always be clear from the context  whether the symbols $\alpha_i,~  \beta_i$,~ $\gamma_i$  and  $c_i$  refer to arcs and curves rather than to integers. We write $(\alpha;\, \beta;\,\gamma,\, c)\in  \mathbb{Z}^{3n+3g-5}\setminus\{0\}$ for the collection of  these integers associated with $\cL$. Let $x^+=\max(x,0)$ throughout the text. 
   \begin{figure}[h!]
\begin{center}
  \labellist
\small\hair 2pt

 \pinlabel {$\scriptstyle{0}$} [ ] at  55 230
  \pinlabel {$\scriptstyle{2}$} [ ] at  85 240

 \pinlabel {$\scriptstyle{2}$} [ ] at  110 230
 \pinlabel {$\scriptstyle{4}$} [ ] at  200 230
  
  \pinlabel {$\scriptstyle{2}$} [ ] at  155 240
 \pinlabel {$\scriptstyle{6}$} [ ] at  155 30
  \pinlabel {$\scriptstyle{2}$} [ ] at  85 30

 \pinlabel {$\scriptstyle{1}$} [ ] at  80 155
 \pinlabel {$\scriptstyle{1}$} [ ] at  150 155
 \pinlabel {$\scriptstyle{1}$} [ ] at  230 155

 \endlabellist  
\includegraphics[scale=.5]{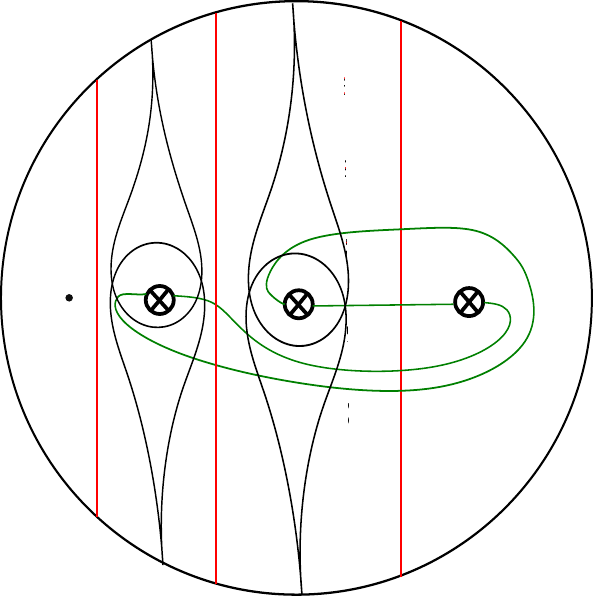}
\caption{A curve  $\cL$ in $N_{3,1}$ with generalized Dynnikov coordinates $\rho(\L)=(-1, -1; 2; 1,1,1)$}\label{fig:example1}
\end{center}
\end{figure}

 Let the function $\rho \colon \mathfrak{L}_{g,n}\to  (\mathbb{Z}^{2(n+g-2)}\times \mathbb{Z}^{ g})\setminus\{0\}$  be defined by 
$$\rho(\mathcal{L}) = (a;\, b;\,t;\, c):=(a_1,\,\ldots,a_{n-1};\,t_1,\ldots, t_{g-1};\, b_1,\, \ldots, b_{n+g-2};\, c_1,\ldots,\,c_g)$$

\noindent where 
\begin{align}\label{dynnikovformula1}  
\displaystyle a_i&=  
 \frac{\alpha_{2i}-\alpha_{2i-1}}{2}  \hspace{30pt} \textrm{;} \hspace{20pt} 1\leq i\leq n-1, \\
 \displaystyle t_{i}&= \frac{\gamma_{2i}-\gamma_{2i-1}}{2} \hspace{32pt} \textrm{;} \hspace{20pt} 1\leq i\leq g-1, \\
b_i&=\frac{\beta_i-\beta_{i+1}}{2} \hspace{42pt}  \textrm{;} \hspace{18pt} 1\leq i\leq n+g-2\end{align}

\noindent We say that $ (a;\, b;\,t;\, c)$  are the generalized Dynnikov coordinates of $\cL \in \mathfrak{L}_{g,n}$.
\begin{notation}Let $\psi_i=\max(c_i^+-|b_{n+i-1}|, 0)$ (the use of this parameter will be explained later) and $\mathcal{S}_{g,n}=\{ (a;\, b;\,t;\, c)\in (\mathbb{Z}^{2(n+g-2)}\times \mathbb{Z}^g)\setminus \left\{0\right\}: |t_i|+\psi_i \quad \text{is even for}\quad 1\leq i\leq g-1\}$.
 \end{notation}
 \begin{remark}
 Note the special case $n=1$ where there is no $a_i$ coordinate, and the special case $g=1$ where there is no $t_i$ coordinate. 
 \end{remark}
 
 The intersection numbers  $(\alpha;\, \beta;\,\gamma,\, c)$ (and hence the multicurve $\cL$) can be recovered from the generalized Dynnikov coordinates  $(a;\, b;\,t;\, c)\in \mathcal{S}_{g,n}$.  Theorem \ref{lem:inverse} gives the inverse of the generalized Dynnikov coordinate function by presenting a formula that describes multicurves from given generalized Dynnikov coordinates.

	\begin{thm}\label{lem:inverse}
Let $ (a;\, b;\,t;\, c)\in \mathcal{S}_{g,n}$. Then $ (a;\, b;\,t;\, c)$ corresponds to a unique multicurve in $\cL \in \mathfrak{L}_{g, n}$ which has 

\begin{align}
\displaystyle \alpha_i&= \left\{     
\begin{array}{lr}
(-1)^i a_{\lceil i/2 \rceil}+\frac{\beta_{\lceil i/2\rceil}}{2}& \mbox{if $b_{\lceil i/2\rceil} \geq 0$,}\label{eq:finalalpha} \\
(-1)^i a_{\lceil i/2\rceil}+\frac{\beta_{1+\lceil i/2\rceil}}{2}& \mbox{if $b_{\lceil i/2\rceil} \leq 0$,}   
\end{array}
 \right.\\
 \displaystyle \gamma_i     &= \left\{     
\begin{array}{lr}
(-1)^i t_{\lceil i/2 \rceil}+\beta_{n+\lceil i/2\rceil-1}+\psi_{\lceil i/2 \rceil}& \mbox{if $b_{n+\lceil i/2\rceil}-1 \geq 0$,} \label{eq:finalgamma}\\
(-1)^i t_{\lceil i/2\rceil}+\beta_{n+\lceil i/2\rceil}+\psi_{\lceil i/2 \rceil}& \mbox{if $b_{n+\lceil i/2\rceil}-1 \leq 0$,}   
\end{array}
 \right.\\
   \beta_i&= Z_i+2\max(0,c_g-\frac{Z_{n+g-1}}{2})\label{eq:finalbeta}
 \end{align}

\noindent  where
\begin{align*}
X_i&=2\left[\left|a_i \right|+\max(b_i,0)+\displaystyle\sum^{i-1}_{k=1} b_k\right]\\
Y_i&=2\left[\max(b_{n+i-1}, 0)+\displaystyle\sum^{n+i-2}_{k=1}b_k\right]+\left|t_i \right|+\psi_i\\
Z_i&=\max_{\substack{1\leq s \leq n-1\\
                  1\leq t\leq g-1}}
        \left\{X_s,Y_t\right\}-2\sum^{i-1}_{k=1}b_k
\end{align*}
 \noindent   Here $\lceil x \rceil$~denotes the smallest integer which is not less than $x$.  
	\end{thm}

	 The main goals of this paper are first to prove Theorem \ref{lem:inverse} (correcting the proof of Theorem 2.14 in \cite{pamukyurttas}); and then give the derivation of the formulae in Theorem \ref{thm:update} which describes  how generalized Dynnikov coordinates change under the action of crosscap transpositions $u_i$ and $u^{-1}_i$ ($1\leq i\leq g-1$). Therefore, Theorem \ref{thm:update} computes for each mapping class $\beta$ written as a word of crosscap transpositions, $\beta:\mathcal{S}_{g,n}\to\mathcal{S}_{g,n}$ given by,
$\beta(a;\,t;\,b;\,c)=\rho\circ\beta\circ\rho^{-1}(a;\,t;\,b;\,c).$ A crosscap transposition $u_i$ is a generator of the mapping class group of $N_{g,n}$ (i.e. group of isotopy classes of homeomorphisms of $N_{g,n}$) exchanging crosscaps~$\otimes_i$ and~$\otimes_{i+1}$ in the {\it counterclockwise} direction  keeping each of the remaining crosscaps fixed \cite{korkmaz, parlak}. We note that while the formulae in Theorem \ref{lem:inverse} and Theorem \ref{thm:update} seem to have a complicated form,  the method we use to obtain them is transparent since it purely relies on algebraic calculations and the properties of multicurves in terms of their associated  intersection numbers $(\alpha;\, \beta;\,\gamma,\, c)$. In addition,  the formulae are ideally suited for computer implementation.

	\begin{notation}\label{not:semiring}
 For computational and notational
convenience, we will work in the {\it max-plus semiring} $(\mathbb{R},\max,+)$ equipped with the additional and multiplicative operations $a+ b = \max(a,b)$ and $a\times b = a+b$  to obtain the formulae in Theorem \ref{thm:update}.  As we  use the normal notation of addition, multiplication, and division  we enclose the formulae in square brackets to indicate that these will be interpreted in the max--plus sense. Therefore, $ \left[a+b\right] = \max(a,b),\, \left[ab\right]= a+b,\, \left[a/b\right] = a-b,\,\left[1\right]=0$. 

\end{notation}

To prove Theorem \ref{thm:update} we shall make use of particular arc systems  called {\it clovers} and {{\it scales}, each of which is associated with an exceptional parameter, certain linear combinations of generalized Dynnikov coordinates, denoted $d_i, e_i, \bar{e}_i, f_i,\bar{f}_i, g_i$ and $\bar{g}_i$.

\begin{notation}
For notational convenience we write $B_j=2b_j$ (i.e. $[B_j]=[b^2_j]$) in Theorem \ref{thm:update}.
\end{notation}

\begin{thm}\label{thm:update}

Let $\cL\in \mathfrak{L}_{g,n}$ have generalized Dynnikov coordinates $(a;\, b;\,t;\, c)$.  Let $(a';\, b';\,t';\, c')$ and  $(a'';\, b'';\,t'';\, c'')$ be the generalized Dynnikov coordinates of $u_i(\cL)$ and $u^{-1}_i(\cL)$ respectively. Then $a'_j=a''_j=a_j$, $b'_j=b''_j=b_j$ for all $1\leq j\leq n$;  $(c'_i, c'_{i+1})=(c_{i+1},c_i)$ and $(c''_i, c''_{i+1})=(c_{i+1},c_i)$ for $1\leq i\leq g-1$; and for $1\leq i\leq g-3$ we have
\leavevmode
%b''_{n+i}&=\left[\frac{\bar{e_i}}{d_i}b^2_{n+i}(\frac{t''_{i}+t'_{i+1}}{t''_{i}})\frac{t_{i+1}}{t_{i}+t_{i+1}}\right]\quad\text{and}\quad b''_{n+i-1}=\left[\frac{b_{n+i}b_{n+i-1}}{b''_{n+i}}\right]

\begin{align}\label{eq:update}
  \begin{aligned}
  t_{i}' &= \tropical{ g_i\big(t_{i}(1+d_iB_{n+i-1})+d_iB_{n+i-1}t_{i+1}\big) }\\     t_{i+1}' &= \tropical{ \frac{t_{i}t_{i+1}B_{n+i}}{f_i\big(t_i(d_{i}+B_{n+i})+ d_it_{i+1})} }\\t_{i}'' &= \left[\frac{t_it_{i+1}}{\bar{g_i}\big(t_{i}d_iB_{n+i-1}+t_{i+1}(d_iB_{n+i-1}+1)\big)}\right]\\ t_{i+1}'' &= \left[\frac{\bar{f_i}\big(t_{i+1}(B_{n+i}+d_i)+d_it_i\big)}{B_{n+i}}\right]
  \end{aligned}
  &&
  \begin{aligned}
   B_{n+i}' &= \tropical{ \frac{e_i}{d_i}B_{n+i}(\frac{t'_{i}+t'_{i+1}}{t'_{i+1}})\frac{t_{i}}{t_{i}+t_{i+1}} } \\       B_{n+i-1}' &= \tropical{ \frac{B_{n+i}B_{n+i-1}}{B'_{n+i}}  }\\B''_{n+i-1}&=\left[\frac{B_{n+i}B_{n+i-1}}{B''_{n+i}}\right],\\ B''_{n+i}&=\left[\frac{\bar{e_i}}{d_i}B_{n+i}(\frac{t''_{i}+t''_{i+1}}{t''_{i}})\frac{t_{i+1}}{t_{i}+t_{i+1}}\right] \end{aligned}
    \end{align}

In the special case the formulae above is interpreted as 

\begin{align}\label{eq:updatespecial}
  \begin{aligned}
t_{g-1}' &= \tropical{g_{g-1}\big(t_{g-1}+d_{g-1}B_{n+g-2}(1+t_{g-1})\big)}\\
t_{g-1}'' &= \tropical{\frac{t_{g-1}}{\bar{g}_{g-1}\big(d_{g-1}B_{n+g-2}(1+t_{g-1})\big)} } \end{aligned} 
 &&
  \begin{aligned}
  B_{n+g-2}' &= \tropical{ \frac{d_{g-1}}{e_{g-1}}B_{n+g-2}\frac{1+t_{g-1}}{t_{g-1}(1+t'_{g-1})} }\\
B_{n+g-2}'' &= \tropical{ \frac{d_{g-2}}{\bar{e}_{g-2}}B_{n+g-2}\frac{t''_{g-1}(1+t_{g-1})}{1+t''_{g-1}}}.
\end{aligned}
\end{align}
In all other cases $t'_j=t_j$, $t''_j=t_j$ and $b'_j=b''_j=b_j$.  

\end{thm}

The paper is organized as follows.  Section \ref{section2} provides background material and contains a detailed study of  generalized Dynnikov coordinates giving proofs of Theorem \ref{lem:inverse} and Theorem \ref{lem:inverse2}. In Section \ref{section3},  we introduce the notions of clovers and scales, certain collections of adjacent arcs in $\mathcal{A}_{g,n}$ and their images under $u_i$ and $u^{-1}_i$  from which we obtain clover and scale  equalities  given in Lemma \ref{lem:clover1}, Lemma \ref{lem:clover2}, Lemma \ref{lem:clover3}  and Lemma \ref{lem:scales1}  which play key roles in the derivation of the formulae in Theorem \ref{thm:update} which we also prove  in Section \ref{section3}.

\section{Constructing multicurves from generalized Dynnikov coordinates}\label{section2}
 In this section we prove Theorem \ref{lem:inverse}  recalling basic properties that  a minimal representative  $L$  satisfies in terms of  the intersection numbers $(\alpha;\, \beta;\,\gamma,\, c)$.  Let $1\leq i\leq n-1$ and $S_i$ denote the region bounded by  the arcs $\beta_{i}$ and $\beta_{i+1}$ containing puncture $i+1$. We denote  by $S_0$ the left most region bounded by $\beta_1$ and the boundary  containing the first puncture.  Now, let $1\leq i \leq g$. Then $S'_i$ denotes the region bounded by  the arcs $\beta_{n+i-1}$ and $\beta_{n+i}$ containing the $i$--th crosscap.  $S'_g$ denotes the right most region bounded by $\beta_{n+g-1}$ and the boundary containing the last crosscap.  Since $L$ is minimal, there are finitely many connected components of $L\cap S_i$ and $L\cap S'_i$ which are depicted in Figure \ref{calculation1a} and Figure \ref{calculation1b}.

 \begin{figure}[h!]
\begin{center}
  \labellist
\small\hair 2pt
 \pinlabel {$\scriptstyle{\alpha_{2i-1}}$} [ ] at  35 190
  \pinlabel {$\scriptstyle{\alpha_{2i}}$} [ ] at  38 -5
   \pinlabel {$\scriptstyle{\alpha_{2i-1}}$} [ ] at  153 190
  \pinlabel {$\scriptstyle{\alpha_{2i}}$} [ ] at  155 -5

   \pinlabel {$\scriptstyle{\beta_{i}}$} [ ] at  2 190
      \pinlabel {$\scriptstyle{\beta_{i+1}}$} [ ] at  70 190
       \pinlabel {$\scriptstyle{\beta_{i}}$} [ ] at  115 190
      \pinlabel {$\scriptstyle{\beta_{i+1}}$} [ ] at  185 190

  \pinlabel {${(a)}$} [ ] at  40 -20
    \pinlabel {${(b)}$} [ ] at  155 -20

 \endlabellist  
\includegraphics[scale=.625]{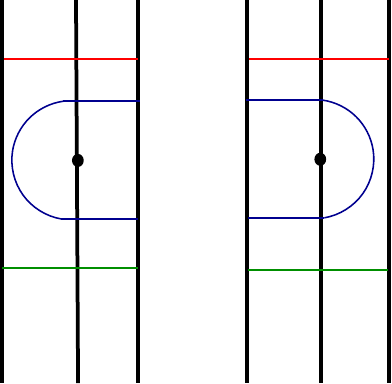}
\vspace{0.1 cm}
\caption{Connected components of $L\cap S_i$}\label{calculation1a}
\end{center}
\end{figure}

 \begin{figure}[h!]
\begin{center}
  \labellist
\small\hair 2pt

     \pinlabel {$\scriptstyle{\beta_{n+i-1}}$} [ ] at  120 195

          \pinlabel {$\scriptstyle{\beta_{n+i}}$} [ ] at  200 195

     \pinlabel {$\scriptstyle{\beta_{n+i-1}}$} [ ] at  0 195
          \pinlabel {$\scriptstyle{\beta_{n+i}}$} [ ] at  76 195

  \pinlabel {$\scriptstyle{\beta_{n+i-1}}$} [ ] at  242 195
    \pinlabel {$\scriptstyle{\beta_{n+i}}$} [ ] at  322 195

     \pinlabel {$\scriptstyle{\beta_{n+i-1}}$} [ ] at 364 195
          \pinlabel {$\scriptstyle{\beta_{n+i}}$} [ ] at 444 195
          
  \pinlabel {$\scriptstyle{\gamma_{2i-1}}$} [ ] at  160 195
    \pinlabel {$\scriptstyle{\gamma_{2i-1}}$} [ ] at  38 195

  \pinlabel {$\scriptstyle{\gamma_{2i-1}}$} [ ] at  282 195

  \pinlabel {$\scriptstyle{\gamma_{2i-1}}$} [ ] at  400 195

   \pinlabel {$\scriptstyle{\gamma_{2i}}$} [ ] at  160 -5
    \pinlabel {$\scriptstyle{\gamma_{2i}}$} [ ] at  35 -5

  \pinlabel {$\scriptstyle{\gamma_{2i}}$} [ ] at  283 -5

  \pinlabel {$\scriptstyle{\gamma_{2i}}$} [ ] at  405 -5
  \pinlabel {${(a)}$} [ ] at  35 -20
    \pinlabel {${(b)}$} [ ] at  160 -20
     \pinlabel {${(c)}$} [ ] at  282 -20
    \pinlabel {${(d)}$} [ ] at  405 -20

 \endlabellist  
\includegraphics[scale=.625]{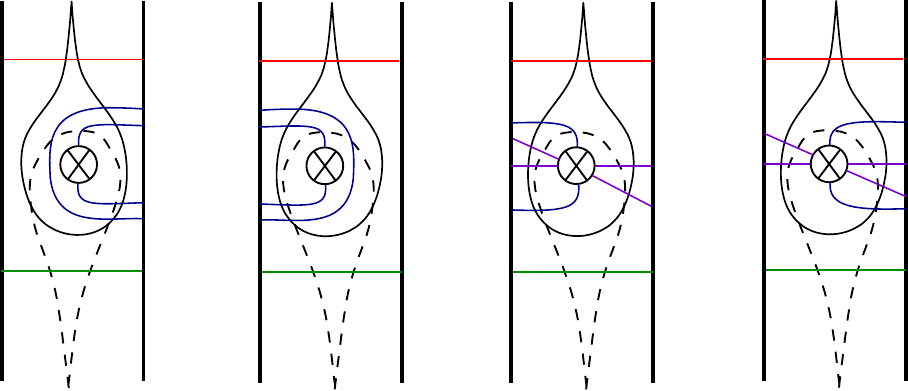}
\vspace{0.1 cm}
\caption{Connected components of $L\cap S'_i$}\label{calculation1b}
\end{center}
\end{figure}

\begin{itemize}	
	\item An {\it{above component}} of  $L\cap S_i$ ($1\leq i \leq n-1$)  has one endpoint on each of the arcs $\beta_{i}$ and $\beta_{i+1}$, and intersects $\alpha_{2i-1}$ but not $\alpha_{2i}$. An {\it{above component}}  of  $L\cap S'_i$ has one endpoint on each of the arcs $\beta_{n}$ and $\beta_{n+1}$, and  intersects $\gamma_{2i-1}$ but not $\gamma_{2i}$.
	
	\vspace{0.1 cm}
\item An {\it{below component}} of  $L\cap S_i$ ($1\leq i \leq n-1$)  has one endpoint on each of the arcs $\beta_{i}$ and $\beta_{i+1}$, and intersects $\alpha_{2i-1}$ but not $\alpha_{2i}$. An {\it{below component}}  of  $L\cap S'_i$  has one endpoint on each of the arcs $\beta_{n+i-1}$ and $\beta_{n+i}$, and  intersects $\gamma_{2i}$ but not $\gamma_{2i-1}$.
\vspace{0.1 cm}

 \item A {\it{left loop component}} of $L\cap S_i$ has both endpoints on $\beta_{i+1}$, and intersects both of  the arcs $\alpha_{2i-1}$ and $\alpha_{2i}$. A  {\it{left loop component}} of $L\cap S'_i$  has  both endpoints on $\beta_{n+i}$, and intersects both of  the arcs $\gamma_{2i-1}$ and $\gamma_{2i}$. If it intersects the core curve, it is called a left core loop. Otherwise it is called a left non-core loop.
 
 \vspace{0.1 cm}

     \item  A {\it{right loop component}} of $L\cap S_i$ has both endpoints on $\beta_{i}$, and intersects both of the arcs $\alpha_{2i-1}$ and $\alpha_{2i}$.   A {\it{right loop component}} of $L\cap S'_i$  has both endpoints on $\beta_{n+i-1}$, and intersects both of  the arcs $\gamma_{2i-1}$ and $\gamma_{2i}$.  If it intersects the crosscap, it is called a right core loop of  $L\cap S'_i$. Otherwise it is called a right non-core loop of $L\cap S'_i$. 
     
      \vspace{0.1 cm}
      
    \item  A {\it{straight component}} of $L\cap S'_i$  has one endpoint on each of the arcs $\beta_{n+i-1}$ and $\beta_{n+i}$, and intersects the core curve and both of  the arcs $\gamma_{2i-1}$ and $\gamma_{2i}$.
      
      \end{itemize}
      
    Above, below and loop components of  $L\cap S_i$ and $L\cap S'_i$  are depicted red, green and blue respectively in Figure  \ref{calculation1a} and Figure  \ref{calculation1b}.   Straight components of $L\cap S'_i$ are depicted purple in Figure \ref{calculation1b}(c) and Figure \ref{calculation1b}(d). Observe that there can only be left loop components in $S_0$ and right loop components in $S'_g$.  The following lemma gives two important equalities which are obvious from Figure  \ref{calculation1a} and Figure  \ref{calculation1b}. 
      
      \begin{lem}\label{lem:equalities}
Let $L$ be a minimal representative of a  multicurve $\mathcal{L}\in \mathfrak{L}_{g,n}$  with intersection numbers $(\alpha;\, \beta;\,\gamma,\, c)$.  Let $\psi_i$ denote the number of straight components of  $L\cap S'_i$. Then, 
\begin{align}\label{eq:relations1a}
\max(\beta_{i}, \beta_{i+1})&=\alpha_{2i-1}+\alpha_{2i}
\end{align}

\begin{align}\label{eq:relation1b}
\max(\beta_{n+i-1}, \beta_{n+i})&=\frac{\gamma_{2i-1}+\gamma_{2i}}{2}-\psi_i
\end{align}
\end{lem}

    \begin{remark}\label{rem:properties1}
    
    	Given a minimal representative $L$ of $\cL\in \mathfrak{L}_{g,n}$ we can initially observe that every component of $L$ intersects each $\beta_i$ and hence each $\alpha_{2i-1}\cup \alpha_{2i}$ and  $\gamma_{2i-1}\cup \gamma_{2i}$ an even number of times.  Therefore $a_i, t_i$ and $b_i$ are integers.\end{remark}

	\begin{lem}\label{lem:loops}
	For each  $1\leq i \leq n+g-2$ let $b_i=\frac{\beta_i-\beta_{i+1}}{2}$. Then there are $|b_i|$  loop components in $S_i$ ( $1\leq i \leq n-1$) and $|b_{n+i-1}|$ loop components in  $S'_i$ ($1\leq i \leq g-1)$. If $b_i>0$ the loop components are right and if $b_i<0$ the loop components are left.

	\end{lem}	

\begin{proof}
We prove the statement for $S_i$ (the argument for $S'_i$ is identical). Let $1\leq i \leq n-1$.  We first note  that there cannot be both left loop and right loop components in $S_i$  since the curves are mutually disjoint. Assume without loss of generality that $\beta_{i+1}\geq \beta_i$. Observe from Figure \ref{calculation1a}($b$)  that the additional intersections on $\beta_{i+1}$ come from left loop components in $S_i$ since above and below components intersect both $\beta_{i}$ and $\beta_{i+1}$ the same number of times. Since each left loop intersects $\beta_{i+1}$ twice it follows that there are  $-b_i=\frac{\beta_{i+1}-\beta_i}{2}$ left loop components in $S_i$. \end{proof}

\begin{rem}\label{rem:endregions}

%Since there are only left loop components of $L\cap S_0$ and right loop components of $L\cap S'_g$, 
The number of loop components of $L\cap S_0$ is given by $\frac{\beta_{1}}{2}$, and the number of right loop components of $L\cap S'_g$ is given by $\frac{\beta_{n+g-1}}{2}$.

\end{rem}
\begin{lem}\label{lem:loops2}
 Let $1\leq i< n+g-2$, and $\lambda_{i}$,  $\lambda_{c_i}$ and $\psi_{i}$ denote the number of non-core loop, core loop and straight components of $L\cap S'_i$. Then,

\begin{align}\psi_{i}&=\max(c_i-|b_{n+i-1}|,0)\\
 \lambda_{i}&=\max(|b_{n+i-1}|-c_i, 0) \quad\text{,} \quad   \lambda_{c_i}= |b_{n+i-1}|-\lambda_i
 \end{align}
\end{lem}

\begin{proof}
%If $c_i\leq 0$, there can only be components of $L\cap S'_i$ which are isotopic to M\"{o}bius curves or non--core loop components of $L\cap S'_i$.  Let $c_i>0$.
There are three possibilities for a connected component of $L\cap S'_i$ which intersects the crosscap. It can be a left core loop or a right core--loop  or  a straight core component.  Observe from Figure \ref{calculation1b} that we have 
\begin{align}\label{eq:straight1}
c_i&=\psi_i+ \lambda_{c_i}
\end{align}
\begin{align}\label{eq:straight2}
|b_{n+i-1}|&=\lambda_{i}+\lambda_{c_i}. 
\end{align}

 If $c_i-|b_{n+i-1}|\geq 0$,  there exist components of $L\cap S'_i$ other than core loop  components which intersect the crosscap. Such components can only be straight components and hence $\psi_i>0$ and $\lambda_i=0$ since non--core loops and straight components cannot exist at the same time. Then, $|b_{n+i-1}|=\lambda_{c_i}$  and hence  $\psi_i=c_i-|b_{n+i-1}|$  by  (\ref{eq:straight1}) and (\ref{eq:straight2}).  If $c_i-|b_{n+i-1}|<0$, then there exist  non--core loop components as well as the core loop components. That is,  $\lambda_i>0$ and hence $\psi_i=0$. Therefore,  $c_i=\lambda_{c_i}$ and hence $\lambda_{i}=|b_{n+i-1}|-c_i$ by  (\ref{eq:straight1}) and (\ref{eq:straight2}).  Therefore, we get 
 $\psi_{i}=\max(c_i-|b_{n+i-1}|,0)$ and $\lambda_i=\max(|b_{n+i-1}|-c_i, 0)$ as required. We immediately get from  (\ref{eq:straight2}) that     $\lambda_{c_i}= |b_{n+i-1}|-\lambda_i$. \end{proof}
 
\begin{remark}\label{rem:endrelation}
  Let $1\leq i< n+g-2$. There are $b^+_i=\max(b_i,\,0)$ right loops  and $(-b_i)^+=\max(-b_i,\,0)$ left loops  about crosscap $i$.  By Remark \ref{rem:endregions}  there are only right loop components of  $L\cap S'_g$ and the number of those is given by $\frac{\beta_{n+1}}{2}$. It immediately follows that there are $c_g$ core loops and hence $\lambda_g=\frac{\beta_{n+g-1}}{2}-c_g$  non--core loops of $L\cap S'_g$.
 \end{remark}

The following Lemma is obvious since each above and below component in $S_i$ intersects $\alpha_{2i-1}$ and $\alpha_{2i}$, and  each above and below component in $S'_i$ intersects $\gamma_{2i-1}$ and $\gamma_{2i}$ respectively (see Figure \ref{calculation1a} and Figure \ref{calculation1b}).

\begin{lem}\label{lem:abovebelow}
Let there be $A_i$ and $B_i$ $(1\leq i \leq n-1)$ above  and  below components of $L\cap S_i$; and $A'_i$ and $B'_i$ $(1\leq i \leq g-1)$ above and  below components of $L\cap S'_i$  respectively. Then,  
\begin{align}
A_i&=\alpha_{2i-1}-|b_i|\\
B_i&=\alpha_{2i}-|b_i|\\
A'_i&=\frac{\gamma_{2i-1}}{2}-\psi_i-|b_{n+i-1}|\\
B'_i&=\frac{\gamma_{2i}}{2}-\psi_i-|b_{n+i-1}|\end{align}

\end{lem}

 \noindent The curve in Figure \ref{fig:example1} has $A'_1=A'_2= B'_1=0; B'_2=2; \lambda_1=\lambda_2=0, \lambda_3=1; \lambda_{c_1}=\lambda_{c_2}=1=\lambda_{c_3}=1$; 
 $\psi_1=\psi_2=0$. These  parameters will frequently be referred to  throughout the paper. 

The generalized Dynnikov coordinate function $\rho:\mathfrak{L}_{g,n}\to \mathcal{S}_{g,n}$ is a bijection:   To describe its inverse, it is sufficient to describe a function from $\mathcal{S}_{g,n}$ to $\mathbb{Z}^{3n+3g-5}\setminus\{0\}$. It is easy to check that this function sends each $(a;\, b;\,t;\, c)\in \mathcal{S}_{g,n}$ to the intersection numbers $(\alpha;\, \beta;\,\gamma,\, c)$ associated with a multicurve $\cL$ with $\rho(\cL)=(a;\, b;\,t;\, c)$.

Next we prove Theorem \ref{lem:inverse}.  As we shall see the proof is completely constructive in the sense that it gives an explicit way of constructing a multicurve in $N_{g,n}$ in finite number of steps.

       \begin{proof}[Proof of Theorem \ref{lem:inverse}]

Let $L$ be a minimal representative of $\cL\in \mathfrak{L}_{g,n}$ with generalized Dynnikov coordinates $\rho(\cL)=(a;\, b;\,t;\, c)$.  Note that $2|a_i|=|\alpha_{2i-1}- \alpha_{2i}|$ and  $2|t_i|=|\gamma_{2i-1}- \gamma_{2i}|$ give  the difference between below and above components in $S_i$ and $S'_i$ respectively by Lemma \ref{lem:abovebelow}. Also $|b_i|$ gives the number of loop components in $S_i$ and $S'_i$  by Lemma \ref{lem:loops}. Let  $m_i$ and $n_i$ be the smaller of above and below components of $L\cap S_i$ and $L\cap S'_i$ respectively. From  Figure \ref{calculation1a} and Figure \ref{calculation1b} it is straightforward to compute $\beta_i$ and $\beta_{n+i-1}$:

\noindent  For $1\leq i\leq n-1$,

\begin{align*}
\beta_i &= \left\{ \begin{array}{ll}
         2m_i+2\left|a_i\right| & \mbox{if $b_i\leq 0$};\\
         2m_i+2\left|a_i\right|+ 2b_i & \mbox{if $b_i \geq 0$}.\end{array} \right.         \end{align*}

\noindent  For $1\leq i\leq g-1$

     \begin{align*}    
         \beta_{n+i-1} &= \left\{ \begin{array}{ll}
         2n_i+\left|t_i\right| +\psi_i& \mbox{if $b_{n+i-1}\leq 0$};\\
         2n_i+\left|t_i\right|+ \psi_i +2b_{n+i-1}& \mbox{if $b_{n+i-1} \geq 0$}.\end{array} \right.
         \end{align*}

from which we get 
\begin{align}
\beta_i&=2\left[\left|a_i \right|+\max(b_i,0)+m_i\right]\label{eq:relations1}\\
\beta_{n+i-1}&=|t_i|+\psi_i +2\left[\max(b_{n+i-1},\,0)+n_i\right]\label{eq:relations2}
\end{align}  

\noindent Since $\beta_{n+i-1}$ is even from Remark \ref{rem:properties1}, equality (\ref{eq:relations2}) implies that $|t_i|+\psi_i$ should be even. That is, $|t_i|+\max(c_i-|b_{n+i-1}|,0) $ is even by Lemma \ref{lem:loops2}. 

%If $\psi=0$ then$|t_i|$ is even Let $\psi_i\neq 0$. From Lemma \ref{eq:straight1}, $\psi_i=c_i^+-|b_{n+i-1}|$.
%.   If $\psi_i$ is odd, then $t_i+\psi_i$ is even. That is, $t_i+c_i^+-|b_{n+i-1}|$ is even. Similarly, if $\psi_i$ is even, then $|t_i|+\psi$ is even. That is, $|t_i|+c_i^+-|b_{n+i-1}|$ is even. Similarly, 

Now, consider a subarc of $L$ which intersects the last crosscap exactly once, has zero intersection with the other  crosscaps and intersects the horizontal diameter of the surface only between the  first puncture and the boundary exactly once as shown in Figure \ref{fig:figR}. Each such arc intersects each $\beta_i$ and $\gamma_i$ twice, and each $\alpha_i$ exactly once. We say that such arcs are {\it{almost boundary parallel}}, and write $R$ for the number of almost boundary parallel arcs.

\begin{figure}[h!]
\begin{center}
\small\hair 2pt
\includegraphics[scale=0.6]{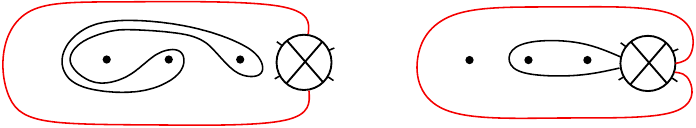}
\caption{Two multicurves  oin $N_{1,3}$ with $R=1$}\label{fig:figR}
\end{center}
\end{figure}

 Using  $\beta_i = \beta_1 - 2\displaystyle\sum^{i-1}_{j=1} b_j$ ($1\leq i \leq n+g-1$) and  subtracting $2R$ from both sides of equation \ref{eq:relations1} and equation \ref{eq:relations2} we get

\begin{align}
\beta_1-2R&=2\left[\left|a_i \right|+\max(b_i,0)+m_i-R+\displaystyle\sum^{i-1}_{j=1}b_j \right]\quad  \text{for}\quad 1\leq i \leq n-1\label{eq:relations3}\\
\beta_1-2R&=2\left[\max(b_{n+i-1}, 0)+n_i-R+\displaystyle\sum^{n+i-2}_{j=1}b_j\right]+\left|t_i \right|+\psi_i \quad  \text{for}\quad 1\leq i \leq g-1.
\end{align}

One crucial fact is that $m_i-R=0$  for some $1\leq i\leq n-1$ or $n_i-R=0$ for some   $1\leq i\leq g-1 $ since otherwise there would be both above and below components in each of the $S_i$ and $S'_i$ except for those which arise from almost boundary arcs, but this would mean $L$ contains  boundary parallel curves which is impossible. Then,
 
When $m_i-R=0$;

\begin{eqnarray*}
\beta_1-2R=2\left[\left|a_i \right|+\max(b_i,0)+\displaystyle\sum^{i-1}_{j=1} b_j \right],
\end{eqnarray*}

When $m_i-R>0$;

\begin{eqnarray*}
\beta_1-2R>2\left[\left|a_i \right|+\max(b_i,0)+\displaystyle\sum^{i-1}_{j=1} b_j\right],
\end{eqnarray*}

When $n_i-R=0$;

\begin{eqnarray*}
\beta_1-2R=2\left[\max(b_{n+i-1}, 0)+\displaystyle\sum^{n+i-2}_{j=1}b_j\right]+\left|t_i \right|+\psi_i,
\end{eqnarray*}

When $n_i-R>0$;

\begin{eqnarray*}
\beta_1-2R>2\left[\max(b_{n+i-1}, 0)+\displaystyle\sum^{n+i-2}_{j=1}b_j\right]+\left|t_i \right|+\psi_i.
\end{eqnarray*}

Therefore, setting
\begin{align*}
X_i&=2\left[\left|a_i \right|+\max(b_i,0)+\displaystyle\sum^{i-1}_{j=1} b_j\right]\\
Y_i&=2\left[\max(b_{n+i-1}, 0)+\displaystyle\sum^{n+i-2}_{j=1}b_j\right]+\left|t_i \right|+\psi_i 
\end{align*}
we get 

\begin{equation}
  \beta_1-2R= \max_{\substack{1\leq s \leq n-1\\
                  1\leq k\leq g-1}}
        \left\{X_s,Y_k\right\}
\end{equation}
and hence 

\begin{equation}\label{eq:finalbeta1}
  \beta_i-2R= \max_{\substack{1\leq s \leq n-1\\
                  1\leq t\leq g-1}}
        \left\{X_s,Y_k\right\}-2\sum^{i-1}_{j=1}b_j
\end{equation}
as required. Next, we compute $R$. Let

 $$Z_i=\max_{\substack{1\leq s \leq n-1\\
                  1\leq k\leq g-1}}
        \left\{X_s,Y_k\right\}-2\sum^{i-1}_{j=1}b_j$$
 By (\ref{eq:finalbeta1}), $\beta_{n+g-1}=Z_{n+g-1}+2R$.  Also by Remark \ref{rem:endrelation} $\beta_{n+g-1}=2c_g+2\lambda_g$. It follows that when $Z_{n+g-1}>2c_g$ we have 
$\beta_{n+g-1}>2c_g$ and hence $\lambda_g\neq 0$. Therefore, $R=0$ since almost  boundary parallel arcs and non-core loop components of $L\cap S'_g$ cannot exist  at the same time, and when $Z_{n+g-1}< 2c_g$ we have $Z_{n+g-1}< \beta_{n+g-1}$ and hence $R>0$ and $\lambda_g=0$ which implies $\beta_{n+g-1}=Z_{n+g-1}+2R=2c_g$ that is  $R=c_g-\frac{Z_{n+g-1}}{2}$. Therefore, $R=\max(0,c_g-\frac{Z_{n+g-1}}{2})$.

To compute $\alpha_i$ and $\gamma_i$ we make use of the equalities in Lemma \ref{lem:equalities}:
\begin{align}\label{eq:relations1a}
\max(\beta_{i}, \beta_{i+1})&=\alpha_{2i-1}+\alpha_{2i}
\end{align}

\begin{align}\label{eq:relation1bx}
\max(\beta_{n+i-1}, \beta_{n+i})&=\frac{\gamma_{2i-1}+\gamma_{2i}}{2}-\psi_i
\end{align}

\noindent Since  $2a_i= \alpha_{2i}-\alpha_{2i-1}$ ($1\leq i\leq n-1$) we get from \ref{eq:relations1a} that

If $\beta_i\geq \beta_{i+1}$ (i.e. $b_i\geq 0$)

\begin{align*}
\alpha_{2i}= a_i+\frac{\beta_i}{2}; \quad \text{and}\quad \alpha_{2i-1}=-a_i+\frac{\beta_i}{2}.
\end{align*}

If $\beta_{i+1}\geq \beta_{i}$ (i.e. $b_i\leq 0$)

\begin{align*}
\alpha_{2i}= a_i+\frac{\beta_{i+1}}{2};  \quad \text{and}\quad \alpha_{2i-1}=-a_i+\frac{\beta_{i+1}}{2}.
\end{align*}

That is to say:
\begin{eqnarray*}
\alpha_i&= \left\{ \begin{array}{ll}
         (-1)^ia _{\lceil i/2 \rceil}+\frac{\beta_{\lceil i/2\rceil}}{2}& \mbox{if $b_{\lceil i/2\rceil} \geq 0$};\\
        (-1)^i a_{\lceil i/2\rceil}+\frac{\beta_{1+\lceil i/2\rceil}}{2}& \mbox{if $b_{\lceil i/2\rceil}\leq 0$}   
         \end{array} \right. 
\end{eqnarray*} 

Similarly, since  $2t_i= \gamma_{2i}-\gamma_{2i-1}$  for each $1\leq i\leq g-1$ we get from \ref{eq:relation1bx} that

If $\beta_{n+i-1}\geq \beta_{n+i}$ (i.e. $b_{n+i-1}\geq 0$)

\begin{align*}
\gamma_{2i}= t_i+\beta_{n+i-1}+\psi_i; \quad \text{and}\quad \gamma_{2i-1}=-t_i+\beta_{n+i-1}+\psi_i,
\end{align*}

If $\beta_{n+i-1}\leq \beta_{n+i}$ (i.e. $b_{n+i-1}\geq 0$)

\begin{align*}
\gamma_{2i}= t_i+\beta_{n+i}+\psi_i; \quad \text{and}\quad \gamma_{2i-1}=-t_i+\beta_{n+i}+\psi_i.
\end{align*}

That is to say:
\begin{align}
 \displaystyle \gamma_i     &= \left\{     
\begin{array}{lr}
(-1)^i t_{\lceil i/2 \rceil}+\beta_{n+\lceil i/2\rceil-1}+\psi_{\lceil i/2 \rceil}& \mbox{if $b_{n+\lceil i/2\rceil-1} \geq 0$,} \\
(-1)^i t_{\lceil i/2\rceil}+\beta_{n+\lceil i/2\rceil}+\psi_{\lceil i/2 \rceil}& \mbox{if $b_{n+\lceil i/2\rceil-1} \leq 0$,}   
\end{array}
 \right.
 \end{align} as required.\end{proof}

\begin{remark}

We note that generalized Dynnikov coordinates for multicurves can be extended in a natural way to  generalized Dynnikov coordinates of measured foliations when the space of measured foliations on $N_{g,n}$ is endowed with its usual topology.
\end{remark}
\subsection{Generalized Dynnikov Coordinates on $N_{g,0}$}

Let $N_{g,0}$ be the standard model of a non--orientable surface of genus $g$ with one boundary component as shown in Figure \ref{arcsproof2} and denote by $\mathfrak{L}_{g,0}$ the set of multicurves on $N_{g,0}$.  Let
$\mathcal{S}_{g,n}=\{t;\,b;\, c)\in (\mathbb{Z}^{2(g-2)}\times \mathbb{Z}^g)\setminus \left\{0\right\}: |t_i|+\psi_i \quad \text{is even}\quad 1\leq i\leq g-2\}$. Let the function $\rho \colon \mathfrak{L}_{g,0}\to \mathcal{S}_{g,0}$  be defined by 
$$\rho(\mathcal{L}) = ( t;\,b;\, c):=(t_1,\ldots, t_{g-2};\, b_1,\, \ldots, b_{g-2};\, c_1,\ldots,\,c_g)$$ 

\noindent where 
\begin{align}\label{dynnikovformula2x}  
 \displaystyle t_{i}&= \frac{\gamma_{2i}-\gamma_{2i-1}}{2} \hspace{32pt} \textrm{;} \hspace{20pt} 1\leq i\leq g-2, \\
b_i&=\frac{\beta_i-\beta_{i+1}}{2} \hspace{42pt}  \textrm{;} \hspace{18pt} 1\leq i\leq g-2\end{align}

  \begin{figure}[h!]
\begin{center}
\labellist
\small\hair 2pt

     \pinlabel {$\scriptstyle{\gamma_{1}}$} [ ] at  85 220
          \pinlabel {$\scriptstyle{\gamma_{2}}$} [ ] at  85 60

                     \pinlabel {$\scriptstyle{\gamma_{2i-1}}$} [ ] at  210 230
                    \pinlabel {$\scriptstyle{\gamma_{2i-2}}$} [ ] at  210 60
                    
                     \pinlabel {$\scriptstyle{\gamma_{2g-5}}$} [ ] at  370 220
                    \pinlabel {$\scriptstyle{\gamma_{2g-4}}$} [ ] at  372 70

     \pinlabel {$\scriptstyle{\beta_{1}}$} [ ] at  50 180
     \pinlabel {$\scriptstyle{\beta_{2}}$} [ ] at  145 200
     
          \pinlabel {$\scriptstyle{\beta_{i}}$} [ ] at  180 200

     \pinlabel {$\scriptstyle{\beta_{i+1}}$} [ ] at  288 200
     
          \pinlabel {$\scriptstyle{\beta_{g-1}}$} [ ] at  410 180

 \endlabellist  
 \includegraphics[scale=0.6]{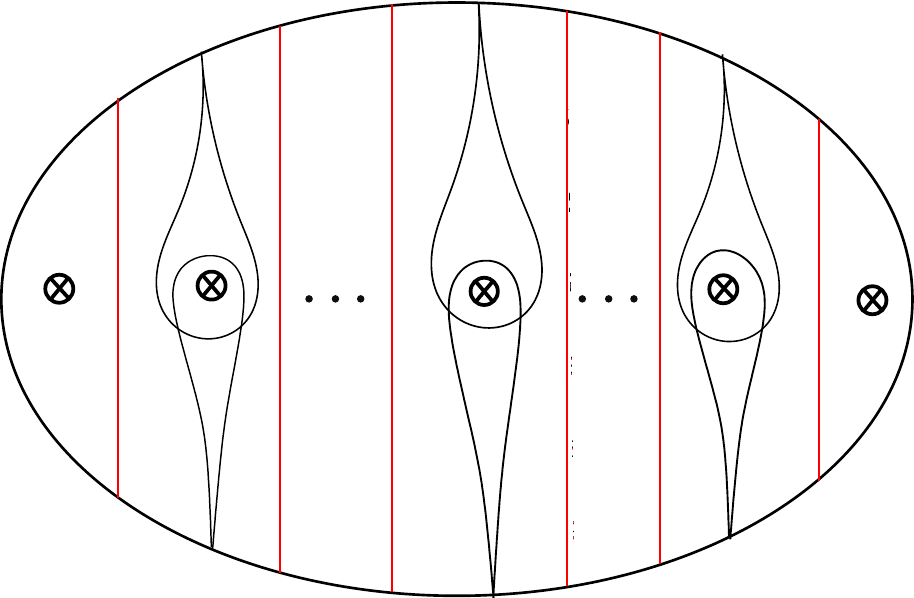}
\caption{The arcs $\gamma_i$, $\beta_i$  and the core curves $\cC_i$ on $N_{g,0}$}\label{arcsproof2}
\end{center}
\end{figure}

\begin{remark}\label{rem:loops3}

For $1\leq i\leq g-1$, $\psi_i$, $\lambda_i$  and $\lambda_{c_i}$ are as given in Lemma \ref{lem:loops2}. For $i=1$ we have  $\psi_1=0$,  $\lambda_1=\frac{\beta_1}{2}-c_1$ and
$\lambda_{c_1}=c_1$. Similarly, for $i=g$ we have  $\psi_g=0$ , $\lambda_g=\frac{\beta_{g-1}}{2}-c_g$ and $\lambda_{c_g}=c_g$ as each component of $L\cap S'_0$ and $L\cap S'_g$ intersecting the core curves should be core loop components.

\end{remark}

The inverse of the coordinate function $\rho$ is described similarly.  However, we need to extend the definition of almost boundary parallel arcs as they could also arise from the first crosscap as shown in Figure \ref{fig:figR1}. That is,  an {\it{almost boundary parallel arc}} associated with crosscap $g$ is a subarc of $L$ which intersects crosscap $g$ exactly once, has zero intersection with crosscaps $2$ through $g-1$;  and intersects  either crosscap 1 or the diameter between crosscap $1$ and the boundary exactly once as shown in Figure \ref{fig:figR1}(a) and Figure \ref{fig:figR1}(b). An {\it{almost boundary parallel arc}} associated with crosscap $1$ is described similarly.   

%An {\it{almost boundary parallel arc}}  associated with both crosscap $1$ and crosscap $g$ has zero intersection with the diameter and intersects only crosscap $1$ and crosscap $g$ exactly once. 

We write $R_1$ and $R_g$ for the number of almost boundary parallel arcs associated with crosscap $1$ and crosscap $g$  respectively.  Observe from Figure \ref{fig:figR1} that there are $R=\max(R_1, R_g)$ almost boundary components in total.  By the same argument as in the proof of Theorem  \ref{lem:inverse} we have  $\beta_i=Z_i+2R$ where

\begin{equation}
Z_i=\max_{1\leq s\leq g-1} \left\{2\max(b_{s}, 0)+ \psi_s +\left|t_s \right|+2\displaystyle\sum^{s-1}_{j=1}b_j\right\}- \displaystyle\sum^{i-1}_{j=1}b_j
 \end{equation}
\begin{figure}[h!]
\begin{center}
\labellist
 \pinlabel {${\scriptstyle{(a)}}$} [ ] at 65 -15
 \pinlabel {${\scriptstyle{(b)}}$} [ ] at  252 -15
 \pinlabel {${\scriptstyle{(c)}}$} [ ] at 452 -15
 \small\hair 2pt
\endlabellist  

  \includegraphics[scale=0.4]{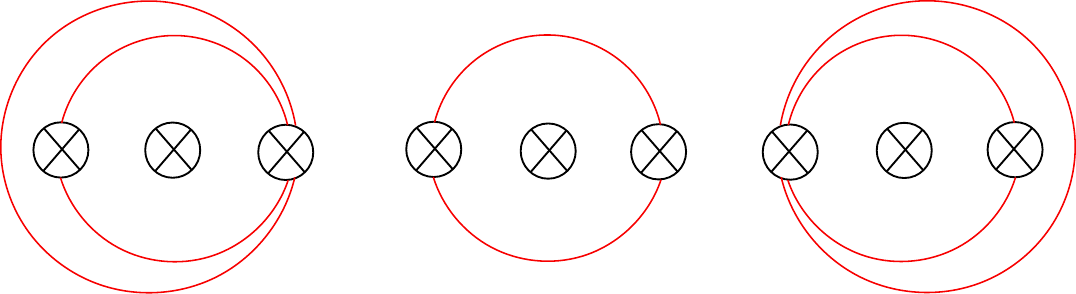}
\caption{Almost boundary arcs on $N_{3,0}$: (a) $R_1=1, R_{3}=2$;  (b) $R_1=R_{3}=1$; (c) $R_1=2, R_{3}=1$}\label{fig:figR1}
\end{center}
\end{figure}

 \noindent as there are no $a_i$ coordinates. Then, $\beta_1-2R= Z_1$ and  $\beta_{g-1}-2R= Z_{g-1}$.   We have $R_1=\max(0, c_1-\frac{Z_1}{2})$ and $R_{g}=\max(0, c_g-\frac{Z_{g-1}}{2})$ yielding $R=\max(0, c_1-\frac{Z_1}{2}, c_g-\frac{Z_{g-1}}{2} )$. The computation of intersection numbers on the arcs $\gamma_i$ is as in the proof of Theorem  \ref{lem:inverse}. Therefore we get Theorem \ref{lem:inverse2} where $\lceil x \rceil$ again denotes the smallest integer which is not less than $x$.

 \begin{remark}
 Note the special case $g=2$ where there are no $t_i$ and  $b_i$ coordinates. 
 \end{remark}

	\begin{thm}\label{lem:inverse2}
Let $ (t;\,b;\, c)\in \mathcal{S}_{g,0}$. Then $ (t;\, b;\, c)$ corresponds to a unique multicurve in $\cL \in \cL_{g, 0}$ which has 

\begin{align}
 \displaystyle \gamma_i  &= \left\{     
\begin{array}{lr}
(-1)^i t_{\lceil i/2 \rceil}+\beta_{n+\lceil i/2\rceil-1}+\psi_{\lceil i/2 \rceil}& \mbox{if $b_{n+\lceil i/2\rceil}-1 \geq 0$,} \label{eq:finalgammax}\\
(-1)^i t_{\lceil i/2\rceil}+\beta_{n+\lceil i/2\rceil}+\psi_{\lceil i/2 \rceil}& \mbox{if $b_{n+\lceil i/2\rceil}-1 \leq 0$,}   
\end{array}
 \right.\\
  \beta_i&= Z_i+2\max(0, c_g-\frac{Z_{g-1}}{2}, c_1-\frac{Z_{1}}{2})\label{eq:finalbetax}
 \end{align}
 
 \noindent  where
\begin{equation*}
Z_i=\max_{1\leq s\leq g-1} \left\{2\max(b_{s}, 0)+ \psi_s +\left|t_s \right|+2\displaystyle\sum^{s-1}_{j=1}b_j\right\}- \displaystyle\sum^{i-1}_{j=1}b_j
 \end{equation*}
 
	\end{thm}
	
\begin{thm}\label{thm:updatespecial}

Let $\cL\in \mathfrak{L}_{g,0}$ have generalized Dynnikov coordinates $( t;\,b;\, c)$.  Let $( t';\,b';\, c')$ and  $(t'';\,b'';\, c'')$ be the generalized Dynnikov coordinates of $u_i(\cL)$ and $u^{-1}_i(\cL)$ respectively. Then   $(c'_{i-1}, c'_{i})=(c_{i},c_{i-1})$ and $(c''_{i-1}, c''_{i})=(c_{i}, c_{i-1})$ for $1\leq i\leq g$; and for $1<i\leq g-2$ $(b';\,t')$ and  $(b'';\,t''')$ are as given in equations (\ref{eq:update}) replacing the subscript $i$ with $i-1$; for $i=g-1$ $t_{g-2}', b_{g-2}', t_{g-2}'', b_{g-2}''$ are as given in equations (\ref{eq:updatespecial}) replacing the subscripts $g-1$ and $n+g-2$ with $g-2$; and for $i=1$ we have

%b''_{n+i}&=\left[\frac{\bar{e_i}}{d_i}b^2_{n+i}(\frac{t''_{i}+t'_{i+1}}{t''_{i}})\frac{t_{i+1}}{t_{i}+t_{i+1}}\right]\quad\text{and}\quad b''_{n+i-1}=\left[\frac{b_{n+i}b_{n+i-1}}{b''_{n+i}}\right]
\begin{align*}\label{eq:update2}
  \begin{aligned}
  t_{1}' &=\left[\frac{t_{1}B_{1}}{f_0(1+B_1)+ t_1}\right]\\
   t_{1}'' &= \left[\frac{\bar{f_0}(t_1(B_1+1)+t_1)}{B_1}\right]\end{aligned} &&
  \begin{aligned}
  B_{1}' &= \left[B_{1}e_0\frac{(1+t'_1)}{t_1'}\frac{1}{(1+t_1)}\right]\\
  B_{1}'' &= \left[{\bar{e}_0}B_1\frac{(1+t''_1)}{t''_1}\frac{t_1}{(1+t_1)}\right]\end{aligned}
\end{align*}
In all other cases $t'_j=t_j$, $t''_j=t_j$ and $b'_j=b''_j=b_j$. 

\end{thm}

\section{Action of crosscap transpositions}\label{section3}
The goal of this section is to prove Theorem \ref{thm:update} which describes  how generalized Dynnikov coordinates change under the action of  $u_i$ and $u^{-1}_i$ ($1\leq i\leq g-1$). The key ingredient for the derivation of the formulae in Theorem \ref{thm:update} is a set of equalities associated with particular arc systems which we call {\it clovers} and {{\it scales}.  These equalities are given in Lemma \ref{lem:clover1}, Lemma \ref{lem:clover2}, Lemma \ref{lem:clover3} and Lemma \ref{lem:scales1}.

\subsection{Clovers and Scales}\label{clovers&scales}

A {\it{clover}} and a {\it{scale}}  about crosscaps $\otimes_i$ and $\otimes_{i+1}$   are  two different configurations consisting of  two vertices $v_1,\, v_2$ at $\partial N_{g,n}$ (identified to the puncture at $\infty$), five arcs $T_1,\,T_2,\,T_3,\,T_4,\,T_5$ with end points  $v_1,\, v_2$ and a curve $\cC_i$ such that the  teardrops $T_1$ and $T_2$ encircle $\otimes_i$, the teardrops $T_3$ and $T_4$ encircle $\otimes_{i+1}$, $T_5$ joins $v_1$ and $v_2$; and  $\cC_i$  is the essential simple closed curve bounding crosscaps $i$ and $i+1$ as shown in Figure \ref{fig:type12}.  We say that the clover has {\it{leaves}} $T_1,\,T_2,\,T_3$ and $T_4$;  {\it diagonal arc} $T_5$ and {\it diagonal curve} $\cC_i$. Also,  the scale has leaves $T_1,\,T_4,\,T_5$ and $\cC_i$; and {\it{diagonal teardrops}} $T_2$ and $T_3$.  

\begin{figure}[h!]
\begin{center}
\labellist
   \pinlabel {${T_1}$} [ ] at  82 210
      \pinlabel {${T_5}$} [ ] at  330 177

      \pinlabel {${T_2}$} [ ] at  5 103
   \pinlabel {$T_4$} [ ] at  82 16
   \pinlabel {${T_3}$} [ ] at   187 103
      \pinlabel {${v_1}$} [ ] at  159 177
   \pinlabel {${T_5}$} [ ] at  85 115

   \pinlabel {${C_i}$} [ ] at   156 107
\pinlabel {${v_2}$} [ ] at  28 38

% \pinlabel {${{X}}$} [ ] at  325 175
      \pinlabel {${{T_1}}$} [ ] at  220 103

      \pinlabel {${{T_2}}$} [ ] at  310 35
           \pinlabel {${{T_3}}$} [ ] at  343 35

   \pinlabel {${C_i}$} [ ] at  330 -5
   \pinlabel {${{T_4}}$} [ ] at   422 103

   %\pinlabel {${X}$} [ ] at  81 107
      \pinlabel {${v_2}$} [ ] at  390 177

 \pinlabel {${v_1}$} [ ] at  260 177

\small\hair 2pt
\endlabellist  

  \includegraphics[scale=0.7]{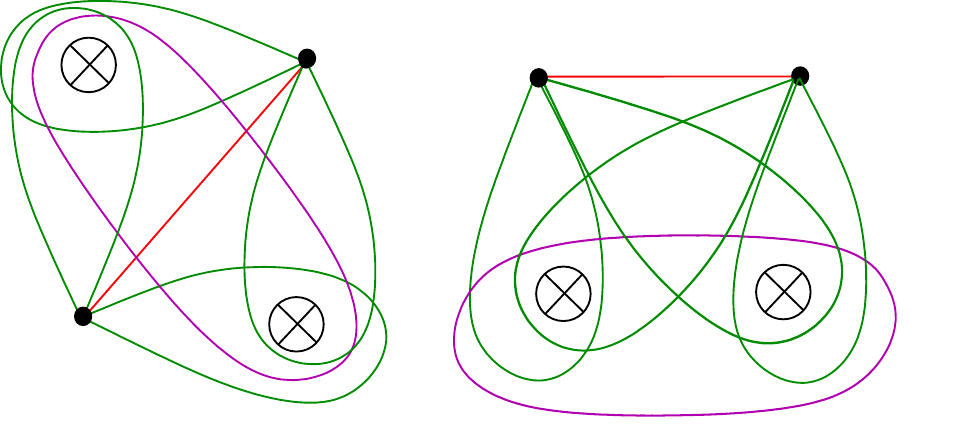}
\caption{A clover and a scale} \label{fig:type12}
\end{center}
\end{figure}

         To compute the action of $u_i$ and its inverse $u^{-1}_i$ ($1\leq i\leq g-1$) in terms of generalized Dynnikov coordinates we will make use of certain equalities associated with {clovers} and {scales} in $S'_i\cup S'_{i+1}$, which we shall call {\it{clover and scale equalities}} throughout.   We shall consider  three types of clovers and four types of scales to obtain clover and scale equalities. These clovers and scales are depicted in Figure \ref{fig:cloverI&II} and  Figure \ref{fig:scaleI&II} respectively.

Clover and scale equalities can be considered as a generalization of a well known equation commonly known as the {\it flip move} which lets us change coordinates from one triangulation to another on punctured orientable surfaces \cite{bell, thurston}. 
Namely, if $Q$ is a rectangle in a surface $S$ and $X_1, X_2, X_3, X_4, X_5$ and $X_6$ are the number of intersections  on the four edges and the diagonals of $Q$ with all of its corners at the punctures and containing no punctures in its interior and $X_{ij}$ denotes the number of arcs joining edge $X_i$ to $X_j$ then there are two possibilities: either $X_{12}$ or $X_{34}$ is zero since curves are non intersecting. This yields the well known equation  $$\max(X_1+X_2, X_3+X_4)=X_5+X_6.$$   The method here will be similar and use case by case analysis for components of $L\cap(S'_i\cup S'_{i+1})$ intersecting the clovers and the scales. 

In what follows we will again denote by $L$ a minimal representative of a multicurve $\cL\in \mathfrak{L}_{g,n}$ with $\rho(\cL)=(a, b; t, c)$ and intersection numbers $(\alpha, \beta; \gamma,  c)$.

%The coordinates related with a single arc have been marked by a ``hat''  in Section \ref{sec:standardarcs} so as to distinguish them from coordinates of $L_i$. Here we use similar notation, and label the parameters given  in Notation \ref{not:deltas} with a ``hat'' in Notation \ref{not:twist} to indicate that they are associated with a single arc $X$ rather than a family of arcs  $L_i$.  

\subsubsection{Components of $L\cap(S'_i\cup S'_{i+1})$}\label{sec:standardarcs}
We can list all topological possibilities for connected components of  $L\cap (S'_i\cup S'_{i+1})$ ($1\leq i \leq g-1$), up to isotopy, making use of their intersections with the core curves $c_i$ and $c_{i+1}$,  and  the arcs $\gamma_l$ ($2i-1\leq l\leq 2i+1$) and $\beta_l$ ($n+l-1\leq l\leq n+l+1$) (and hence from their generalized Dynnikov coordinates). Let $\otimes_i$ and $\otimes_{i+1}$ denote crosscap $i$ and crosscap $i+1$ respectively.   Given a connected component $X$ of  $L\cap (S'_i\cup S'_{i+1})$ we associate with it a  {\it{signature vector}} $\widehat{v_i}=(\widehat{c_i}, \widehat{c}_{i+1}; \widehat{b}_{n+i-1}, \widehat{b}_{n+i})\in \mathbb{Z}_{\geq0}^2\times\mathbb{Z}^2$ such that  $\widehat{c_i}$ and $\widehat{c}_{i+1}$; give the number of intersections between $X$ and the  core curves of $\otimes_i$ and $\otimes_{i+1}$ respectively, and for $j=i, i+1$

\begin{align} \widehat{b}_j=\frac{\widehat{\beta}_{n+j-1}-\widehat{\beta}_{n+j}}{2}\end{align}

\noindent where $\widehat{\beta_{k}}$  denote the number of intersections of $X$ with $\beta_k$ $(k=n+j-1, n+j)$.
\begin{remark} Note that each signature vector $\widehat{v_i}$ must  satisfy equation (\ref{eq:relation1b}) of Lemma \ref{lem:equalities}, and the inequality $\sgn( \widehat{b}_{n+i-1})\sgn( \widehat{b}_{n+i})\leq 0$ since  $X$ is a connected component of  $L\cap (S'_i\cup S'_{i+1})$.\end{remark} Then each connected component of $L\cap (S'_i\cup S'_{i+1})$ is either a simple closed curve  supported in $S'_i\cup S'_{i+1}$, or  one of the following arcs described in Notation \ref{not:components}. 

\begin{figure}[h!]
\begin{center}
\labellist
 \pinlabel {${\scriptstyle{\textcolor{red}{1}}}$} [ ] at  140 568

  \pinlabel {${\scriptstyle{\textcolor{red}{2}}}$} [ ] at 140 535
   \pinlabel {${\scriptstyle{\textcolor{red}{3}}}$} [ ] at  230 475
 \pinlabel {${\scriptstyle{\textcolor{red}{4}}}$} [ ]  at  235 440

 \pinlabel {${\scriptstyle{\textcolor{red}{5}}}$} [ ] at   530 590

 \pinlabel {${\scriptstyle{\textcolor{red}{6}}}$} [ ]  at  540 550
 \pinlabel {${\scriptstyle{\textcolor{red}{7}}}$} [ ]  at  500 525

 \pinlabel {${\scriptstyle{\textcolor{red}{8}}}$} [ ]   at  550 390

 \pinlabel {${\scriptstyle{\textcolor{red}{9}}}$} [ ]  at  140 210
 \pinlabel {${\scriptstyle{\textcolor{red}{10}}}$} [ ]  at  140 180
 \pinlabel {${\scriptstyle{\textcolor{red}{11}}}$} [ ]  at  158 140
 \pinlabel {${\scriptstyle{\textcolor{red}{12}}}$} [ ]  at  440 225
 \pinlabel {${\scriptstyle{\textcolor{red}{13}}}$} [ ]   at  350 130

  \pinlabel {${\scriptstyle{\textcolor{red}{14}}}$} [ ]  at  390 95
       
\pinlabel {${\scriptstyle{\textcolor{red}{15}}}$} [ ]   at  560 100
 \pinlabel {${\scriptstyle{\textcolor{red}{16}}}$} [ ]  at 445 55

 \pinlabel {${\scriptstyle{\gamma_{2i-1}}}$} [ ] at  75 290
  \pinlabel {${\scriptstyle{\gamma_{2i+1}}}$} [ ] at  180 290

 \pinlabel {${\scriptstyle{\gamma_{2i}}}$} [ ] at  80 4
  \pinlabel {${\scriptstyle{\gamma_{2i+2}}}$} [ ] at  185 4
   \pinlabel {${\scriptstyle{\gamma_{2i}}}$} [ ] at  405 4
  \pinlabel {${\scriptstyle{\gamma_{2i+2}}}$} [ ] at  505 4
  
   \pinlabel {${\scriptstyle{\gamma_{2i-1}}}$} [ ] at  75 600
  \pinlabel {${\scriptstyle{\gamma_{2i+1}}}$} [ ] at  180 600
 \pinlabel {${\scriptstyle{\gamma_{2i}}}$} [ ] at  80 320
  \pinlabel {${\scriptstyle{\gamma_{2i+2}}}$} [ ] at  182 320
  
    \pinlabel {${\scriptstyle{\gamma_{2i-1}}}$} [ ] at  400 605
  \pinlabel {${\scriptstyle{\gamma_{2i+1}}}$} [ ] at  502 605
 \pinlabel {${\scriptstyle{\gamma_{2i}}}$} [ ] at  403 320
  \pinlabel {${\scriptstyle{\gamma_{2i+2}}}$} [ ] at  510 320

 \pinlabel {${\scriptstyle{\gamma_{2i-1}}}$} [ ] at  400 290
  \pinlabel {${\scriptstyle{\gamma_{2i+1}}}$} [ ] at  500 290

   \pinlabel {${\scriptstyle{\beta_{n+i-1}}}$} [ ] at  10 290
  \pinlabel {${\scriptstyle{\beta_{n+i}}}$} [ ] at  120 290
  
    \pinlabel {${\scriptstyle{\beta_{n+i+1}}}$} [ ] at  250 290

   \pinlabel {${\scriptstyle{\beta_{n+i-1}}}$} [ ] at  340 290
  \pinlabel {${\scriptstyle{\beta_{n+i}}}$} [ ] at  450 290
  \pinlabel {${\scriptstyle{\beta_{n+i+1}}}$} [ ] at  590 290
  
    \pinlabel {${\scriptstyle{\beta_{n+i-1}}}$} [ ] at  340 608
  \pinlabel {${\scriptstyle{\beta_{n+i}}}$} [ ] at  450 608
  \pinlabel {${\scriptstyle{\beta_{n+i+1}}}$} [ ] at  590 608

 \pinlabel {${\scriptstyle{\beta_{n+i-1}}}$} [ ] at  10 608
  \pinlabel {${\scriptstyle{\beta_{n+i}}}$} [ ] at  120 608
  \pinlabel {${\scriptstyle{\beta_{n+i+1}}}$} [ ] at  250 608

\small\hair 2pt
\endlabellist  

  \includegraphics[scale=0.49]{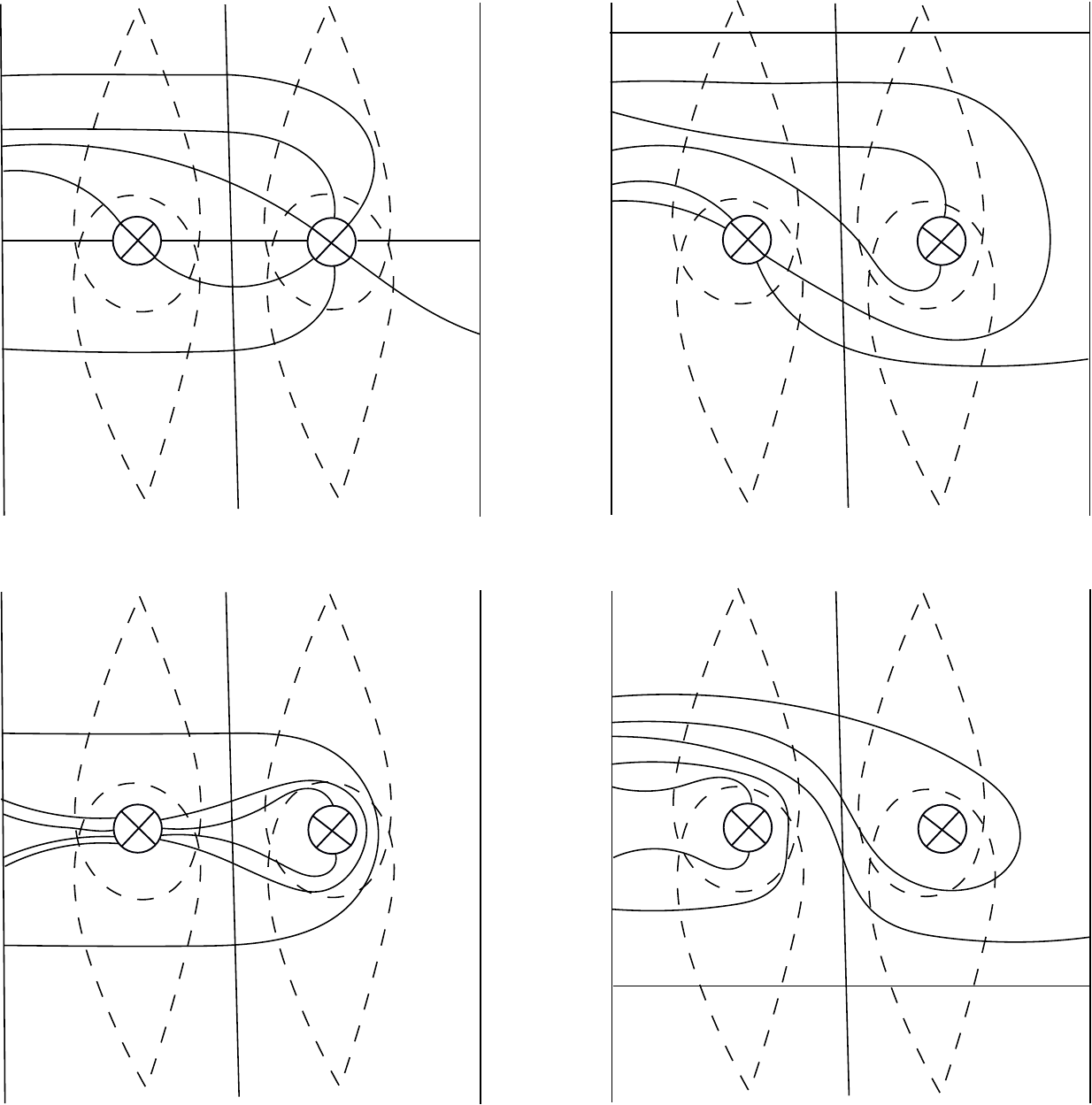}
\caption{$\textcolor{red}{1}:X^{(1,1;0,1)}_{2i-1; i}$, 
$\textcolor{red}{2}:X^{(0,1; 0,1)}_{2i-12i}$,
$\textcolor{red}{3}:X_{ii+1}$,
$\textcolor{red}{4}:X_{2i-1; \,i+1}$,
$\textcolor{red}{5}: X_{2i-12i+1}$,
$\textcolor{red}{6}: X^{(1,0; 0,1)}_{2i-1; i}$,
$\textcolor{red}{7}: X^{(0,1; 0,1)}_{2i-12i-1}$,
$\textcolor{red}{8}: X_{2i+2;i}$,
$\textcolor{red}{9}: X_{2i-12i}^{(0,0; 0,1)}$,
$\textcolor{red}{10}: X_i^{(2,0; 0,1)}$,
$\textcolor{red}{11}: X_i^{(2,1; 0,1)}$,
$\textcolor{red}{12}: X^{(0,0; 0,1)}_{2i-12i-1}$,
$\textcolor{red}{13}: X_i$,
$\textcolor{red}{14}: X_{2i-12i}^{(0,0;1,0)}$,
$\textcolor{red}{15}: X_{2i-12i+2}$,
$\textcolor{red}{16}: X_{2i2i+2}$} \label{fig:standardarcs}
\end{center}
\end{figure}
 %with signature vector $\widehat{v_i}=(\widehat{c_i}, \widehat{c}_{i+1}; \widehat{b}_{n+i-1}, \widehat{b}_{n+i})$

\begin{notation}\label{not:components}
Let $2i-1 \le l\le  m \le  2i+2 $ and $k=i, i+1$.  
\begin{itemize}
\item[1.] ${X}^{v_i}_{lm}$:
\begin{itemize} 
\item[(a)] If $(l, m)\neq (2i-1,2i), (2i+1,2i+2)$ it passes above (resp. below) $\otimes_i$ if $l=2i-1$ (resp. $l=2i$), and it passes above (resp. below) $\otimes_{i+1}$ if $m=2i+1$ (resp. $m=2i+2$). It has one end point on  $\beta_{n+i-1}$ and the other on  $\beta_{n+i+1}$.

\item[(b)]  If $(l, m)= (2i-1, 2i)$  it passes both above  and below $\otimes_i$, and it has both end points on $\beta_{n+i-1}$. The case  $(l, m)= (2i+1, 2i+2)$ is described similarly. 

\item[(c)]  If $l=m$  it passes above (resp. below) $\otimes_i$ if $l=2i-1$ (resp. $l=2i$)  and has both end points on $\beta_{n+i-1}$. The case  $l=m=2i+1$ is described similarly. 

\end{itemize}

 \item[2.]${X}^{v_i}_{l; k}$: it passes above (resp. below) $\otimes_i$ if $l=2i-1$ (resp. $l=2i$). It has  $\psi_i=1$ and both end points on  $\beta_{n+i-1}$. The cases $l=2i+1$ and $l=2i+2$ are described similarly. 

\item[3.]  ${X}_i$,  ${X}_i^{v_i}$,  ${X}_{i+1}$, ${X}_{i+1}^{v_i}$ and ${X}_{ii+1}$:   none of these components pass above or below $\otimes_i$ and $\otimes_{i+1}$.  ${X}_i$ and ${X}_i^{v_i}$ have both end points on $\beta_{n+i-1}$ with $\psi_i=0$ and $\psi_i\neq 0$  respectively.  ${X}_{i+1}$ and ${X}_{i+1}^{v_i}$ are described similarly. ${X}_{ii+1}$ has one end point on  $\beta_{n+i-1}$ and the other on  $\beta_{n+i+1}$.
\end{itemize}

See for example the arcs 5,15,16 for item 1(a); 2, 9, 14 for item 1(b);  7, 12 for item 1(c);  1,4, 6, 8  for item 2; and 3, 10,11,13 for item 3.  in Figure \ref{fig:standardarcs}.
\end{notation}

\begin{notation}\label{not:omit}

We shall omit the superscript $\widehat{v}_i$ for when  $\widehat{b}_{n+i-1}=\widehat{b}_{n+i}=0$. See for example the arcs 3, 4, 5, 8, 15 and 16 in Figure \ref{fig:standardarcs}. We write  $[P_i]$ for the set of simple closed curves in $S'_i\cup S'_{i+1}$, and  $[{X}^{v_i}_{lm}]$, $[{X}^{v_i}_{l; k}]$ and $[{X}_i]$ for  the set of arcs described in item 1, item 2. and item 3 respectively. \end{notation}
		%  Observe that $A'_k$ and $B'_k$~($k=i, i+1$) can be determined from the subscripts of $X^{v_i}_{l; k}$ and $X_{l, m}^{v_i}$. If $A'_k\neq 0$ then $A'_k=1$. Similar argument holds for $B'_k$. 

\begin{notation}\label{not:compatibleL}
Here and in what follows we write $L_i$ to denote the arc system $L\cap(S'_i\cup S'_{i+1})$ for convenience.
\end{notation}
 \begin{defn}\label{def:standard&twisted}

 Let $X$ be a component  of $L_i$ with signature vector $\widehat{v_i}=(\widehat{c_i}, \widehat{c}_{i+1}; \widehat{b}_{n+i-1}, \widehat{b}_{n+i})$.  We say that $X$ is a {\it standard} arc if at least one of $\widehat{b}_{n+i-1}$ and $ \widehat{b}_{n+i}$ equals zero (Figure \ref{fig:standardarcs}). We say that $X$ is {\it twisted} if  $\widehat{b}_{n+i-1}<0$  and $\widehat{b}_{n+i}>0$. \end{defn}

   %This is equivalent to saying that each twisted component  is the image of some standard component under a mapping class supported in  $S'_i\cup S'_{i+1})$ \cite{}. 

%\begin{remark} Let $X$ be a connected component of $L\cap (S'_i\cup S'_{i+1})$. Then, 
%\begin{itemize}
%\item It satisfies (2.2) of Lemma \ref{lem:equalities}. 
%\item  If it is standard at least one of $b_i$ and $b_{i+1}$ equals zero.
%\item $A'_k$ and $B'_k$~($k=i, i+1$) can be determined from the subscripts of $X^{v_i}_{l; k}$ and $X_{l, m}^{v_i}$. If $A'_k\neq 0$ then $A'_k=1$. Similar argument holds for $B'_k$. 

%\end{itemize}

%\end{remark}

 \begin{defn}\label{def:positive}
 We say that a component of $L_i$ is a {\it positive} arc if it satisfies $t_{i}-t_{i+1}>0$, it is {\it negative} if  $t_{i}-t_{i+1}<0$ and that it is neutral if $t_{i+1}-t_i=0$. \end{defn}
 
Figure \ref{fig:standardarcs} depicts examples for negative and neutral standard components of $L_i$. Other standard components can be obtained by symmetry, reflecting them in the arc  $\beta_{n+i}$, or the diameter of the surface. 

 \begin{remark}\label{rem:curvet}
 Since $\widehat{b}_{n+i-1}<0$  and $\widehat{b}_{n+i}>0$, and $t_{i}-t_{i+1}=0$ for a simple closed curve in $S'_i\cup S'_{i+1}$ we regard each such curve neutral and twisted. 
\end{remark}

\begin{defn}\label{defn:halftwisted}
A twisted component of $L_i$ is called a  {\it negatively (resp. positively) half twisted}  arc if it is the $u^{-1}_i$  (resp. $u_i$)  image of a standard arc.  A negative (resp. positive)  twisted component of $L_i$ which is not half twisted is called {\it negatively (resp. positively) highly twisted}. Each simple closed curve and twisted neutral component of $L_i$ is called {\it neutrally twisted}.  
\end{defn}
 %Note that  a half twisted component may still be standard.  For example the standard arc $6$ is the $u^{-1}_i$ image of the standard arc labeled 2  in Figure \ref{fig:standardarcs}.

 \begin{notation}\label{not:lambdas} Let $\lambda_k$ denote the number of non--core loop components of $L\cap S'_k$ ($k=i, i+1$). We denote

\begin{equation*}
\lambda^+_{i+1}= \begin{cases}
             \lambda_{i+1}  & \text{if } b_{n+i} >0\\
             0 & \text{if } b_{n+i} <0
       \end{cases} \quad \text{and} \quad
\lambda^-_{i}= \begin{cases}
             \lambda_{i}  & \text{if } b_{n+i-1} <0\\
             0 & \text{if } b_{n+i-1} >0
       \end{cases}
\end{equation*}
We describe $\lambda^+_{c_{i+1}} $ and $\lambda^-_{c_{i}}$ similarly. Elsewhere $x^+$ denotes $\max(x, 0)$. \end{notation}

  Let $A'_i$ and $B'_i$ be the number of above and below components of $L_i$ given in Lemma \ref{lem:abovebelow}.

\begin{notation}\label{not:deltas}
Let  $A'_{i,i+1} = \min(A'_i, A'_{i+1})$ \text{and} $B'_{i,i+1} = \min(B'_i, B'_{i+1})$.
We write 
 \begin{align} \Delta_i(A)=A'_i-A'_{i,i+1} \qquad \text{and}  \qquad  \Delta_{i+1}(A)=A'_{i+1}-A'_{i,i+1}\\
 \Delta_i(B)=B'_i-B'_{i,i+1} \qquad \text{and}  \qquad  \Delta_{i+1}(B)=B'_{i+1}-B'_{i,i+1}
 \end{align}
 
%We write $\widehat{\Delta}_k(A)$ and $\widehat{\Delta}_k(B)$  ($k=i, i+1$) to associate these parameters with a given single arc $X$ rather than a family of arcs of $L\cap (S'_i\cup S'_{i+1})$. 

\end{notation}
 
\begin{remark}\label{rem:deltas}
Geometrically, $A'_{i, i+1} $ and $B'_{i,i+1}$ give the number of components of $L_i$ which  lie entirely above and below the diameter of the surface respectively. Therefore, $A'_{i, i+1}=X_{2i-12i+1}$ and $B'_{i, i+1}=X_{2i2i+2}$. Then, $\Delta_i(A)$ (resp. $\Delta_i(B)$) is the number of above (resp. below) components of $L\cap S'_{i}$ which are not  contained in the arcs $X_{2i-1, 2i+1}$ (resp.  $X_{2i, 2i+2}$).    $\Delta_{i+1}(A)$ and  $\Delta_{i+1}(B)$ are interpreted similarly. 
\end{remark}
Figure \ref{fig:motivation} illustrates these parameters which we will refer to later to describe other parameters. 

\begin{figure}[h!]
\begin{center}
\labellist
 \pinlabel {${\scriptstyle{\textcolor{red}{A'_{i}=1}}}$} [ ] at  50 270
 \pinlabel {${\scriptstyle{\textcolor{red}{A'_{i+1}=2}}}$} [ ] at  190 270
 \pinlabel {${\scriptstyle{\textcolor{red}{\lambda'_{i}=2}}}$} [ ] at  50 80

 \pinlabel {${\scriptstyle{\textcolor{red}{B'_{i}=0}}}$} [ ] at  50 40
 \pinlabel {${\scriptstyle{\textcolor{red}{B'_{i+1}=2}}}$} [ ] at  190 40

 \pinlabel {${\scriptstyle{\textcolor{red}{\psi_{i+1}=1}}}$} [ ] at  240 155
 \pinlabel {${\scriptstyle{\textcolor{red}{\lambda^+_{c_{i+1}}=1}}}$} [ ] at  240 200

 \pinlabel {${\scriptstyle{\textcolor{red}{\lambda^-_{c_{i}}=1}}}$} [ ] at  62 173

\small\hair 2pt
\endlabellist  

  \includegraphics[scale=0.49]{motivation}
\caption{$A'_{i,i+1}=1, B'_{i,i+1}=0, \Delta_i(A)= \Delta_i(B)=0, \Delta_{i+1}(A)=1, \Delta_{i+1}(B)=2$} \label{fig:motivation}
\end{center}
\end{figure}

\begin{remark}\label{rem:posneg}
The condition $t_{i}-t_{i+1}>0$ in Definition \ref{def:positive}  implies that either $\Delta_i(B)\neq 0$ or $\Delta_{i+1}(A)\neq 0$ holds for a positive component. Then,  ${X}^{v_i}_{lm}$ is positive  if $l=2i+1$ or $m=2i$. Similarly, ${X}^{v_i}_{l; k}$ ($k=i, i+1$) is positive if $l=2i+1$ or $l=2i$. Similar arguments hold for negative components. Similarly, the condition $t_{i+1}-t_i=0$ implies that a neutral component in item 1(a) (i.e. $X_{2i-12i+1}$ and $X_{2i2i+2}$)  and  item 2. has $\Delta_k(A)=\Delta_k(B)=0$~($k=i, i+1$). A neutral component in item 1(b) has either $\Delta_i(A)\neq 0$ and $\Delta_i(B)\neq 0$ (i.e. $X^{v_i}_{2i-12i}$)  or $\Delta_{i+1}(A)\neq 0$ and $\Delta_{i+1}(B)\neq 0$  (i.e. $X^{v_i}_{2i+12i+2}$). \end{remark}

\begin{defn}\label{def:compatible}
 We say that two arcs in $S'_i\cup S'_{i+1}$ are compatible if they can be embedded disjointly in  $S'_i\cup S'_{i+1}$.  \end{defn}

A positive and a negative arc are compatible only when they form the arc systems {\it scissors},  {\it anchors} or  {\it ribbons}. 
   
\begin{figure}[h!]
\begin{center}
\labellist

 \pinlabel {${\scriptstyle{\beta_{n+i-1}}}$} [ ] at  4 288
  \pinlabel {${\scriptstyle{\beta_{n+i}}}$} [ ] at  118 288

 \pinlabel {${\scriptstyle{\beta_{n+i+1}}}$} [ ] at  260 288
 
 \pinlabel {${\scriptstyle{\beta_{n+i-1}}}$} [ ] at  350 288
  \pinlabel {${\scriptstyle{\beta_{n+i}}}$} [ ] at 465 288

 \pinlabel {${\scriptstyle{\beta_{n+i+1}}}$} [ ] at  580 288
 
  \pinlabel {${\scriptstyle{\beta_{n+i-1}}}$} [ ] at  680 288
  \pinlabel {${\scriptstyle{\beta_{n+i}}}$} [ ] at 780 288

 \pinlabel {${\scriptstyle{\beta_{n+i+1}}}$} [ ] at  900 288
  \pinlabel {(a)} [ ] at  120 -10
    \pinlabel {(b)} [ ] at  470 -10
      \pinlabel {(c)} [ ] at  780 -10
 \small\hair 2pt
\endlabellist  

  \includegraphics[scale=0.3]{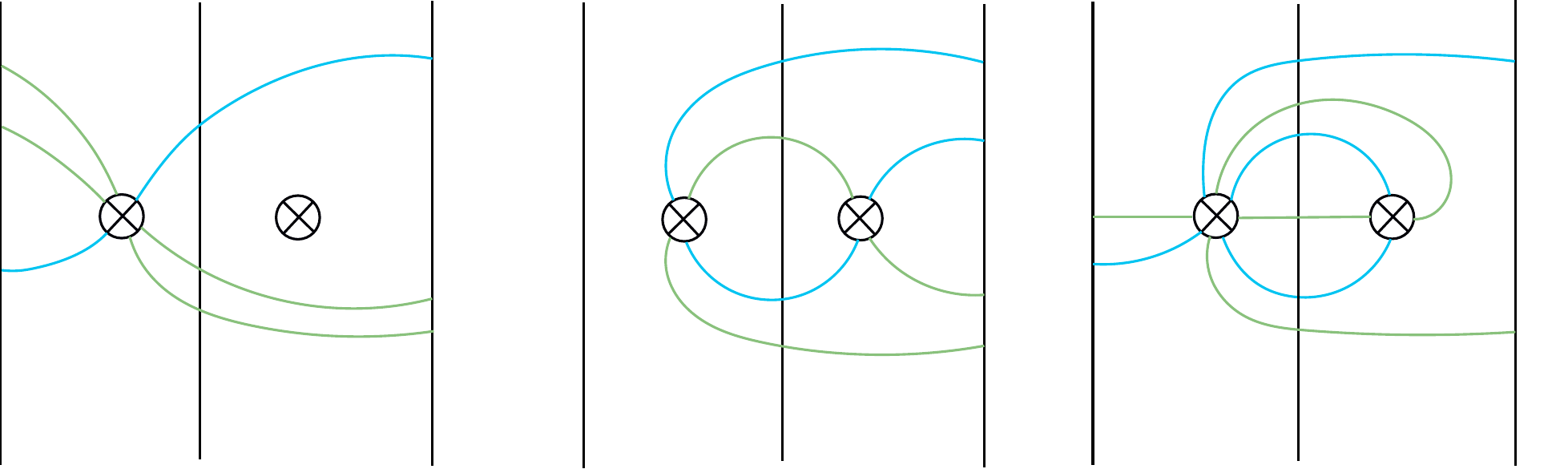}
\caption{Scissors  at $\otimes_i$ with $X_{2i+2;\, i}=2$  and $X_{2i+1;\, i}=1$ , a left anchor with  $X^{(1,1; -1,0)}_{2i+2,i+1}=1$  and $X^{(1,1; -1,0)}_{2i+1,i+1}=1$; a left ribbon with $X^{(2,1; -1,1)}_{2i+2; i}=1$ and $X^{(2,1; -1,1)}_{2i+1; i}=1$} \label{fig:compatiblesxxx}
\end{center}
\end{figure}

 \begin{defn}\label{def:scissorsanchors} {\it Scissors} at $\otimes_i$ consists of the arcs  $X_{2i+2;\, i}$  and $X_{2i+1;\, i}$ (Figure  \ref{fig:compatiblesxxx}(a)), and {{scissors}} at  $\otimes_{i+1}$ consists of the arcs $X_{2i;\, i+1}$ and  $X_{2i-1;\, i+1}$. 
 
 A {{left anchor}}  consists of  $X^{(1,1; -1,0)}_{2i+2,i+1}$  and $X^{(1,1; -1,0)}_{2i+1,i+1}$ (see Figure  \ref{fig:compatiblesxxx}(b)), and a {{right anchor}}  consists of $X^{(1,1; 0,1)}_{2i; i}$ and $X^{(1,1; 0,1)}_{2i-1; i}$.
   
 			A  {{left ribbon}}  consists of arcs from the sets $[X_{2i+2;  i}]$ and $[X_{2i+1; i}]$ (arcs of a ribbon are twisted and can have different signature vectors). And a {{right ribbon}}   consists of elements from the sets $[X_{2i; i+1}]$ and $[X_{2i-1; i+1}]$.  Note that scissors, ribbons and anchors may contain multiple copies of the same arc.
			
			The positive and negative arcs of given scissors are respectively called positive and negative arms for the scissors. For example,   $X_{2i+2;\, i}$ and $X_{2i+1;\, i}$ are  negative and positive arms of the scissors respectively. Positive and  negative arms of anchors and ribbons are described similarly. 
\end{defn}

    \begin{remark}\label{rem:poss} 
    Note that scissors, anchors and ribbons are not compatible with each other.  In fact if a component of $L_i$ is compatible with any  scissors , anchors or ribbons it must be a neutral arc. 
 \end{remark}
%none of these three arc systems are compatible with each other.  For the same reason scissors at $\otimes_i$  and scissors at $\otimes_{i+1}$ are not compatible. Similarly,  a right anchor is not compatible with a left anchor as a right ribbon is not compatible with a left ribbon. 

\begin{defn}\label{def:compatible}

If $L_i$  is comprised of standard arcs it is called {\it standard}, and if it contains at least one twisted arc it is called  {\it twisted}. If it contains positive arcs possibly with neutral arcs it  is called positive. The case when $L_i$ is negative is described similarly. If $L_i$ contains only neutral arcs it is called {\it neutral}, and if it contains both positive and negative arcs it is called {\it mixed}.    \end{defn}  

\noindent Therefore $L_i$ is mixed if  and only if it contains either scissors or an anchor or a ribbon, and that the only case $L_i$ is both twisted and mixed is when it contains ribbons.

         \subsubsection{Clover Equalities}\label{sec:clovers}
          Let  $u^{-1}_i(\gamma;\,\beta)=(\gamma';\,\beta')$ and $u_i(\gamma;\,\beta)=(\gamma'';\,\beta'')$.   A clover of type I has leaves $\gamma_{2i-1}, \gamma_{2i}, \gamma_{2i+1}$, $\gamma_{2i+2}$, the diagonal arc $\beta_{n+i}$  and the diagonal curve $\cC_i$  (Figure \ref{fig:cloverI&II}(a)).    A clover of type II and a clover of type III are the images of a clover of type I under the mapping classes $u^{-1}_i$ and $u_i$ respectively, and hence are as depicted in Figure \ref{fig:cloverI&II}(b) and Figure  \ref{fig:cloverI&II}(c) respectively.
%         
         %and hence it has leaves $\gamma_{2i}$, $\gamma_{2i+1}$, $\gamma'_{2i}$, $\gamma'_{2i+1}$; the diagonal arc $\beta'_{n+i}$ and the diagonal curve $\cC_i$ (note that $\cC_i$ is fixed under the action of $u_i$ and $u_i^{-1}$). 

          We present equalities associated with a clover of type I, type II and type III given in Lemma \ref{lem:clover1},  Lemma \ref{lem:clover2} and Lemma \ref{lem:clover3} respectively. Here and in what follows we abuse notation again using the symbols in Notation \ref{not:components} to denote the number of corresponding components  of $L_i$, and the symbols $\gamma;\,\beta, \gamma';\,\beta', \gamma'';\,\beta''$ to denote the number of intersections.  First, we  fix some notation which will be necessary in the proof of Lemma \ref{lem:clover1}. 
          
           \begin{figure}[h!]
\begin{center}
\labellist
 \pinlabel {${\scriptstyle{\gamma_{2i-1}}}$} [ ] at  50 230
  \pinlabel {${\scriptstyle{\gamma_{2i+1}}}$} [ ] at  260 230
 \pinlabel {${\scriptstyle{\gamma_{2i}}}$} [ ] at  70 40
  \pinlabel {${\scriptstyle{\gamma_{2i+2}}}$} [ ] at  260 40
     \pinlabel {\begin{turn}{-90}$\scriptstyle{\beta_{n+i}}$\end{turn}} [ ] at  140 250
      \pinlabel {${\scriptstyle{C_{i}}}$} [ ] at  280 150
 \pinlabel {${\scriptstyle{\gamma'_{2i}}}$} [ ] at  500 250
% \pinlabel {${\scriptstyle{\beta_{n+i-1}}}$} [ ] at  378 150
 %\pinlabel {${\scriptstyle{\beta_{n+i+1}}}$} [ ] at  780 150
 \pinlabel {${\scriptstyle{\gamma_{2i}}}$} [ ] at  500 45
 \pinlabel {${\scriptstyle{\gamma_{2i+1}}}$} [ ] at  620 240
      \pinlabel {\begin{turn}{-90}$\scriptstyle{\beta'_{n+i}}$\end{turn}} [ ] at  410 60

 \pinlabel {${\scriptstyle{\gamma'_{2i+1}}}$} [ ] at  667 85
 
       \pinlabel {\begin{turn}{-90}$\scriptstyle{\gamma''_{2i+2}}$\end{turn}} [ ] at   1115 60

% \pinlabel {${\scriptstyle{\beta_{n+i-1}}}$} [ ] at  378 150
 %\pinlabel {${\scriptstyle{\beta_{n+i+1}}}$} [ ] at  780 150
 \pinlabel {${\scriptstyle{\gamma_{2i+2}}}$} [ ] at  965 45
 \pinlabel {${\scriptstyle{\gamma_{2i-1}}}$} [ ] at  940 240
      \pinlabel {\begin{turn}{-90}$\scriptstyle{\beta''_{n+i}}$\end{turn}} [ ] at  1060 60

 \pinlabel {${\scriptstyle{\gamma''_{2i-1}}}$} [ ] at  760 240

\small\hair 2pt
\endlabellist  
 \includegraphics[scale=0.4]{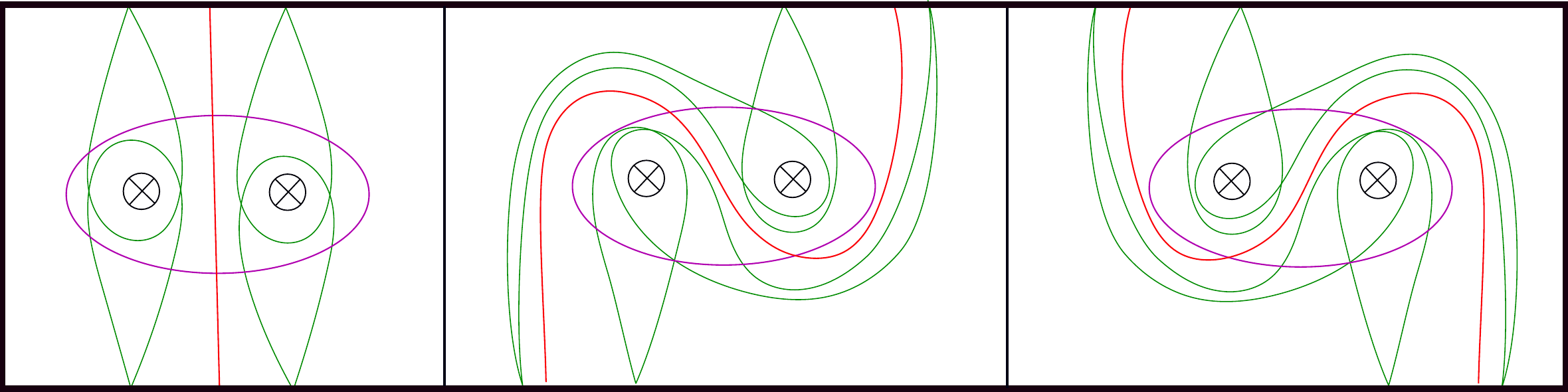}
\caption{(a)~A clover of type I, (b)~a clover of type  II and (c)~a clover of type III}\label{fig:cloverI&II}
\end{center}
\end{figure}

\begin{notation}\label{not:clover1}
Given scissors at $\otimes_i$ and scissors at $\otimes_{i+1}$ let $s_i=\min(X_{2i+1; \,i}, X_{2i+2; \,i})$ and $s_{i+1}=\min(X_{2i-1;\,i+1}, X_{2i;\, i+1})$ respectively.  Given a right and a left anchor in   $S'_i\cup S'_{i+1}$  let  $z_{i}=\min(X^{(1,1; 0,1)}_{2i-1, i}, X^{(1,1;0,1)}_{2i, i})$ and $z_{i+1}=\min(X^{(1,1; -1,0)}_{2i+1, i+1}, X^{(1,1; -1,0)}_{2i+2, i+1})$ respectively. Given a left and a right ribbon in $S'_i\cup S'_{i+1}$  let $r_i=\min(X^{v_i}_{2i+1; i}, X^{v'_i}_{2i+2; i})$ and  $r_{i+1}=\min(X^{v_i}_{2i-1; i+1}, X^{v'_i}_{2i; i+1})$ respectively. Finally, let  $N_i$ and $N_{i+1}$ denote  the sum of all neutral  arcs from the sets $[X_{2i+12i+2}]$ and $[{X}_{2i-12i}]$ respectively described in 1(b) in Notations \ref{not:components} except for the standard components $X^{(0,0;0,1)}_{2i-12i}$ and ${X}^{(0,0;-1,0)}_{2i+12i}$ (i.e. those with zero intersection with the crosscaps).  For each $1\le  i\le g-1$ we set  $d_i=\epsilon_i+\epsilon_{i+1}$ where  $\epsilon_i=s_i+z_{i+1}+w_i+N_i$ and $\epsilon_{i+1}=s_{i+1}+z_{i}+w_{i+1}+N_{i+1}$. \end{notation}
 
 \begin{remark}
 
 If $\epsilon_i\neq 0$ then $\epsilon_{i+1}=0$, and that only one of $s_i, z_{i+1}$ and $w_i$ can be different than zero by Remark \ref{rem:poss}. Similarly,  if $\epsilon_{i+1}\neq 0$ then $\epsilon_{i}=0$ and only one of $s_{i+1}, z_{i}$ and $w_{i+1}$ is  nonzero.
 \end{remark}

\begin{lem}[Equality for a clover of type I]\label{lem:clover1} 
Given a clover of type I in $S'_i\cup S'_{i+1}$ we have
\begin{eqnarray}\label{eq:clover1}
2\beta_{n+i}+C_i=\max(\gamma_{2i-1}+\gamma_{2i+2}, \gamma_{2i}+\gamma_{2i+1})+2d_i.
\end{eqnarray}
 \noindent where $d_i$ is  as given  in Notation \ref{not:clover1}.
\end{lem}
  \begin{proof}
  Let $[X_{2i-12i}]=\{X^{v_i}_{2i-12i}\}$,  $[X_{2i+12i+2}]=\{X^{v_i}_{2i+12i+2}\}$ and  $ [X_i]$ be the sets of neutral components of $L_i$ described in item 1(b) and item 3.  of Notation \ref{not:components} respectively.  For $2i-1 \le l\le  m \le  2i+2$ with $(l, m)\neq (2i-1,2i), (2i+1,2i+2)$ and $k\in\{i, i+1\}$ write $[X_{l}]=\{X^{v_i}_{l, m}, \, X^{v_i}_{l; k}\}$ for  the set of components described in item 1(a), item 1(c) and item 2. of Notation \ref{not:components}, and  $[P_i]$ for set of simple closed curves supported in $S'_i\cup S'_{i+1}$.  By Remark \ref{rem:poss} we have the following cases:
\begin{itemize}
\item [Case 1:]  $L_i$ is either positive or negative  or neutral.  
\vspace{0.1 cm}

\begin{itemize}
\item[a)] $[{X}_{2i-1}]\cup [{X}_{2i+2}] \neq \emptyset$ and $[{X}_{2i}]\cup [{X}_{2i+1}]=\emptyset$; 
\item[b)] $[{X}_{2i-1}]\cup [{X}_{2i+2}] = \emptyset$ and $[{X}_{2i}]\cup [{X}_{2i+1}]\neq \emptyset$;
\item[c)] $[{X}_{2i-1}]\cup [{X}_{2i+2}] = \emptyset$ and $[{X}_{2i}]\cup [{X}_{2i+1}]= \emptyset$;
\end{itemize}

\item [Case 2:] $L_i$ is mixed that is $L_i$ contains  either scissors or an anchor or a ribbon. There are 2 subcases 
\vspace{0.1 cm}

\begin{itemize}
\item[a)] $[{X}_{2i+1}]\cup [{X}_{2i+2}] \neq \emptyset$ and $[{X}_{2i-1}]\cup [{X}_{2i}]=\emptyset$; 
\item[b)] $[{X}_{2i-1}]\cup [{X}_{2i+2}] = \emptyset$ and $[{X}_{2i}]\cup [{X}_{2i+1}]\neq \emptyset$;
\end{itemize}

 %That is, any other element of such a set may belong to $\{X_{2i-12i+1}, \,X_{2i2i+2}\}$, $[P_i]$, $[X_i]$, $[{X}_{2i-12i-1}]$, $[{X}_{2i+22i+2}]$, $[X_{2i-12i}]$ or $[X_{2i+12i+2}]$ 
  
  		 For Case 1 we have $s_k=z_{k}=r_k=0$ ($k=i, i+1$)  since $L_i$ is not mixed.    Suppose first that Case 1a) holds. Then,  $L_i$ is negative by Remark \ref{rem:posneg}  and Definition \ref{def:compatible} (see for example Figure \ref{fig:compatiblesx}(i)); and we have $\gamma_{2i-1}+\gamma_{2i+2}> \gamma_{2i}+\gamma_{2i+1}$.   First assume that $L_i$ is standard.   It is easy to check that each negative standard component in $[{X}_{2i-1}]\cup [{X}_{2i+2}]$, each curve in $[P_i]$  and each neutral component except for $X^{(0,1;0,1)}_{2i-12i}$ and  $X^{(0,1;0,1)}_{2i+12i+2}$  satisfies  \begin{align}\label{eq:proof1}2\beta_{n+i}+C_i=\gamma_{2i-1}+\gamma_{2i+2}=\max(\gamma_{2i-1}+\gamma_{2i+2}, \gamma_{2i}+\gamma_{2i+1}).\end{align}
 Suppose that  $X^{(0,1;0,1)}_{2i-12i} \neq 0$ and hence $X^{(0,1;0,1)}_{2i+12i+2}= 0$ (note that elements of $[X_{2i-12i}]$ and  $[X_{2i+12i+2}]$  are not compatible).                                    We check that $X^{(0,1;0,1)}_{2i-12i}$ satisfies  $\gamma_{2i-1}+\gamma_{2i+2}=\gamma_{2i}+\gamma_{2i+1}$, and  $$2\beta_{n+i}+C_i-(\gamma_{2i-1}+\gamma_{2i+2})=2X^{(0,1;0,1)}_{2i-12i}.$$ That is, 
    \begin{align}\label{eq:proof2}2\beta_{n+i}+C_i=\max(\gamma_{2i-1}+\gamma_{2i+2}, \gamma_{2i}+\gamma_{2i+1})+2 X^{(0,1;0,1)}_{2i-12i}.\end{align}
    
Similar argument holds for the case $X^{(0,1;0,1)}_{2i-12i}=0$ and $X^{(1,0;-1,0)}_{2i+12i+2}\neq 0$, and we get    
 \begin{align}\label{eq:proof2b}2\beta_{n+i}+C_i=\max(\gamma_{2i-1}+\gamma_{2i+2}, \gamma_{2i}+\gamma_{2i+1})+2 X^{(1,0;-1,0)}_{2i+12i+2}\end{align}

 	 Equalities (\ref{eq:proof1}),  (\ref{eq:proof2}) and   (\ref{eq:proof2b})  are also satisfied for  corresponding twisted components: given a standard component $X_{l, m}$ of $L_i$ with intersection numbers  $\gamma_{j}~(2i-1\leq j\leq 2i+2)$, $\beta_{n+i}$ and $\cC_i$ the twisted component $X^{v_i}_{l, m}$  has the same  number of intersections on $\cC_i$ and increases the number of intersections on each $\gamma_{j}~(2i-1\leq j\leq 2i+2)$ and $2\beta_{n+i}$ (and hence $2\beta_{n+i}+C_i$) by the same amount determined by its signature vector (Figure \ref{fig:compatiblesx}(ii)). Therefore  equality (\ref{eq:proof1}) is also satisfied for each $X^{v_i}_{l, m}$. Similar argument holds for other twisted components. Then we have $d_i=\epsilon_i=N_i$, and

    $$2\beta_{n+i}+C_i=\max(\gamma_{2i-1}+\gamma_{2i+2}, \gamma_{2i}+\gamma_{2i+1})+2d_i$$ 
 \noindent as required.   Case 1b)  follows by symmetry.  For Case 1c), we note that since $L_i$ consists of only neutral components and each neutral component satisfies $\gamma_{2i-1}+\gamma_{2i+2}= \gamma_{2i}+\gamma_{2i+1}$ except for $X^{(0,1;0,1)}_{2i-12i}$ and  $X^{(0,1;0,1)}_{2i+12i+2}$ as shown above, we get either $\epsilon_i=N_i\neq 0$ or  $\epsilon_{i+1}=N_{i+1}\neq 0$ since $X^{(0,1;0,1)}_{2i-12i}$ and  $X^{(0,1;0,1)}_{2i+12i+2}$ are not compatible. Therefore, $d_i=\epsilon_i+\epsilon_{i+1}$, and
 
 \begin{align}2\beta_{n+i}+C_i=\max(\gamma_{2i-1}+\gamma_{2i+2}, \gamma_{2i}+\gamma_{2i+1})+2 d_i.\end{align}
 \noindent as required.
    
     \begin{figure}[h!]
\begin{center}
\labellist

 \pinlabel {${\scriptstyle{\gamma_{2i}}}$} [ ] at  85 -3
  \pinlabel {${\scriptstyle{\gamma_{2i+2}}}$} [ ] at  185 -3
  
   \pinlabel {${\scriptstyle{\gamma_{2i}}}$} [ ] at  375 -3
  \pinlabel {${\scriptstyle{\gamma_{2i+2}}}$} [ ] at  490 -3

 \pinlabel {${\scriptstyle{\gamma_{2i-1}}}$} [ ] at  80 290
  \pinlabel {${\scriptstyle{\gamma_{2i+1}}}$} [ ] at  185 290

 \pinlabel {${\scriptstyle{\gamma_{2i-1}}}$} [ ] at  370 290
  \pinlabel {${\scriptstyle{\gamma_{2i+1}}}$} [ ] at  490 290

   \pinlabel {${\scriptstyle{\gamma_{2i}}}$} [ ] at  680 -3
  \pinlabel {${\scriptstyle{\gamma_{2i+2}}}$} [ ] at  795 -3

   \pinlabel {${\scriptstyle{\gamma_{2i-1}}}$} [ ] at  675 290
  \pinlabel {${\scriptstyle{\gamma_{2i+1}}}$} [ ] at  790 290

%   \pinlabel {${\scriptstyle{\beta_{n+i-1}}}$} [ ] at  10 290
  \pinlabel {${\scriptstyle{\beta_{n+i}}}$} [ ] at  130 290
  
  %  \pinlabel {${\scriptstyle{\beta_{n+i+1}}}$} [ ] at  250 290

 %  \pinlabel {${\scriptstyle{\beta_{n+i-1}}}$} [ ] at  330 290
  \pinlabel {${\scriptstyle{\beta_{n+i}}}$} [ ] at  430 290
    \pinlabel {${\scriptstyle{\beta_{n+i}}}$} [ ] at  740 290

%  \pinlabel {${\scriptstyle{\beta_{n+i+1}}}$} [ ] at  570 290
     \pinlabel {(i)} [ ] at  130 -20
    \pinlabel {(ii)} [ ] at  435 -20
      \pinlabel {(iii)} [ ] at  750 -20

 \small\hair 2pt
\endlabellist  

  \includegraphics[scale=0.35]{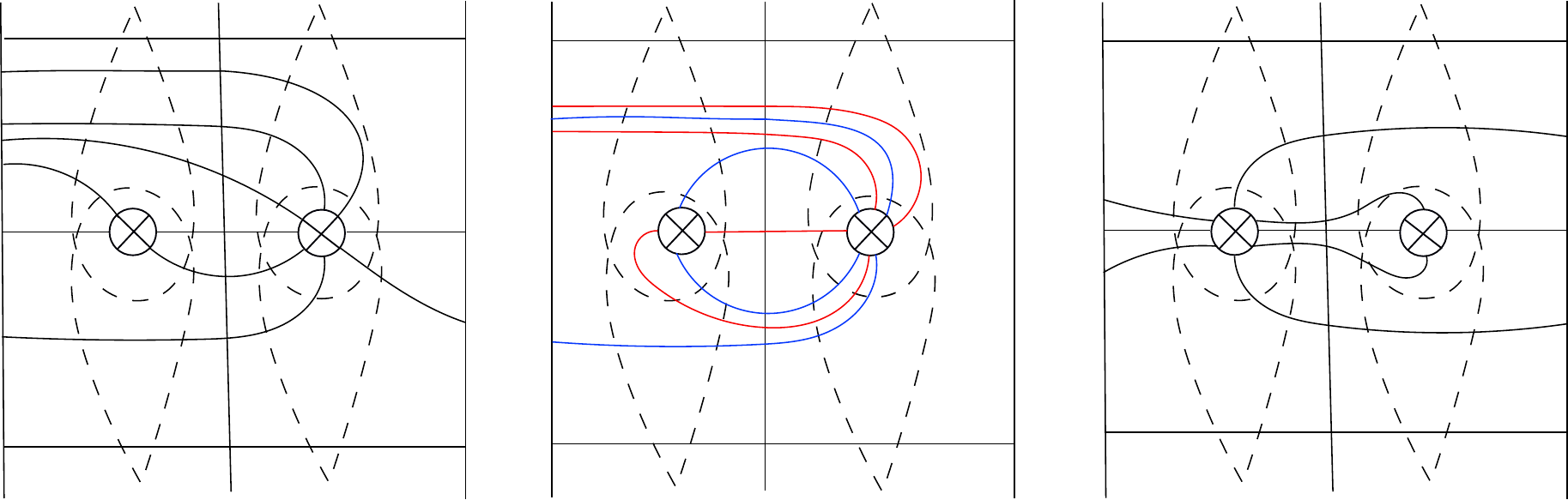}
\caption{ (i) Case 1a) where $L_i$ is standard;   (ii) Case 1a) where $L_i$ is twisted and (iii) Case 1c)} \label{fig:compatiblesx}
\end{center}
\end{figure}

\end{itemize}

Case 2a) is divided into 3 subcases. Either $L_i$ contains scissors at $\otimes_i$ or a left anchor or a left ribbon  by Definition \ref{def:scissorsanchors} and Remark \ref{rem:poss}. Assume that  $L_i$ contains scissors at $\otimes_i$. Then  we have $X_{2i+1; i}\neq 0$, $X_{2i+2; i}\neq 0$ and hence $d_i=\epsilon_i=s_i=\min(X_{2i+1; \,i}, X_{2i+2; \,i})+N_i$ where $N_i= X^{(1,0;-1,0)}_{2i+12i+2}$ (twisted neutral components are not compatible with scissors).  We obtain

\begin{equation}
\label{eq:scissorsaa}
  \begin{aligned}
 2\beta_{n+i}+C_i-2X_{2i+1; i}-2X^{(1,0;-1,0)}_{2i+12i+2}&=\gamma_{2i-1}+\gamma_{2i+2}
     \end{aligned}
\end{equation}

\begin{equation}
\label{eq:scissorsbb}
  \begin{aligned}
 2\beta_{n+i}+C_i-2X_{2i+2; i}-2X^{(1,0;-1,0)}_{2i+12i+2}&=\gamma_{2i}+\gamma_{2i+1}
     \end{aligned}
\end{equation}
 as shown in Figure \ref{fig:compatiblesy}(i).  Observe that $\gamma_{2i-1}+\gamma_{2i+2}\geq \gamma_{2i}+\gamma_{2i+1}$ if and only if $X_{2i+2; i}\geq X_{2i+1; i }$. Then from equation (\ref{eq:scissorsaa}) and equation (\ref{eq:scissorsbb})  we get

\begin{align*} 2\beta_{n+i}+C_i&=\max(\gamma_{2i-1}+\gamma_{2i+2}, \gamma_{2i}+\gamma_{2i+1})+2\min(X_{2i+2; i}, X_{2i+1; i})+2X^{(1,0;-1,0)}_{2i+12i+2}\\
&=\max(\gamma_{2i-1}+\gamma_{2i+2}, \gamma_{2i}+\gamma_{2i+1})+2s_i+2X^{(1,0;-1,0)}_{2i+12i+2}.
\end{align*}

\noindent as required. Similarly, if $L_i$ contains a left anchor it contains no scissors, no ribbons and no right anchor.  We have $d_i=z_{i+1}+N_i$ and equation (\ref{eq:clover1}) is verified analogously. The case when there is a left ribbon is proved similarly.  Case 2b) is also divided into 3 subcases: Either  $L_i$ contains scissors at $\otimes_{i+1}$ or a right anchor or a right ribbon by Definition \ref{def:scissorsanchors} and Remark \ref{rem:poss}.  Then Case 2b)  follows immediately from Case 2a) by symmetry.  \end{proof}

%$X_{2i-1; i+1}\neq 0 $ and $X_{2i; i+1}\neq 0$ is similar  and we get 
%\begin{align*}2\beta_{n+i}+C_i=\max(\gamma_{2i-1}+\gamma_{2i+2}, \gamma_{2i}+\gamma_{2i+1})+2s_{i+1}+2X^{(0,1;0,1)}_{2i-12i}
%\end{align*}
%
%
%That is,  $2\beta_{n+i}+C_i=\max(\gamma_{2i-1}+\gamma_{2i+2}, \gamma_{2i}+\gamma_{2i+1})+2d_i$ as required.  The case for when there are ribbons is proved analogously.
%
%
%
 \begin{figure}[h!]
\begin{center}
\labellist

 \pinlabel {${\scriptstyle{\gamma_{2i}}}$} [ ] at  70 -3
  \pinlabel {${\scriptstyle{\gamma_{2i+2}}}$} [ ] at  180 -3
  
   \pinlabel {${\scriptstyle{\gamma_{2i}}}$} [ ] at  370 -5
  \pinlabel {${\scriptstyle{\gamma_{2i+2}}}$} [ ] at  490 -5

 \pinlabel {${\scriptstyle{\gamma_{2i-1}}}$} [ ] at  60 290
  \pinlabel {${\scriptstyle{\gamma_{2i+1}}}$} [ ] at  170 290

 \pinlabel {${\scriptstyle{\gamma_{2i-1}}}$} [ ] at  360 290
  \pinlabel {${\scriptstyle{\gamma_{2i+1}}}$} [ ] at  480 290

   %\pinlabel {${\scriptstyle{\beta_{n+i-1}}}$} [ ] at  10 290
  \pinlabel {${\scriptstyle{\beta_{n+i}}}$} [ ] at  115 290
  
  %  \pinlabel {${\scriptstyle{\beta_{n+i+1}}}$} [ ] at  250 290

   %\pinlabel {${\scriptstyle{\beta_{n+i-1}}}$} [ ] at  330 290
  \pinlabel {${\scriptstyle{\beta_{n+i}}}$} [ ] at  420 290
 % \pinlabel {${\scriptstyle{\beta_{n+i+1}}}$} [ ] at  590 290
     \pinlabel {(i)} [ ] at  120 -15
    \pinlabel {(ii)} [ ] at  425 -15

 \small\hair 2pt
\endlabellist  

  \includegraphics[scale=0.4]{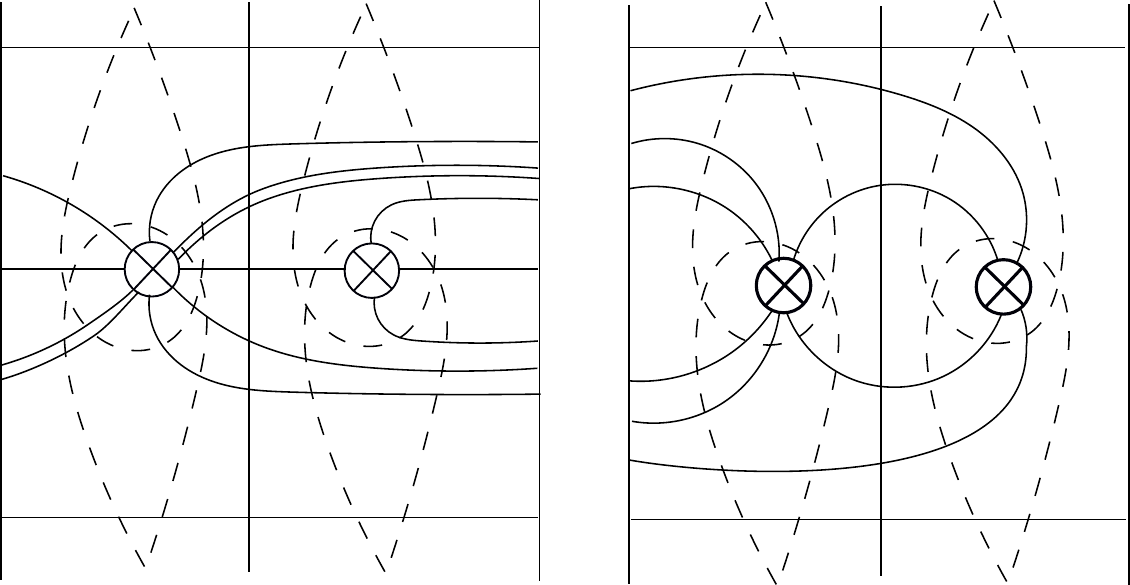}
\caption{Examples for $L_i$ for Case 2a) and Case 2b)} \label{fig:compatiblesy}
\end{center}
\end{figure}

The proof of Lemma \ref{lem:clover1} shows that each component of $L_i$ except for certain arcs and arc systems satisfies
$$2\beta_{n+i}+C_i=\max(\gamma_{2i-1}+\gamma_{2i+2}, \gamma_{2i}+\gamma_{2i+1}).$$ We shall call these arcs and arc systems {\it exceptional arcs and arc systems with respect to a clover type I}; and $d_i$ the {\it exceptional parameter for a clover of type I}. 
\begin{defn}\label{def:extype1}
Each neutral arc in $[X_{2i+12i+2}]\setminus\{X^{(0,0;-1,0)}_{2i+12i+2}\}$ and $[X_{2i-12i}]\setminus \{X^{(0,0;0,1)}_{2i-12i}\}$  is called an {\it exceptional arc with respect to a clover type I at $\otimes_i$ and $\otimes_{i+1}$} respectively (for example in Figure \ref{fig:standardarcs} the arc labeled 2 is an exceptional arc at $\otimes_i$ while the arc labeled 9 is not); scissors, ribbons and anchors  in $S'_i\cup S'_{i+1}$ are called {\it exceptional arc systems with respect to a clover of type I}.

\end{defn}
 
 Let $\Delta_{k}(A)$ and $\Delta_{k}(B)$~$(k=i, i+1)$ be as described in Notation  \ref{not:deltas}. The proof of Lemma \ref{cor:except} follows immediately from the definition of exceptional arcs and arc systems.
 
\begin{lem}\label{cor:except}
 $\Delta_{i+1}(A)\neq 0$ and $\Delta_{i+1}(B)\neq 0$ if and only if there exists at least one of the following arc or arc systems in $S'_i\cup S'_{i+1}$:  an arc from the set $[X_{2i+12i+2}]$, scissors at $\otimes_i$, a left anchor or a left ribbon.  Similarly, $\Delta_{i}(A)\neq 0$ and $\Delta_{i}(B)\neq 0$ if and only if at least one of the following exists in $S'_i\cup S'_{i+1}$: an arc from the set $[X_{2i-12i}]$, scissors at $\otimes_{i+1}$, a right anchor or a right ribbon. 
\end{lem}

Then, we can compute the exceptional parameter $d_i=\epsilon_i+\epsilon_{i+1}$ in terms of generalized Dynnikov coordinates.  We first give the following preliminary lemma.

\begin{lem}\label{lem:lambda2}
Let $X_{2i+12i+2}^{(0,0;-1,0)}$ and $X_{2i-12i}^{(0,0;0,1)}$  be as given in Notation \ref{not:components}. Then,
 \begin{align} X_{2i+12i+2}^{(0,0;-1,0)}&=\min(\Delta_{i+1}(A), \, \Delta_{i+1}(B),\,\lambda^-_{{i}})\\
 X_{2i-12i}^{(0,0;0,1)}&=\min(\Delta_{i}(A), \, \Delta_{i}(B),\,\lambda^+_{{i+1}})
\end{align}
\end{lem}

\begin{proof}
The proof follows from Remark \ref{rem:deltas} and Figure \ref{fig:deltas}. \end{proof}

\begin{figure}[h!]
\begin{center}
\labellist
 \pinlabel {${\scriptstyle{A'_{i,i+1}}}$} [ ] at  500 235
  \pinlabel {${\scriptstyle{\Delta_i(A)}}$} [ ] at  370 200
  
   \pinlabel {${\scriptstyle{B'_{i,i+1}}}$} [ ] at  500 40
  \pinlabel {${\scriptstyle{\Delta_i(B)}}$} [ ] at  370 80
  % \pinlabel {${\scriptstyle{\lambda_{c_{i+1}}}}$} [ ] at  472 170
      \pinlabel {${\scriptstyle{\lambda^+_{{i+1}}}}$} [ ] at  535 170

     % \pinlabel {${\scriptstyle{\lambda_{c_i}}}$} [ ] at  95 150
      \pinlabel {${\scriptstyle{\lambda^-_{i}}}$} [ ] at 45 150

  \pinlabel {${\scriptstyle{\Delta_{i+1}(A)}}$} [ ] at  190 200
 \pinlabel {${\scriptstyle{A'_{i, i+1}}}$} [ ] at  80 235

  \pinlabel {${\scriptstyle{\Delta_{i+1}(B)}}$} [ ] at  185 80
 \pinlabel {${\scriptstyle{B'_{i, i+1}}}$} [ ] at  80 40
 \pinlabel {${\scriptstyle{\beta_{n+i-1}}}$} [ ] at  15 260
\pinlabel {${\scriptstyle{\beta_{n+i}}}$} [ ] at  120 260
\pinlabel {${\scriptstyle{\beta_{n+i+1}}}$} [ ] at  255 260
 \pinlabel {${\scriptstyle{\beta_{n+i-1}}}$} [ ] at  315 260
\pinlabel {${\scriptstyle{\beta_{n+i}}}$} [ ] at  450 260
\pinlabel {${\scriptstyle{\beta_{n+i+1}}}$} [ ] at  575 260
\pinlabel {(a)} [ ] at  130 -13

\pinlabel {(b)} [ ] at  450 -13

 \small\hair 2pt
\endlabellist  

  \includegraphics[scale=0.47]{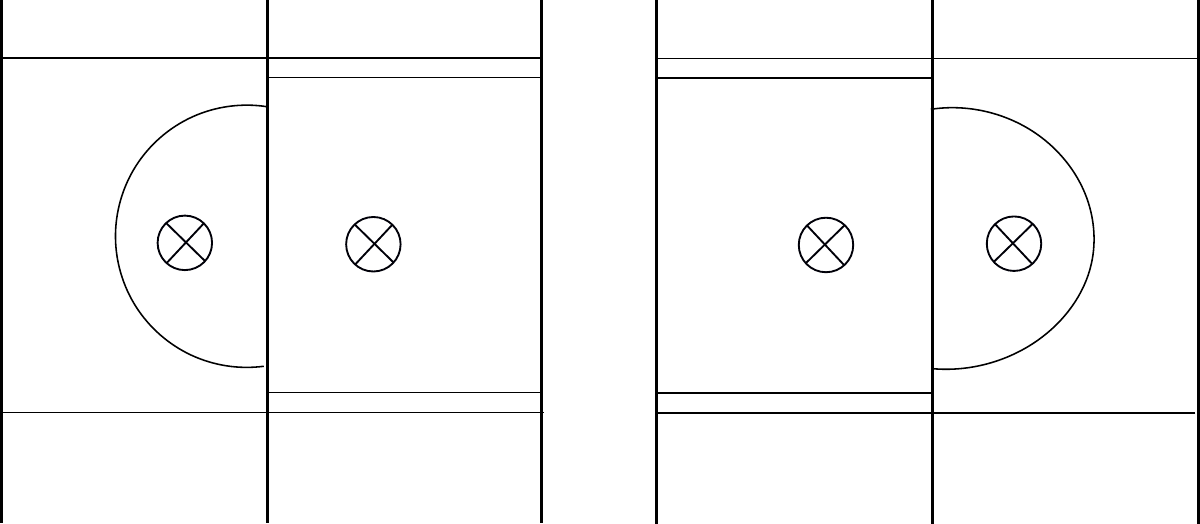}
\caption{Computation for $X_{2i+12i+2}^{(0,0; -1,0)}$ and  $X_{2i-12i}^{(0,0; 0,1)}$}\label{fig:deltas}
\end{center}
\end{figure}

\begin{lem}\label{lem:di}

  Let $\epsilon_k$~$(k=i, i+1$) be as described in  Notation \ref{not:clover1}. Then $d_i=\epsilon_i+\epsilon_{i+1}$ where

\begin{align}\label{eq:di}
\epsilon_i&=\min(\Delta_{i+1}(A), \Delta_{i+1}(B))-X_{2i+12i+2}^{(0,0;-1,0)}\\
\epsilon_{i+1}&=\min(\Delta_{i}(A), \Delta_{i}(B))-X_{2i-12i}^{(0,0;0,1)}
\end{align}
and $X_{2i-12i}^{(0,0;0,1)}$ and $X_{2i+12i+2}^{(0,0;-1,0)}$ are as given in Lemma \ref{lem:lambda2}.
\end{lem}
\begin{proof}

We compute $\epsilon_i$. $\epsilon_{i+1}$  is computed analogously reflecting in the arc $\beta_{n+i}$. If at least one of $\Delta_{i+1}(A)$ and $ \Delta_{i+1}(B)$ is equal to zero then $\epsilon_i=0$ by Lemma \ref{cor:except}. So suppose that $\Delta_{i+1}(A)\neq 0$ and $ \Delta_{i+1}(B)\neq 0$.  

 \begin{itemize}

\item[Case 1.] There exists no exceptional arc system with respect to a clover type I  in which case there must be exceptional arcs  from the set $[X_{2i+12i+2}]$.

% and hence $d_i=\epsilon_i=N_i$ where $N_i$ denotes the number of exceptional arcs of $L\cap(S'_i\cup S'_{i+1})$. 

\item[Case 2.]  There exists an exceptional arc system with respect to a clover type I which can be either scissors at $\otimes_i$ or a left anchor or a left ribbon possibly with compatible exceptional arcs  from the set $[X_{2i+12i+2}]$. 

\end{itemize} 

We have two subcases in Case 1: (a) $\Delta_{i+1}(A)\leq \Delta_{i+1}(B)$ and (b) $\Delta_{i+1}(A)\geq \Delta_{i+1}(B)$. Assume that we are in Case 1(a). Then there exist negative but no positive arcs in $S'_i\cup S'_{i+1}$  since only negative components increase the difference $\Delta_{i+1}(B)-\Delta_{i+1}(A)$ by 
Remark \ref{rem:posneg} and that there is no exceptional arc system with respect to a clover type I in  $S'_i\cup S'_{i+1}$  by assumption. Since each element of $[X_{2i+12i+2}]$ increases $\Delta_{i+1}(A)$ and $\Delta_{i+1}(B)$ by $1$ and $N_i$ of those are exceptional we have

\begin{align*}
\Delta_{i+1}(A)&=N_i+X_{2i+12i+2}^{(0,0;-1,0)}\leq \Delta_{i+1}(B)
\end{align*}

\noindent (recall that $X_{2i+12i+2}^{(0,0;-1,0)}$ is the only element of $[X_{2i+12i+2}]$ which is not exceptional). Since $d_i=\epsilon_i=N_i$ we get $d_i=\min(\Delta_{i+1}(A), \Delta_{i+1}(B))-X_{2i+12i+2}^{(0,0;-1,0)}$ as required. Case 1(b) is proved similarly. Now assume that  we are in Case 2.   Assume that there exists scissors at $\otimes_i$ (Figure \ref{fig:compatiblesy}). Since the only element of $[X_{2i+12i+2}]$  which is compatible with the scissors is the standard exceptional arc $X_{2i+12i+2}^{(1,0;-1,0)}$ we have $N_i=X_{2i+12i+2}^{(1,0;-1,0)}$, and hence $d_i=\epsilon_i=s_i+N_i$. By Remark \ref{rem:poss} every other arc compatible with the scissors is neutral and has no affect on $\Delta_{i+1}(A)$ and $\Delta_{i+1}(B)$. Therefore,

\begin{align*}
\Delta_{i+1}(A)&=X_{2i+1; i}+X_{2i+12i+2}^{(1,0;-1,0)}\\
\Delta_{i+1}(B)&=X_{2i+2; i}+X_{2i+12i+2}^{(1,0;-1,0)}
\end{align*}
Since $s_i=\min(X_{2i+1; i},X_{2i+2; i})$,  $N_i=X_{2i+12i+2}^{(1,0;-1,0)}$ and $X_{2i+12i+2}^{(0,0;-1,0)}=0$ we obtain  $$d_i=\min(\Delta_{i+1}(A), \Delta_{i+1}(B))=\min(\Delta_{i+1}(A), \Delta_{i+1}(B))+X_{2i+12i+2}^{(0,0;-1,0)}$$ 
\noindent as required. In the cases when there is a left anchor and a left ribbon we note that there may exist both twisted and standard exceptional arcs in  $[X_{2i+12i+2}]$. Also,  $X_{2i+12i+2}^{(0,0;-1,0)}= 0$ if there exists a left ribbon and $X_{2i+12i+2}^{(0,0;-1,0)}\geq 0$ if there exists a left anchor. Again, since each positive (resp. negative) arm of a left anchor or a left ribbon increases $\Delta_{i+1}(A)$ (resp. $\Delta_{i+1}(B)$) by $1$,  and only neutral arcs are compatible with exceptional arc and arc systems the proof follows similarly.   \end{proof}

\begin{figure}[h!]
\begin{center}
\labellist
 \pinlabel {${\scriptstyle{\beta_{n+i-1}}}$} [ ] at  4 288
  \pinlabel {${\scriptstyle{\beta_{n+i}}}$} [ ] at  118 288

 \pinlabel {${\scriptstyle{\beta_{n+i+1}}}$} [ ] at  260 288
 
 \pinlabel {${\scriptstyle{\beta_{n+i-1}}}$} [ ] at  395 288
  \pinlabel {${\scriptstyle{\beta_{n+i}}}$} [ ] at 510 288

 \pinlabel {${\scriptstyle{\beta_{n+i+1}}}$} [ ] at  640 288
 
   \pinlabel {(a)} [ ] at  120 -10
    \pinlabel {(b)} [ ] at  505 -10

 \small\hair 2pt
\endlabellist  

  \includegraphics[scale=0.3]{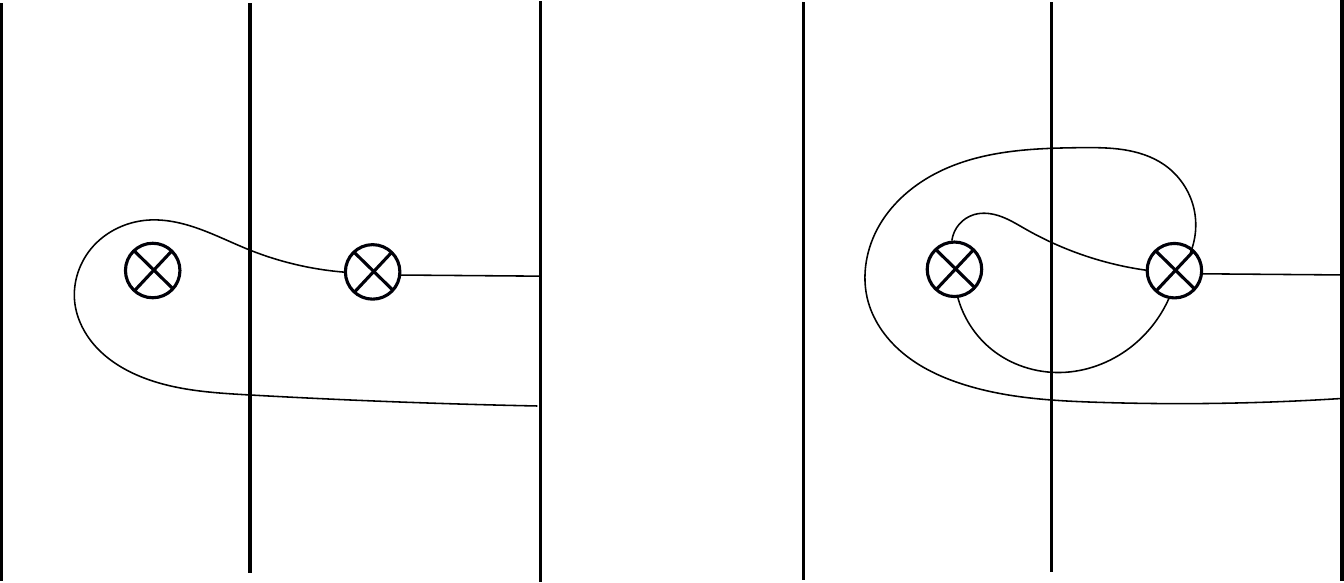}
\caption{$u^{-1}_i$ images of exceptional arcs of type I give exceptional arcs of type II}\label{fig:htexep} 
\end{center}
\end{figure}

\begin{figure}[h!]
\begin{center}
\labellist
 \pinlabel {${\scriptstyle{\beta_{n+i-1}}}$} [ ] at  4 288
  \pinlabel {${\scriptstyle{\beta_{n+i}}}$} [ ] at  118 288

 \pinlabel {${\scriptstyle{\beta_{n+i+1}}}$} [ ] at  260 288
 
 \pinlabel {${\scriptstyle{\beta_{n+i-1}}}$} [ ] at  350 288
  \pinlabel {${\scriptstyle{\beta_{n+i}}}$} [ ] at 465 288

 \pinlabel {${\scriptstyle{\beta_{n+i+1}}}$} [ ] at  580 288
 
  \pinlabel {${\scriptstyle{\beta_{n+i-1}}}$} [ ] at  680 288
  \pinlabel {${\scriptstyle{\beta_{n+i}}}$} [ ] at 780 288

 \pinlabel {${\scriptstyle{\beta_{n+i+1}}}$} [ ] at  900 288
  \pinlabel {(a)} [ ] at  120 -13
    \pinlabel {(b)} [ ] at  465 -13
      \pinlabel {(c)} [ ] at  780 -13

 \small\hair 2pt
\endlabellist  

  \includegraphics[scale=0.3]{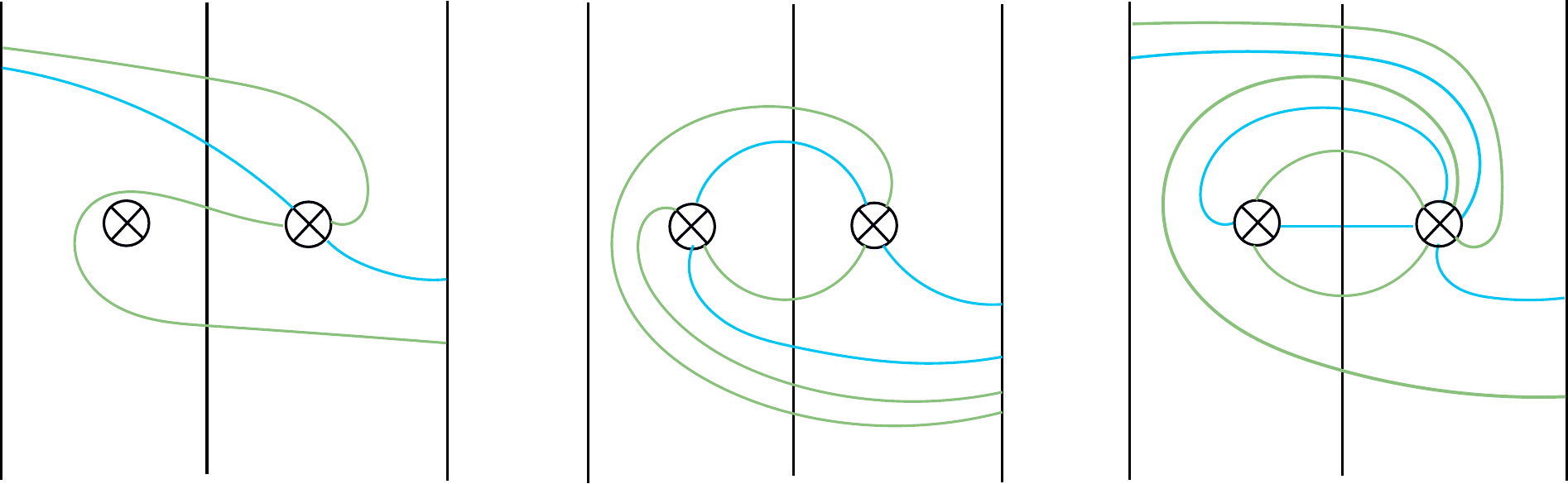}
\caption{ $u^{-1}_i$ images of scissors, anchors and ribbons are exceptional arc systems with respect to a clover of type II} \label{fig:halftwisted}
\end{center}
\end{figure}

  Taking the $u_i^{-1}$ images of exceptional arcs with respect to a clover type I in $S'_i\cup S'_{i+1}$ we obtain elements of $[X_{2i+2; i+1}]$ and $[X_{2i-1; i}]$ which are called {\it exceptional arcs with respect to a clover of type II} in $S'_i\cup S'_{i+1}$. Similarly, taking the $u_i^{-1}$ images of scissors, anchors and ribbons in $S'_i\cup S'_{i+1}$ we get negatively half twisted scissors, anchors and ribbons which are  {\it exceptional arc systems with respect to a clover type II}  (see Figure \ref{fig:htexep} and Figure \ref{fig:halftwisted} for some examples). This leads us to the equality in Lemma \ref{lem:clover2}. But first we introduce some notation for the parameters associated with exceptional arc and arc systems with respect to a clover type II.

 \begin{notation}\label{not:cloverII}
Let $e_i=\epsilon'_i+\epsilon'_{i+1}$ where $\epsilon'_i=s'_i+z'_i+r_i'+N'_i$ where $N'_i$ and $N'_{i+1}$  denote the number of exceptional arcs of type II which are from the sets $[X_{2i-1; i}]$ and $[X_{2i+2; i}]$ respectively, and 
\begin{align*}
s'_i&=\min\big(u^{-1}_i(X_{2i+1; i}), u_i^{-1}(X_{2i+2; i})\big), \hspace{1.1 cm}   s'_{i+1}=\min\big(u^{-1}_i(X_{2i-1; i+1}), u^{-1}_i(X_{2i; i+1})\big),\\  
z'_{i}&=\min\big(u^{-1}_i(X^{(1,1; -1,0)}_{2i+1, i+1}), u^{-1}_i(X^{(1,1; -1,0)}_{2i+2, i+1})\big), \hspace{0.1 cm}  z'_{i+1}=\min\big(u_i^{-1}(X^{(1,1; 0,1)}_{2i-1, i}), u_i^{-1}(X^{(1,1;0,1)}_{2i, i})\big),\\
 r'_i&=\min\big(u_i^{-1}(X^{v_i}_{2i+1, i}), u_i^{-1}(X^{v'_i}_{2i+2, i})\big),   \hspace{1.1 cm}   r'_{i+1}=\min\big(u_i^{-1}(X^{v_i}_{2i-1, i+1}), u_i^{-1}(X^{v'_i}_{2i, i+1})\big).\end{align*}
$e_i$ is  called  the exceptional parameter for a clover of type II.
\end{notation}

\begin{lem}[Equality for a clover of type II]\label{lem:clover2} 
Given a clover of type II we have
 \begin{eqnarray}\label{eq:clover2}
2\beta'_{n+i}+C_i=\max(\gamma_{2i+1}+\gamma_{2i}, \gamma'_{2i+1}+\gamma'_{2i})+2e_i.
\end{eqnarray}
where the exceptional parameter $e_i$ is as given in Notation \ref{not:cloverII}. 
\end{lem}
We require the following parameters to compute $e_i$ and the other exceptional parameters in terms of generalized Dynnikov coordinates.

 \begin{notation}\label{not:twist}    Let $X$ be a component  of $L_i$ with $\widehat{v_i}=(\widehat{c_i}, \widehat{c}_{i+1}; \widehat{b}_{n+i-1}, \widehat{b}_{n+i})$.  Write
 \begin{align*}\widehat{\chi}_{i}(A)&=(\widehat{b}_{n+i})^{+}-\widehat{\Delta}_{i}(A), \quad \widehat{\chi}_{i+1}(A)=(-\widehat{b}_{n+i-1})^{+}-\widehat{\Delta}_{i+1}(A)\\
   \widehat{\chi}_{i}(B)&=(\widehat{b}_{n+i})^{+}-\widehat{\Delta}_{i}(B),\quad \widehat{\chi}_{i+1}(B)=(-\widehat{b}_{n+i-1})^{+}-\widehat{\Delta}_{i+1}(B)
 \end{align*} \end{notation}

  \begin{defn}\label{rem:together}
   We define $\chi(X)$ as follows: $$
\chi(X)=\begin{cases}
\min( \widehat{\chi}_{i}(B), \widehat{\chi}_{i+1}(A))& \text{if } X \text{ is positive} \\
\min( \widehat{\chi}_{i}(A), \widehat{\chi}_{i+1}(B)), & \text{if }X\text{ is negative }\\
\min( \widehat{\chi}_{k}(A), \widehat{\chi}_{k}(B)), & \text{if }X\text{ is neutral}
\end{cases}
$$We note that  if $X$ is neutral $\widehat{\chi}_{k}(A)=\widehat{\chi}_{k}(B)$ for each $k=i, i+1$. %$\widehat{\chi}_{i}(A), \widehat{\chi}_{i}(B), \widehat{\chi}_{i+1}(A)$ and $\widehat{\chi}_{i+1}(B)$  %Also note that   $\chi(X)=0$ if $X$ is a neutral component since such components are not twisted. 	
  \end{defn}
Geometrically, $\chi(X)$  gives information about the amount of twist of $X$, and reveals whether  $X$ positive or negative or neutral by Remark \ref{rem:posneg}. Observe that a standard component $X$ of $L_i$ either has $\chi(X)=0$ or $\chi(X)=-1$. The possibilities for the latter case is given in Remark \ref{rem:twarc}.

\begin{remark}\label{rem:twarc}Let $X$ be standard. When $X$ is   negative $\chi(X)=-1$ if and only if $$X \in \{X^{(0,0;0,1)}_{2i-12i-1},\, X^{(0,1;0,1)}_{2i-12i-1},\, X^{(0,0; -1,0)}_{2i+22i+2},\, X^{(1,0;-1,0)}_{2i+22i+2},\, X_{2i-12i+2}\}$$ See for instance $l_2=X^{(0,0;-1,0)}_{2i+22i+2},\, l_3=X^{(1,0;-1,0)}_{2i+22i+2},\, l^{\star}_3=X_{2i-12i+2}$ in Figure \ref{fig:illustration2}. Similarly, when  $X$ is   positive    $\chi(X)=-1$ if and only if   $$X\in \{X^{(0,0; -1,0)}_{2i+12i+1},\, X^{(1,0;-1,0)}_{2i+12i+1},\, X^{(0,0; 0,1)}_{2i2i},\, X^{(0,1;0,1)}_{2i2i},\, X_{2i2i+1}\}.$$   \end{remark}

\begin{remark}\label{rem:notcompatible}
Note that a standard arc cannot be compatible with a highly twisted arc $X$ with $\chi(X)>1$. Furthermore, if an arc $X$ has $\psi_{i+1}\neq 0$ it is either a standard arc or a highly twisted arc from the sets $[X_{2i+2;i+1}]$ or $[X_{2i+1;i+1}]$. Similar argument holds for an arc $X$ with $\psi_{i}\neq 0$.  \end{remark}

\begin{defn}\label{def:twset}
Suppose that $L_i$ is not mixed.  We define $\chi(i)$ for $L_i$ by replacing $X$ with $L_i$ and removing all hats from the symbols given in Notation \ref{not:twist}.     \end{defn}

\begin{notation}\label{not:pis}Let $X_{2i-12i}=X_{2i-12i}^{(0,0;0,1)}$ and $X_{2i+12i+2}=X_{2i+12i+2}^{(0,0;-1,0)}$ be as given in Lemma \ref{lem:lambda2}. For notational simplicity we shall denote  ${\Lambda}_i=\lambda^-_i-X_{2i+12i+2},\, {\Lambda}_{i+1}= \lambda^+_{i+1}-X_{2i-12i}$ and
\begin{align*}
  \bar{\Delta}_{i}(B)&= \Delta_{i}(B)-X_{2i-12i},  \quad   \bar{\Delta}_{i}(A)= \Delta_{i}(A)-X_{2i-12i}\\
 \bar{\Delta}_{i+1}(B)&= \Delta_{i+1}(B)-X_{2i+12i+2}, \quad \bar{\Delta}_{i+1}(A)= \Delta_{i+1}(A)-X_{2i+12i+2}
  \end{align*}
 We also introduce the following components for twisted components of $L_i$:
 \begin{align*}
p_1=2(-b_{n+i-1})^+-\Delta_{i+1}(B), \quad p_2=2(b_{n+i})^+-\Delta_{i}(A)\\
p_3=2(-b_{n+i-1})^+-\Delta_{i+1}(A), \quad p_4=2(b_{n+i})^+-\Delta_{i}(B) \end{align*}
Geometrically, $p_1$ denotes the number of loop components of $L\cap S'_i$ which are not contained in below components of $L\cap S'_{i+1}$ for a twisted component of $L_i$. The interpretation of the other parameters $p_k$ is similar.
\end{notation}
To understand these parameters better let us consider Figure \ref{fig:motivation} again where $L_i$ consists of three components: $l_1=X_{2i+12i+2},\,  l_2=X_{2i+1, 2i+2}$ and 
$ l_3=X^{(1,2; -1,1)}_{2i+2; i+1}$.  Since $l_1$ and $l_2$ are neutral $\chi(l_1)=\chi(l_2)=0$. Also, $l_3$ is a negative twisted component with  $\chi(l_3)=\min(1, 1)=1$ and 
$p_1=3$.

\begin{lem}\label{lem:clover2coordinates}
Let $L_i$ be negative. Then, $e_i=\epsilon'_i+\epsilon'_{i+1}$ where
\begin{equation}\label{eq:eps'}
\epsilon'_{i}=\min({\Lambda}_i, \psi_{i+1},\bar{\Delta}_{i+1}(B), p_1) \quad\text{and}\quad \epsilon'_{i+1}=\min({\Lambda}_{i+1}, \psi_{i}, \bar{\Delta}_{i}(A), p_2)
\end{equation}

\end{lem}

%That is to say, $I=\{X_{2i-12i+1}, X_{2i2i+2}, X_{i+1}\}$. 
\begin{proof}
Here we only prove $\epsilon'_{i}$.  The formula for $\epsilon'_{i+1}$ is obtained similarly. We first note that each exceptional arc in $[X_{2i+2; i+1}]$  and negatively half twisted scissors at $\otimes_{i}$, left anchor and left ribbon increases each $\lambda^-_i, \psi_{i+1}, \Delta_{i+1}(B)$ by $1$ (see Figure \ref{fig:htexep} and  Figure \ref{fig:halftwisted}).    Therefore if at least one of $\lambda^-_i, \psi_{i+1}, \Delta_{i+1}(B)$ equals zero for $L_i$  we get $\epsilon'_i=0$.  For convenience, let us say that a subset of $L_i$ has property $P$ if it satisfies  $\lambda^-_i\neq 0$, $\psi_{i+1}\neq 0$ and  $\Delta_{i+1}(B)\neq 0$. The proof is based on constructing all possible configurations of arcs (i.e. compatible sets) satisfying  property $P$, and verifying that equation (\ref{eq:eps'}) holds for each such collection of arcs. Let $L_i$ have property $P$. The constraint provided by property $P$ and Remark  \ref{not:twist} imply that each element of $L_i$ belongs to one of the two sets described as follows: The first set $I$ is the subset of $L_i$  whose elements affect none of  the values $\lambda^-_i$, $\psi_{i+1}$ and $\Delta_{i+1}(B)$ yet compatible with those satisfying property $P$. The second  set $S$ contains negative components affecting at least one of $\lambda^-_i$, $\psi_{i+1}$, $\Delta_{i+1}(B)$. Furthermore, $S$ is partitioned into two subsets $S_{0}$ and $S_{1}$ such that if $X\in S_0$ then $\chi(X)\leq 0$ and $X$ is one of the following arcs depicted in Figure \ref{fig:illustration2}; and  if $X\in S_1$ then $\chi(X)\geq 1$ and $X$ is compatible with an arc with $\psi_{i+1}\neq 0$. In particular, if  $X\in S_1$ with $\psi_{i+1}=0$ it  has $\chi(X)=1$ and it is one of the arcs $l_{15}, l_{16}, l_{17}, l_{18}$ depicted in Figure \ref{fig:illustration3}; and if $\psi_{i+1}\neq 0$ it is a highly twisted  exceptional arc from the set $[X_{2i+2;i+1}]$ such as $l_{18}$ and $l_{19}$ as depicted in Figure \ref{fig:illustration3}.   Let $\mathcal{C}_k$ be the family of $k$--element subsets of $L_i$ (i.e. possible configurations of exactly $k$ arcs from the sets $S_0$ and $S_1$) satisfying property $P$. Again, by abuse of notation the symbols $l_i$ we use to indicate these arcs will also denote the number of corresponding  arcs. 
%Also $X^{\star}_{2i+2}$ used in the formula is denoted $l_3$ for notational convenience. 

 %The arcs $l_{14}=X_{2i+12i+2}^{(0,0;-1,0)}$ and $l_2=X^{(1,0;-1,0)}_{2i+22i+2}$ which are used in the formulae are labeled with a star in Figure \ref{fig:illustration2} and Figure \ref{fig:illustration3}.

\begin{figure}[h!]
\begin{center}
\labellist
 \pinlabel {${\scriptstyle{l_0}}$} [ ] at  90 350

 \pinlabel {${\scriptstyle{l_1}}$} [ ] at  300 350
 \pinlabel {${\scriptstyle{l_2}}$} [ ] at  550 350
  \pinlabel {${\scriptstyle{l_3}}$} [ ] at   800 350
 \pinlabel {${\scriptstyle{l_{{\star}}}}$} [ ] at  1000 350

 \pinlabel {${\scriptstyle{l_4}}$} [ ] at   1230 350
 \pinlabel {${\scriptstyle{l_5}}$} [ ] at  1500 350
  \pinlabel {${\scriptstyle{l_6}}$} [ ] at  1730 350
\pinlabel {${\scriptstyle{l_7}}$} [ ] at  100 -15
 \pinlabel {${\scriptstyle{l_8}}$} [ ] at  330 -15
 \pinlabel {${\scriptstyle{l_9}}$} [ ] at  580 -15
  \pinlabel {${\scriptstyle{l_{10}}}$} [ ] at  780 -15
    \pinlabel {${\scriptstyle{l_{11}}}$} [ ] at  1000 -15

  \pinlabel {${\scriptstyle{l_{12}}}$} [ ] at  1280 -15
  \pinlabel {${\scriptstyle{l_{13}}}$} [ ] at  1480 -15

  \pinlabel {${\scriptstyle{l_{14}}}$} [ ] at  1720 -15

 \small\hair 2pt
\endlabellist  

  \includegraphics[scale=0.2]{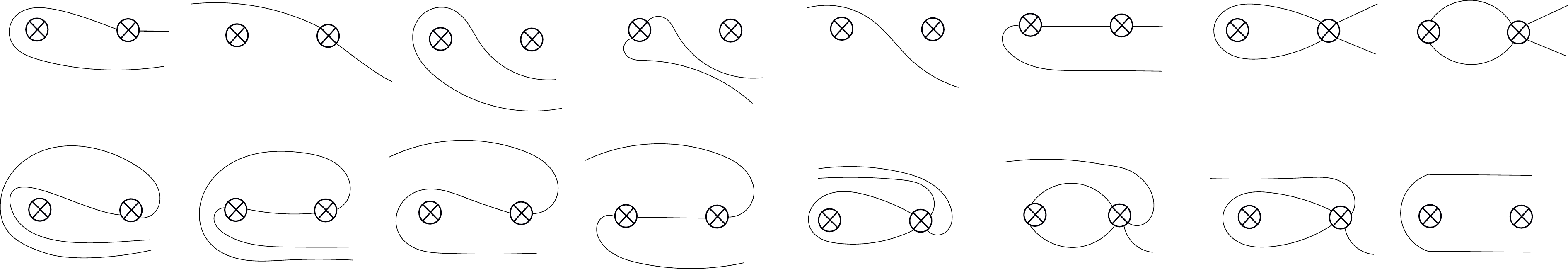}
\caption{The set $S_0$}\label{fig:illustration2}
\end{center}
\vspace{0.3 cm}

\begin{center}
\labellist
\pinlabel {${\scriptstyle{l_{15}}}$} [ ] at  100 -10
 \pinlabel {${\scriptstyle{l_{16}}}$} [ ] at  300 -10
 \pinlabel {${\scriptstyle{l_{17}}}$} [ ] at  550 -10
 \pinlabel {${\scriptstyle{l_{18}}}$} [ ] at  750 -10
 \pinlabel {${\scriptstyle{l_{19}}}$} [ ] at  1000 -10

 \pinlabel {${\scriptstyle{l_{20}}}$} [ ] at  1250 -10

 \small\hair 2pt
\endlabellist  

  \includegraphics[scale=0.2]{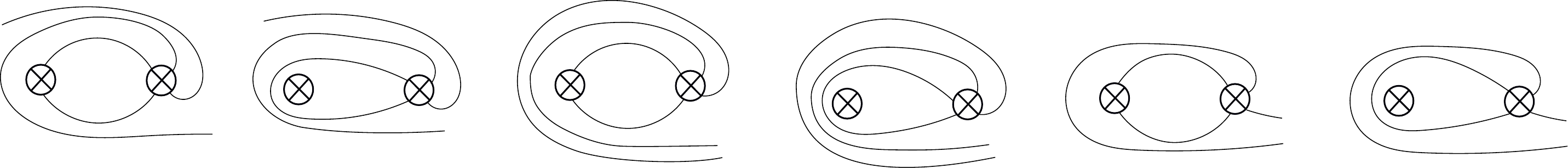}
\caption{The case where each arc $X$ in $S_{1}$ has $\chi(X)=1$}\label{fig:illustration3}
\end{center}
\end{figure}

\noindent First, observe from Figure \ref{fig:htexep} and Figure \ref{fig:halftwisted} that $l_0$ is a standard exceptional arc with respect to a clover type II.  Also, $l_{18}, l_{19}$ are examples for twisted exceptional arcs; and $\{l_1, l_9\}$, $\{l_4, l_8\}$ and $\{l_{12}, l_{15}\}$ form  exceptional arc systems of with respect to a clover type II.

\begin{table}[h!]
\begin{center}
\labellist
\small\hair 2pt

     \pinlabel {$\tiny{l_0}$} [ ] at  250 520
     \pinlabel {$\tiny{(1, 1, 1, 1)}$} [ ] at  240 450

         % \pinlabel {$\tiny{(0, 1)}$} [ ] at  240 350

      \pinlabel {$\tiny{l_1}$} [ ] at  650 520
     \pinlabel {$\tiny{(0, 1, 0, 0)}$} [ ] at  640 450
       %   \pinlabel {$\tiny{(0, 1)}$} [ ] at  640 350

         \pinlabel {$\tiny{l_2}$} [ ] at  1050 520
     \pinlabel {$\tiny{(1, 0, 2, 0)}$} [ ] at  1040 450
         % \pinlabel {$\tiny{(0, 0)}$} [ ] at  1040 350

 \pinlabel {$\tiny{l_3}$} [ ] at  1450 520
     \pinlabel {$\tiny{(0, 0, 2, 0)}$} [ ] at  1440 450
        %  \pinlabel {$\tiny{(1,0)}$} [ ] at  1440 350
        
         \pinlabel {$\tiny{l_{{\star}}}$} [ ] at  1850 520
     \pinlabel {$\tiny{(0, 0, 1,-1)}$} [ ] at  1880 450

 \pinlabel {$\tiny{l_4}$} [ ] at  2300 520
     \pinlabel {$\tiny{(0, 1, 1,1)}$} [ ] at   2300 450
         % \pinlabel {$\tiny{(1, 1)}$} [ ] at  1850 350

 \pinlabel {$\tiny{l_5}$} [ ] at   2700 520
     \pinlabel {$\tiny{(1, 2, 0, 2)}$} [ ] at  2700 450
         % \pinlabel {$\tiny{(0, 2)}$} [ ] at  2300 350

 \pinlabel {$\tiny{l_6}$} [ ] at  3180 520
     \pinlabel {$\tiny{(0, 2, 0, 2)}$} [ ] at  3180 450
          %\pinlabel {$\tiny{(1, 2)}$} [ ] at  2700 350

      \pinlabel {$\tiny{l_7}$} [ ] at  250 240
     \pinlabel {$\tiny{(2, 0, 2, 2)}$} [ ] at  240 150
         % \pinlabel {$\tiny{(0, 1)}$} [ ] at  240 50

      \pinlabel {$\tiny{l_8}$} [ ] at  650 240
     \pinlabel {$\tiny{(1, 0, 2, 2)}$} [ ] at  640 150
         % \pinlabel {$\tiny{(1, 1)}$} [ ] at  640 50

         \pinlabel {$\tiny{l_9}$} [ ] at  1050 240
     \pinlabel {$\tiny{(1, 0, 1, 1)}$} [ ] at  1040 150
        %  \pinlabel {$\tiny{(0, 1)}$} [ ] at  1040 50

 \pinlabel {$\tiny{l_{10}}$} [ ] at  1450 240
     \pinlabel {$\tiny{(0, 0, 1, 1)}$} [ ] at  1440 150
         % \pinlabel {$\tiny{(1, 1)}$} [ ] at  1440 50

 \pinlabel {$\tiny{l_{11}}$} [ ] at  1850 240
     \pinlabel {$\tiny{(1, 0, 0, 2)}$} [ ] at  1850 150
          %\pinlabel {$\tiny{(0, 2)}$} [ ] at  1850 50

 \pinlabel {$\tiny{l_{12}}$} [ ] at  2300 240

     \pinlabel {$\tiny{(0, 1, 0, 2)}$} [ ] at  2300 150
      %    \pinlabel {$\tiny{(1, 2)}$} [ ] at  2300 50

 \pinlabel {$\tiny{l_{13}}$} [ ] at  2700 240
     \pinlabel {$\tiny{(1, 1, 0, 2)}$} [ ] at  2700 150
      %   \pinlabel {$\tiny{(0, 2)}$} [ ] at  2700 50

 \pinlabel {$\tiny{l_{14}}$} [ ] at  3180 240
     \pinlabel {$\tiny{(1, 0, 1, 1)}$} [ ] at  3180 150

%  \pinlabel {$\tiny{(0, 1)}$} [ ] at  240 350
%    \pinlabel {$\tiny{(0, 1)}$} [ ] at  640 350
     %\pinlabel {$\tiny{(0, 0)}$} [ ] at  1040 350
        %\pinlabel {$\tiny{(1,0)}$} [ ] at  1440 350
            %    \pinlabel {$\tiny{(0,0)}$} [ ] at  1850 350

         %   \pinlabel {$\tiny{(1, 1)}$} [ ] at  2300 350
            %     \pinlabel {$\tiny{(0, 2)}$} [ ] at  2700 350
               %                   \pinlabel {$\tiny{(1, 2)}$} [ ] at  3200 350

                  %\pinlabel {$\tiny{(0, 2)}$} [ ] at  240 50
                     %   \pinlabel {$\tiny{(1, 1)}$} [ ] at  640 50
                        %  \pinlabel {$\tiny{(0, 1)}$} [ ] at  1040 50
                           %\pinlabel {$\tiny{(1, 1)}$} [ ] at  1440 50
             %  \pinlabel {$\tiny{(0, 2)}$} [ ] at  1850 50
    %       \pinlabel {$\tiny{(1, 2)}$} [ ] at  2300 50
       %    \pinlabel {$\tiny{(0, 2)}$} [ ] at  2700 50

          % \pinlabel {$\tiny{(0, 0)}$} [ ] at  3200 50

     \endlabellist

 \includegraphics[scale=0.13]{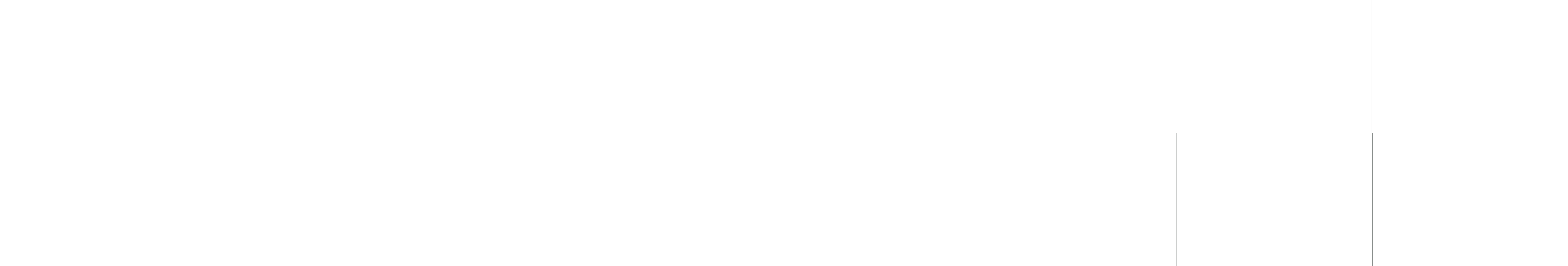}
\caption{$(\lambda_i, \psi_{i+1}, \Delta_{i+1}(B),  p_1)$ for the arcs of $S_0$} \label{fig:table}
\end{center}
\end{table}

%$p_1=c_{i+1}-(c_i-X^{\star}_{2i+2}) \quad \text{and}\quad and  the pairs $(c_i, c_{i+1})$  

\begin{itemize}
\item[Case 1.]  $\chi(i)\leq 0$: Each component of $L_i\setminus I $ belongs to $S_{0}$ where $\chi(l_2)= \chi(l_3)=\chi(l_{\star})=-1$;  and $\chi(l_k)=0$ for every other $l_k$ in $S_0$. As to be explained later in Remark \ref{rem:explanation}, we need another parameter $p_1=2(-b_{n+i-1})^+-\Delta_{i+1}(B)$. For simplicity, we list the 4--tuples $(\lambda_i, \psi_{i+1}, \Delta_{i+1}(B), p_1)$  in Table \ref{fig:table} corresponding to the arcs $l_k$ in $S_0$.  Observe from Table \ref{fig:table} that $l_0$ is the only element satisfying property $P$ alone hence $\mathcal{C}_1=\{\{l_0\}\}$. Furthermore, it  increases  each $\lambda_i, \psi_{i+1}$, $ \Delta_{i+1}(B)$ and $p_1$ by $1$, yielding $\epsilon'_i=N_i=l_0=\min\big(M_1,   p_1\big)$ as required.  In order to construct $\mathcal{C}_k$~$(k>1$) we make use of another set $I_k$ which is the set of $k$-tuples $(i_1, i_2,\dots, i_k)$ for the compatible components $l_{i_1}, l_{i_2},\dots, l_{i_k}$ where $i_j\in \{\star, k:0\leq k\leq 14\}$. We chose to use the star symbol here to indicate that $l_{\star}$ is only arc that is not compatible with any exceptional arc or arc system (in fact it is compatible with only $l_1, l_2$ and $l_3$). We get,  $I_2=I^{(1)}_2\cup I^{(2)}_2\cup I^{(3)}_2$, where $$I^{(1)}_2=\{(0, j): 1\leq j\leq 14\}, \quad I^{(2)}_2=\{(1,\star), (1, j): 1< j< 14, j\neq  7, 8\} \quad \text{and}$$  
\begin{align*}
 I^{(3)}_2=&\big\{(2,3), (2,{\star}), (2,7), (2,9), (2,14), (3,\star), (3,4), (3, 7), (3, 8), (3,9), (3, 10), (3,14), (4, 6), (4, {8}), \\&(4, {10}), (4,12), (4,14),  (5, 6), (5,13), (5,14), (6,12), (6,13), (6,14), (7, 8), (7,{9}),(7, 10), (7,14), \\&(8, 9), (8, {10}), (8,14), (9, {10}), (9, 11), (10,12),(11,13), (12,13)\big\}.\end{align*}
 \noindent and obtain that  $\mathcal{C}_2=\mathcal{C}^{(1)}_2\cup \mathcal{C}^{(2)}_2$ where $$ \mathcal{C}^{(1)}_2=\big\{\{l_0, l_j\}: 1\leq j\leq 14\big\}\quad \text{and}\quad  \mathcal{C}^{(2)}_2=\big\{ \{l_1, l_2\}, \{l_1, l_{9}\}, \{l_4, l_8\}, \{l_4, l_{14}\}, \{l_5, l_{14}\}, \{l_6, l_{14}\} \big\}.$$

\noindent First consider $\mathcal{C}^{(1)}_2$. Recall that to compute the parameters  associated with a compatible set we simply add the corresponding components of arcs. For example,  from Table \ref{fig:table} we compute that $(M_1, p_1)=(1,1)$ for $(l_0, l_1)$.  We immediately check that for each $\{l_0, l_j\}\in \mathcal{C}^{(1)}_2$ we have $l_0=N'_i=M_1\leq p_1$, and  and therefore $\epsilon'_i=N'_i=l_0=\min\big(M_1, p_1\big)$ (note that taking multiple copies of arcs does not change the formula).   We continue with $\mathcal{C}^{(2)}_2$. The set $ \{l_1, l_2\}$ contains no exceptional arc or arc systems with respect to a clover type II, and observe from Table \ref{fig:table} as above that $p_1=0$. Therefore,  $\epsilon'_i=0=\min(M_1, p_1)$ as required.  Similarly, we check $\{l_1, l_9\}$. This set contains negatively half twisted scissors (Figure \ref{fig:halftwisted}(a)) but no exceptional arcs hence we have $\epsilon'_i=s'_i$. We check that $l_1=\psi_{i+1}$ and $ l_9=\lambda_i=\Delta_{i+1}$. Therefore, $s'_i=M_1=\min (l_1, l_9)$. Furthermore, $p_1$ is increased  by $1$ by both $l_1$ and $l_9$. This implies $M_1\leq p_1$ yielding  $\epsilon'_i=s'_i=\min\big(M_1, p_1\big)$ as required. Similarly,  $\{l_4, l_8\} $ is a negatively half twisted left anchor (Figure \ref{fig:halftwisted}(b))  and hence $\epsilon'_i=z'_i$. We check that $z'_i=\min(l_4, l_8)=M_1\leq p_1$ yielding $\epsilon'_i=z'_i=\min\big(M_1, p_1\big)$ as required. Finally, none of $(l_j, l_{14})\in \mathcal{C}^{(2)}_2$ contains an exceptional arc or arc system. Since $l_{14}=X_{2i+12i+2}^{(0,0;-1,0)}$ we get $M_1=0$ yielding $\epsilon'_i=0$ as expected.
The formula can be verified similarly for elements of $\mathcal{C}_3$ as follows. 

		We have $I_3=I^{(1)}_3\cup I^{(2)}_3\cup I^{(3)}_3$ where $I^{(1)}_3=\{(0, j, k): (j, k)\in I_2, j, k\neq \star\}$ and $I^{(2)}_3=\{(1, j, k):  (j, k)\in I_2, j, k\notin \{7, 8, 14\}\}$ and

\begin{align*}
 I^{(3)}_3=\big\{&(2,3,7), (2,3,9),  (2,3,\star),  (2,3,14), (2,7,9), (3, 4, 8), (3,4,10),  (3,4,14), (3,7,8), (3,7,9), (3,7,10),  \\&(3,7,14), (3,8,9), (3,8,10), (3,8,14), (3,9,10),  (4, 6,12), (4, 6,14), 
 (4, 8,10), (4, 8,14), (4, 10,12), \\&  (5, 6, 13), (5, 6, 14), (6,12, 13), (7, 8, 9), (7, 8, 10), (7, 8, 14),  (7, 9, 10), (8, 9, 10)\big\}.
\end{align*}
\noindent Hence, $\mathcal{C}_3=\mathcal{C}^{(1)}_3\cup \mathcal{C}^{(2)}_3$ where $  \mathcal{C}^{(1)}_3=\big\{\{l_0, l_j, l_k\}: (j, k)\in I_2, j, k\neq \star\big\}$ and
 %If $(j, k)\in I^{(1)}_2$ we obtain $\{l_0, l_k\}$  for $1\leq k\leq 13$ and the proof is similar to that of $\mathcal{C}^{(1)}_2$.
\begin{align*}\mathcal{C}^{(2)}_3=\big\{ &\{l_1, l_2, l_3\},\{l_1, l_2, l_{\star}\}, \{l_1, l_2, l_9\}, \{l_1, l_3, l_9\},  \{l_1, l_9, l_{10}\}, \{l_1, l_9, l_{11}\}, \{l_3, l_4, l_8\}, \{l_3, l_4, l_{14}\}, \\
&\{l_4, l_6, l_{14}\}, \{l_4, l_8, l_{10}\}, \{l_4, l_8, l_{14}\},  \{l_5, l_6, l_{14} \}\big\}. \end{align*}

%
%
%Observe from Table \ref{fig:table} that $l_0=M_1$ where $M_1=\psi_{i+1}-l_1$ for $k=2, 9$;   $M_1=\psi_{i+1}-l_1=\lambda_i$ for $k=3,10$; $M_1=\lambda_i$ for $k=4$; $M_1=\Delta_{i+1}(B)$ for $k=5, 13$; $M_1=\lambda_i=\Delta_{i+1}(B)$ for $k=6, 12$ and  $M_1=\psi_{i+1}-l_1=\Delta_{i+1}(B)$ for $k=11$.  Furthermore, for each such $\{l_0, l_1, l_k\}$ we have $c_{i+1}-2[\chi_{i+1}(B)]^+\geq c_{i+1}-(c_i-X^{\star}_{2i+2}+[\chi_{i+1}(B)]^+)\geq  M_1$. 

%that $(4,8)$ is the only pair such that $\{l_4, l_8\}$ satisfies property $P$. Therefore, $M_1=0$ for each  $\{l_j, l_k\}$ with $(j, k)\neq (4,8)$. We immediately check from Table \ref{fig:table} that for each such $\{l_0, l_1, l_k\}$ we have 
%$$M_1\leq c_{i+1}-(c_i-X^{\star}_{2i+2}+[\chi_{i+1}(B)]^+)\leq c_{i+1}-2[\chi_{i+1}(B)]^+$$ that is $l_0=\min(M_1, c_{i+1}-(c_i-X^{\star}_{2i+2}+\chi_{i+1}(B)), c_{i+1}-2[\chi_{i+1}(B)]^+)$ as required. If $(j, k)=(4,8)$ we have $\{l_0, l_4, l_8\}$, and similar arguments give that $l_0=c_{i+1}-c_i-X^{\star}_{2i+2}+[\chi_{i+1}(B)]^+)$ which equals $$\min(M_1, c_{i+1}-(c_i-X^{\star}_{2i+2}+[\chi_{i+1}(B)]^+, c_{i+1}-2[\chi_{i+1}(B)]^+)$$
%\marginnote{Affect etmeyenler $M_1$ property P saglarlarsa mutlaka $X^{\star}_{2i-1}$ vardir sistemde bu da $l_0=0$ verir. Ornegin-fotografa bak}
%

%We have $\{l_0, l_j, l_k\}\in \mathcal{C}^{(1)}_3$.

First consider $ \mathcal{C}^{(1)}_3$.  If $(j, k)\in I^{(2)}_1$ we get $\{l_0, l_0, l_k\}~(1\leq k\leq 14)$, and the proof is similar to that of $\mathcal{C}^{(1)}_2$. If $(j, k)\in I^{(2)}_2$ we get $\{l_0, l_1, l_k\}$~($1\leq k\leq 13$, $k\neq 7, 8$). We compute from Table \ref{fig:table} that for each $k$ with $k\neq 2, 9$ we have that $\epsilon'_i=l_0=N'_i=M_1=\min (M_1, p_1)$.  For $k=2$ we get that $\epsilon'_i=N_i=l_0=p_1=\min(M_1, p_1)$. Similarly for $k=9$ we compute that $\epsilon'_i=s'_i+N'_i=\min(l_1, l_9)+l_0=\min(M_1, p_1)$ ($L_i$ contains half twisted scissors $\{l_1, l_9\}$).  If $(j, k)\in I^{(3)}_2$ we have $\epsilon'_i=l_0=M_1=\min(M_1, p_1)$ since $M_1=0$ for each  corresponding $(l_j, l_k)$; and for $\{l_0, l_4, l_8\}$ we have $\epsilon'_i=z'_i+N'_i=\min(l_4, l_8)+l_0=M_1=\min(M_1, p_1)$ ($L_i$ contains half twisted anchor $\{l_4, l_8\}$). Finally for  $\{l_j, l_k, l_m\}\in \mathcal{C}^{(2)}_3$ we similarly verify from Table \ref{fig:table} that $\epsilon'_i=p_1=0$ for $ (i, j, k)\in \{(1, 2, 3),(1, 2, {\star})\}$; $\epsilon'_i=s'_i=\min(M_1, p_1)$ for $ (i, j, k)\in \{(1, 2, 9), (1, 3, 9),  (1, 9, {10}), (1, 9, 11)\}$; $\epsilon'_i=M_1=\min(M_1, p_1)$ for  $\{l_3, l_4, l_8\},  \{l_4, l_8, l_{10}\}, \{l_4, l_8, l_{14}\}$ and $\epsilon'_i=M_1=0=\min(M_1, p_1)$  
for $ \{l_3, l_4, l_{14}\}, \{l_4, l_6, l_{14}\}, \{l_5, l_6, l_{14} \}$. Similarly,  we  have $I_4=I^{(1)}_4\cup I^{(2)}_4\cup I^{(3)}_4$ where 
$$I^{(1)}_4=\{(0, j, k, m): (j, k, m)\in I_3, j, k, m \neq\star\}, I^{(2)}_4=\{(1, j, k,  m):  (j, k, m)\in I_3, j, k, m\notin \{7, 8\}\}$$ and $I^{(3)}_4=\big\{(2,3,7,9),  (3, 4, 8,10), (3,7,8, 9), (3,7, 8,10),( 3,8, 9,10),  (7, 8,9,10)\big\}$.

We get  $\mathcal{C}_4=\mathcal{C}^{(1)}_4\cup \mathcal{C}^{(2)}_4$ where $  \mathcal{C}^{(1)}_4=\big\{\{l_0, l_j, l_k, l_m\}:(j, k, m)\in I_3, j, k, m \neq\star\big\}$ 
%First let $\{l_0, l_j, l_k, l_m\}\in \mathcal{C}^{(1)}_4$. If $(j, k, m)\in I^{(1)}_3$ we obtain $\{l_0, l_0, l_k, l_m\}$; and
\begin{align*}\mathcal{C}^{(2)}_4=\big\{\{&l_1, l_2, l_3, l_9\}, \{l_1, l_2, l_3, l_\star\},  \{l_1, l_3, l_9, l_{10}\},   \{l_3, l_4, l_8, l_{10}\},  \{l_3, l_4, l_8, l_{14}\}\big\}. \end{align*}  Also,  $\mathcal{C}_5=\mathcal{C}^{(1)}_5\cup \mathcal{C}^{(2)}_5$ where $  \mathcal{C}^{(1)}_5=\big\{\{l_0, l_j, l_k, l_m, l_n\}: (j, k, m, n)\in I_4,  j, k, m, n \neq\star \big\}$ and $\mathcal{C}^{(2)}_5=\emptyset$ since there is no $5$ element compatible set which doesn't contain $l_0$ but satisfies property $P$. And finally,  $\mathcal{C}_6=\{l_0, l_3, l_7, l_8, l_9, l_{10}\}$. We note that there is no $\mathcal{C}_k$ with $k\geq 6$. The verification of the formula for $k=4, 5, 6$ is analogous. 

\item[Case 2.]   $\chi_i>0$: At least a component of $L_i\setminus I$ belongs to the set $S_{1}$. There are 2 subcases depending on whether or not 
$L_i$ contains a highly twisted arc $X$ with $\chi(X)>1$.
 \begin{itemize}

\item[(a)] No component  $X$ of $L_i$ has $\chi(X)>1$: Then $L_i$ contains a highly twisted arc $X$ which has $\chi(X)=1$ (Figure \ref{fig:illustration3}). Let $\mathcal{A}_k$ denote the set of arcs  which are compatible with the arc $l_k\in S_1$. Then, 
\begin{align*}\mathcal{A}_{15}&=\{l_0, l_8, l_9, l_{10}, l_{12}, l_{16}, l_{19}, l_{20}\}, \mathcal{A}_{16}=\{l_0, l_9, l_{11}, l_{10}, l_{15}, l_{20}\},\mathcal{A}_{17}=\{l_0, l_8, l_9, l_{18}, l_{19}\},\\
\mathcal{A}_{18}&=\{l_0, l_7, l_9, l_{17}, l_{20}\},\mathcal{A}_{19}=\{l_0, l_4, l_6, l_8, l_{10}, l_{17}, l_{15}\},\mathcal{A}_{20}=\{l_0, l_5, l_6, l_{11}, l_{15}, l_{18}\}\end{align*} 

 Any compatible set containing  $l_k\in S_1$ is constructed from elements of $\mathcal{A}_{k}$ such that property $P$ is satisfied.  Therefore,  a compatible set $\cC$ can contain the standard exceptional arc $l_0$ (which is compatible with each element of $S_1$) and the twisted exceptional arcs $ l_{19}$ and $l_{20}$. The only exceptional arc system $\cC$ can contain is the negatively half twisted ribbon which is the arc system  $\{l_{12}, l_{15}\}$. Consider  for example the compatible sets containing $l_{15}$. Each such set  is constructed from $\mathcal{A}_{15}$ in such a way that  it contains at least one of $l_0, l_{12}, l_{19}$ and $l_{20}$ so that $\psi_{i+1}\neq 0$ (the other two assumptions $\lambda_i\neq 0$ and $\Delta_{i+1}(B)\neq 0$ are satisfied by each element in $S_1$).   We immediately check that  for each such compatible set we have $p_1\geq M_1$ and $M_1=N'_i+r'_i$ where $r'_i=\min(l_{12}, l_{15})$.  Therefore, $\epsilon'_i=N'_i+r'_i=\min(M_1, p_1)$ as required.

 \item[(b)]  Some component  $X$ of $L_i$ has $\chi(X)>1$: Then $l_0=0$, $s'_i=z'_i=r'_i=0$ by Remark \ref{rem:notcompatible} since each exceptional arc system with respect to a clover type II contains a standard arc. Since $\psi_{i+1}\neq 0$ by assumption there exists a highly twisted exceptional arc $X\in [X_{2i+2; i+1}]$ (such arcs are the only highly twisted arcs with $\psi_{i+1}\neq 0$ and satisfying property $P$)  each of which increases $M_1$ by $1$. Since $M_1=0$ for any other arc compatible with $X$ we get   $\epsilon'_i=N'_i=\min(M_1, p_1)$ as required.\end{itemize}  \end{itemize} \end{proof}
 
 \begin{remark}\label{rem:explanation} The proof of Lemma \ref{lem:clover2coordinates} shows that there exists compatible sets satisfying property $P$ yet containing no exceptional arcs or arc systems. Such arc systems either contain $\{l_1,l_2\}$ or the arc $X_{2i+12i+2}=l_{14}$ together with an arc that satisfies $\psi_{i+1}\neq 0$ such as $\{l_5, l_{14}\}$ (note that $l_{14}$ contributes to both $\lambda_i$ and $\Delta_{i+1}(B)$). Using parameter $p_1$ and subtracting $X_{2i+12i+2}$ from $\lambda_i$ and $\Delta_{i+1}(B)$ rules out such arc systems  giving a way to compute only exceptional arcs or arc systems.\end{remark}
% In particular, using $p_1$ guarantees that $\epsilon'_i=0$ for compatible sets containing $\{l_1,l_2\}$, and  subtracting $X_{2i+12i+2}$ from $\lambda_i$ and $\Delta_{i+1}(B)$ guarantees that $\epsilon'_i=0$ for compatible sets containing $l_{14}$

  \begin{remark}\label{reflectionhori}~Reflection in the horizontal diameter of the surface conjugates each
  crosscap transposition ~$u_i$ to $u_i^{-1}$. Therefore a clover of type III is the reflection of a clover of type II along the horizontal diameter, 
and the corresponding transformation of generalized Dynnikov coordinates in max-plus notation is given by $[t_i; b_i]\mapsto[1/t_i; b_i]$. For example for $n>0$
\[(t_1,\ldots, t_{g-1},b_1,\ldots,b_{n+g-2}) \mapsto
\tropical{(1/t_1,\ldots,1/t_{g-1},b_1,\ldots,b_{n+g-2})}.\] 

\end{remark} 
  
  By Remark \ref{reflectionhori} we conclude that exceptional arcs  with respect to a clover of type III can be obtained by reflecting exceptional arcs with respect to a clover of type II in the horizontal diameter. Therefore, replacing $u^{-1}_i$ with $u_i$ in Notation \ref{not:cloverII}, we obtain the exceptional parameter $\bar{e_i}=\epsilon''_i+\epsilon''_{i+1}$ for a clover of type III as given in Lemma \ref{lem:clover3coordinates}  and Lemma \ref{lem:clover3}.

%Similarly,  taking the $u_i$ images of exceptional arcs  and arc systems with respect to a clover type I in $S'_i\cup S'_{i+1}$ we obtain {\it exceptional arcs} (i.e. elements of $[X_{2i ; i}]$ and  $[X_{2i+1; i+1}]$) and  {\it exceptional arc systems} with respect to a clover of type III (i.e. positively half twisted scissors, anchors and ribbons). 

\begin{lem}\label{lem:clover3coordinates}

Let $L_i$ be positive.  Then $\bar{e}_i=\epsilon''_{i}+ \epsilon''_{i+1}$ where 
\begin{align}\label{eq:eibar}
\epsilon''_{i}&=\min({\Lambda}_i, \psi_{i+1}, \bar{\Delta}_{i+1}(A), p_3);\\ \epsilon''_{i+1}&=\min({\Lambda}_{i+1}, \psi_{i}, \bar{\Delta}_{i}(B), p_4)
\end{align}

\end{lem}
\begin{lem}[Equality for a clover of type III]\label{lem:clover3} 
Given a clover of type III we have
 \begin{eqnarray}\label{eq:clover3}
2\beta''_{n+i}+C_i=\max(\gamma_{2i-1}+\gamma_{2i+2}, \gamma''_{2i-1}+\gamma''_{2i+2})+2\bar{e}_i.
\end{eqnarray}\end{lem}

 \subsection{Scale equalities}\label{sec:scales}
 
 Let  $u^{-1}_i(\gamma;\,\beta)=(\gamma';\,\beta')$ and $u_i(\gamma;\,\beta)=(\gamma'';\,\beta'')$.  A scale of type I has leaves $\gamma_{2i+1}$, $\gamma_{2i}$, $\cC_i$, $\beta_{n+i+1}$; diagonals $\gamma'_{2i+1}$ and $\gamma_{2i+2}$; and a scale of type II has leaves $\gamma_{2i+1}$, $\gamma_{2i}$, $\cC_i$, $\beta_{n+i-1}$; and diagonals $\gamma'_{2i}$ and $\gamma_{2i-1}$. Reflecting these two scales along the horizontal diameter we respectively obtain a scale of type III and a scale of type IV (Figure \ref{fig:scaleI&II}). Observe that this is natural since a scale of type III and a scale of type IV are the $u_i$ images of a  scale of type I and a scale of type II respectively (see Remark \ref{reflectionhori}).

 %We present equalities of a scale of type I and type II given in Lemma \ref{lem:scales1} and Lemma  \ref{lem:scales2} respectively. We begin with a scale of type
 %I as shown in Figure \ref{fig:scaleI&II}(a). 

 \begin{figure}[h!]
\begin{center}
\labellist

 \pinlabel {${\scriptstyle{\gamma_{2i+1}}}$} [ ] at  138 250
      \pinlabel {\begin{turn}{-90}$\scriptstyle{\gamma'_{2i+1}}$\end{turn}} [ ] at  240 220

 % \pinlabel {${\scriptstyle{\gamma'_{2i+1}}}$} [ ] at  255 200
 \pinlabel {${\scriptstyle{\gamma_{2i+2}}}$} [ ] at  147 35
 \pinlabel {${\scriptstyle{C_{i}}}$} [ ] at  110 217
     \pinlabel {\begin{turn}{-90}$\scriptstyle{\beta_{n+i+1}}$\end{turn}} [ ] at  298 200
 \pinlabel {${\scriptstyle{\gamma_{2i}}}$} [ ] at  50 38
 
 \pinlabel {\begin{turn}{-90}$\scriptstyle{\beta_{n+i-1}}$\end{turn}} [ ] at  370 200

  \pinlabel {${\scriptstyle{\gamma'_{2i}}}$} [ ] at  420 40

  \pinlabel {${\scriptstyle{\gamma_{2i+1}}}$} [ ] at  605 240
  \pinlabel {${\scriptstyle{\gamma_{2i-1}}}$} [ ] at  500 240
% \pinlabel {${\scriptstyle{C_{i}}}$} [ ] at  650 70
  \pinlabel {${\scriptstyle{\gamma_{2i}}}$} [ ] at  489 40
% \pinlabel {${(b)}$} [ ] at  610 -20
% \pinlabel {${(a)}$} [ ] at  190 -20

 \pinlabel {\begin{turn}{-90}$\scriptstyle{\beta_{n+i+1}}$\end{turn}} [ ] at  968 200

 \pinlabel {\begin{turn}{-90}$\scriptstyle{\gamma''_{2i+2}}$\end{turn}} [ ]   at  900 80
      \pinlabel {${\scriptstyle{\gamma_{2i+1}}}$} [ ] at  890 250

  \pinlabel {${\scriptstyle{\gamma_{2i-1}}}$} [ ] at  700 250
% \pinlabel {${\scriptstyle{C_{i}}}$} [ ] at  650 70
  \pinlabel {${\scriptstyle{\gamma_{2i+2}}}$} [ ] at  800 40
    \pinlabel {${\scriptstyle{\gamma_{2i+2}}}$} [ ] at  1210 40
        \pinlabel {${\scriptstyle{\gamma_{2i}}}$} [ ] at  1110 40

 \pinlabel {\begin{turn}{-90}$\scriptstyle{\beta_{n+i-1}}$\end{turn}} [ ] at  1025 200

      \pinlabel {\begin{turn}{-90}$\scriptstyle{\gamma''_{2i-1}}$\end{turn}} [ ] at  1090 250

  \pinlabel {${\scriptstyle{\gamma_{2i-1}}}$} [ ] at  1180 250

% \pinlabel {${\scriptstyle{C_{i}}}$} [ ] at  650 70

% \pinlabel {${(b)}$} [ ] at  610 -20
% \pinlabel {${(a)}$} [ ] at  190 -20

\small\hair 2pt
\endlabellist  
 \includegraphics[scale=0.35]{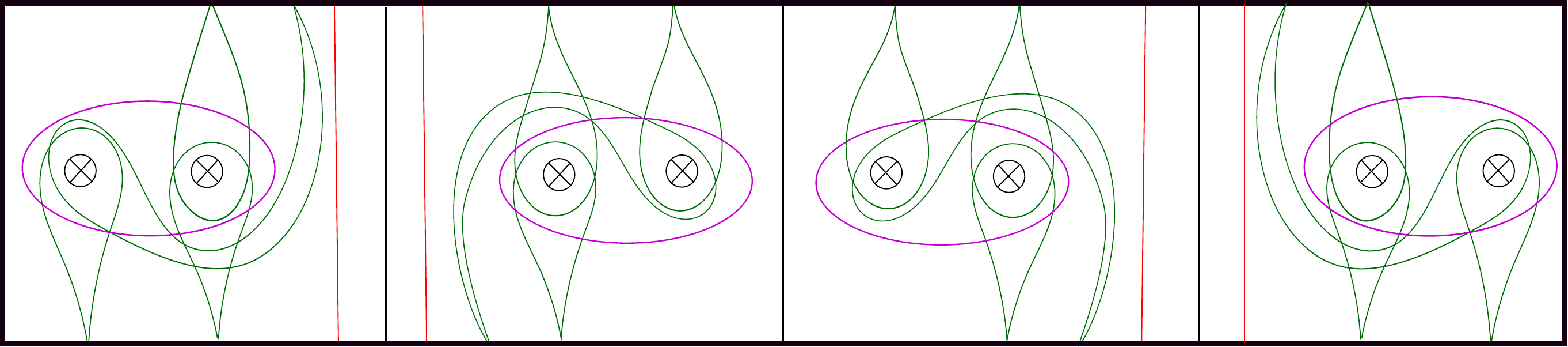}

\caption{A scale of type I, type II,  type III and  type IV from left to right} \label{fig:scaleI&II}
\end{center}
\end{figure}

\begin{notation}\label{not:Xpar}Let $X$ be a component of $L_i$.  In what follows we shall write $\lambda(X), \lambda_{c_k}(X), \psi_k(X)$ ($k=i, i+1$) to denote the number of non--core, core and straight components of a given component $X$ of $L\cap S'_k$.  
\end{notation}

The key idea in the proof of Lemma \ref{lem:scales1} is that it is easy to find out  which standard arcs satisfy equality  (\ref{eq:scale1a}), (\ref{eq:scale1b}), (\ref{eq:scale1c}) and (\ref{eq:scale1d}) since there are only finitely many standard arcs to check. Similarly,  we call  arcs which don't satisfy equality  (\ref{eq:scale1a}),   (\ref{eq:scale1b}),   (\ref{eq:scale1c})  and   (\ref{eq:scale1d})   \emph{exceptional with respect to a scale of type I, type II, type III and  type IV} respectively. 
 \begin{align}\label{eq:scale1a}
\gamma'_{2i+1}+\gamma_{2i+2}=\max(\gamma_{2i+1}+\gamma_{2i}, 	C_i+2\beta_{n+i+1}) \end{align}
 \begin{align}\label{eq:scale1b}
\gamma'_{2i}+\gamma_{2i-1}=\max(\gamma_{2i+1}+\gamma_{2i}, 	C_i+2\beta_{n+i-1}) \end{align}
\begin{align}\label{eq:scale1c}
\gamma''_{2i+2}+\gamma_{2i+1}=\max(\gamma_{2i+2}+\gamma_{2i-1}, 	C_i+2\beta_{n+i+1}) \end{align}
\begin{align}\label{eq:scale1d}
\gamma''_{2i-1}+\gamma_{2i}=\max(\gamma_{2i-1}+\gamma_{2i+2}, 	C_i+2\beta_{n+i-1}) \end{align}

We say that $X$ is straight in $S'_k$ ($k=i, i+1$) if $\psi_k(X)\neq 0$. Similarly, we say that $X$ is $u_i$--straight in $S_{k}$ if $\psi_k(u_i(X))\neq 0$ for some $k=i, i+1$.  Definition for  $u^{-1}_i$--straight arc in $S_{k}$ is similar.

 Analysis of exceptional arcs with respect to a scale of type I shows that each standard component $X$ of $L_i$ which isn't straight  in $S'_{i+1}$ (i.e. those with $\psi_{i+1}=0$) satisfies equality (\ref{eq:scale1a}).  An analogous statement for a twisted component in the class of $X$ is also true since each such component has the same  number of intersections with $\cC_i$ as $X$, and increases the number of intersections on each $\gamma'_{2i+1}$, $\gamma_{j}~(2i\leq j\leq 2i+2)$ and $2\beta_{n+i}$ (and hence $C_i+2\beta_{n+i+1}$)  by the same amount.  Also the only standard straight arcs in   $S'_{i+1}$ which satisfy equality (\ref{eq:scale1a}) are $X_{2i+2; i+1}^{(0,1;-1,0)}$ and $X_{2i-1; i+1}$ (Figure \ref{fig:nonexceptionals}(a)). That is, each standard straight arc in $S'_{i+1}$  apart from $X_{2i+2; i+1}^{(0,1;-1,0)}$ and $X_{2i-1; i+1}$ is  exceptional with respect to a scale of type I (top row of Figure \ref{fig:illustration4}). Also the only twisted arcs which are straight in $S'_{i+1}$ are in the class of $[X_{k; i+1}]$ ($2i-1\leq k\leq 2i+2$) and $[X_i]$ (examples of which are as given on the bottom row of  Figure \ref{fig:illustration4}). Each  such exceptional arc $X$ satisfies $$\gamma'_{2i+1}+\gamma_{2i+2}=\max(\gamma_{2i+1}+\gamma_{2i}, 	C_i+2\beta_{n+i+1})+2\psi_{i+1}(X)$$

\begin{figure}
\begin{center}

  \includegraphics[scale=0.2]{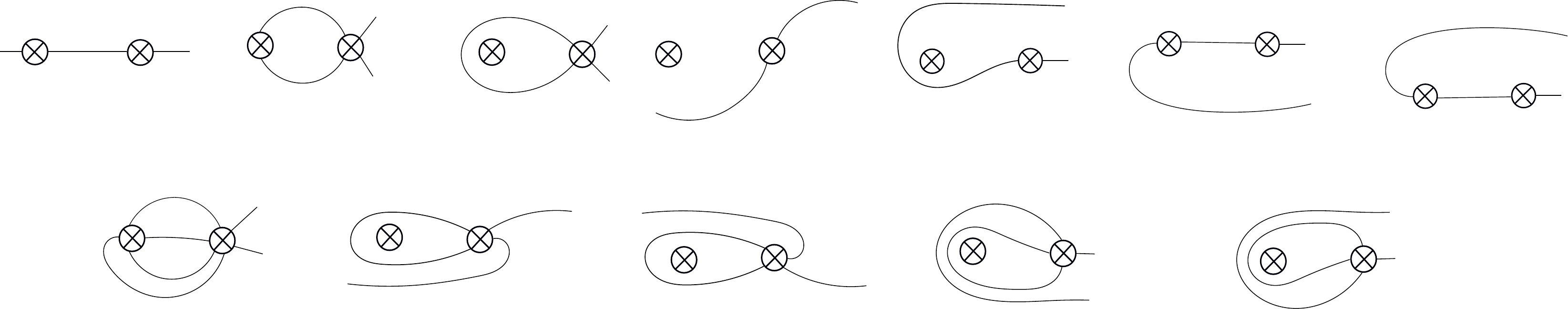}
\caption{Top row shows standard and bottom row shows examples for twisted arcs that are straight  in $S'_{i+1}$ }\label{fig:illustration4}

\end{center}
\begin{center}

  \includegraphics[scale=0.2]{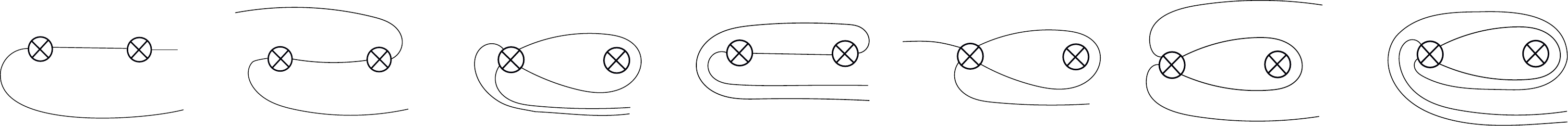}
\caption{Examples for $u_i$--straight arcs in $S'_{i+1}$}\label{fig:illustration4b}
\end{center}

\begin{center}
\labellist
\small\hair 2pt

\endlabellist  
  \includegraphics[scale=0.3]{nonexceptionals}
\caption{Straight  arcs and $u_i$--straight arcs in $S'_{i+1}$ which aren't exceptional with respect to a scale of type I}\label{fig:nonexceptionals}
\end{center}
\end{figure}
Since a scale of type III is the $u_i$ image of a  scale of type I it follows from Remark \ref{reflectionhori} that  each arc $X$ with $\psi_{i+1}(X)\neq 0$  apart from $X_{2i+1; i+1}^{(0,1;-1,0)}$ and $X_{2i; i+1}$ is exceptional with respect to a scale of type III. For the same reason, each $u_i$--straight arc $X$ in $S'_{i+1}$ (Figure \ref{fig:illustration4b}) apart from $X^{(1,0;-1,0)}_{2i+12i+2}$ and $X_{2i+2 ; i}$    (Figure \ref{fig:nonexceptionals}(b)) is exceptional with respect to a scale of type I. Let us write $\psi_{k}(u_i(X))=\psi_{k}(X')$ and  $\psi_{k}(u^{-1}_i(X))=\psi_{k}(X'')$. Then for each $u_i$--straight arc $X$ in $S'_{i+1}$ we get
 $$\gamma'_{2i+1}+\gamma_{2i+2}=\max(\gamma_{2i+1}+\gamma_{2i}, 	C_i+2\beta_{n+i+1})+2\psi_{i+1}(X')$$
  Write  $\psi_{k}(u_i(L))=\psi'_{k}$ and  $\psi_{k}(u^{-1}_i(L))=\psi''_{k}$ for a given $L\in \mathcal{L}_n$. 

Then, setting \begin{align}f_i=\psi_{i+1}-(X_{2i+2; i+1}^{(0,1;-1,0)}+X_{2i-1; i+1})+\psi'_{i+1}-(X^{(1,0;-1,0)}_{2i+12i+2}+X_{2i+2 ; i}) \end{align}
\noindent we obtain equality (\ref{eq:scale1aa}) given in Lemma \ref{lem:scales1}. We call $f_i$ the exceptional parameter for a scale of type I. The exceptional parameters for a scale of type II, type III and type IV follow from symmetry:

 \begin{figure}[h!]
\begin{center}
\labellist
\small\hair 2pt
\pinlabel {${\scriptstyle{X_{2i-1;i+1}}}$} [ ] at  -50 150
    \pinlabel {${\scriptstyle{X_{2i;i+1}}}$} [ ] at  -70 20
    
 \pinlabel {${\scriptstyle{X_{2i+2;i}}}$} [ ] at  640 140
    \pinlabel {${\scriptstyle{X_{2i+1;i}}}$} [ ] at  640 20
\endlabellist  
\includegraphics[scale=0.22]{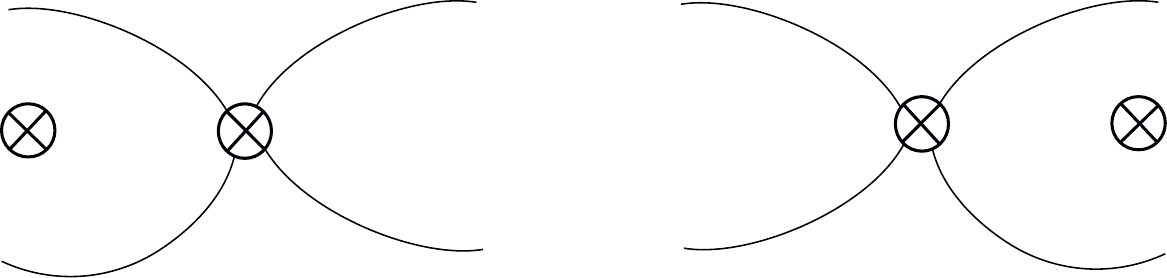}
%\caption{ }\label{fig:illustration4}
\end{center}
\vspace{1 cm}
\begin{center}
\labellist
\small\hair 2pt
 \pinlabel {${\scriptstyle{X^{(0,1;-1,0)}_{2i+2; i+1}}}$} [ ] at  30 100
 \pinlabel {${\scriptstyle{X^{(0,1;-1,0)}_{2i+1; i+1}}}$} [ ] at  300 100
 \pinlabel {${\scriptstyle{X^{(1,0;0,1)}_{2i; i}}}$} [ ] at  500 100
 \pinlabel {${\scriptstyle{X^{(1,0;0,1)}_{2i-1; i}}}$} [ ] at  700 100
\endlabellist  
\includegraphics[scale=0.3]{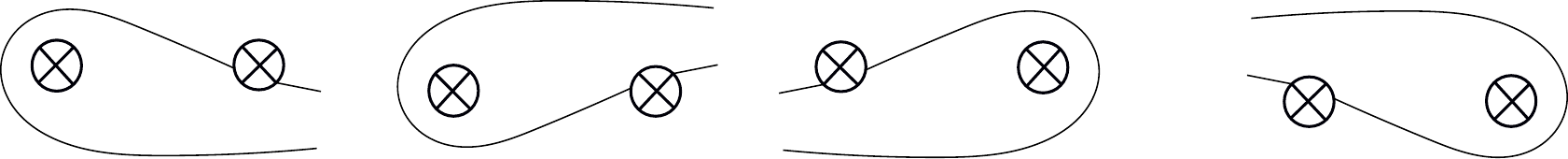}
\caption{ }\label{fig:alone}
\end{center}
\end{figure}

 \begin{align}\label{eq:gi}
g_i&=\psi_{i}-(X_{2i-1; i}^{(1,0;0,1)}+X_{2i+2; i})+\psi'_{i}-(X^{(0,1;0,1)}_{2i-12i}+X_{2i-1; i+1})
 \end{align}
 \begin{align}\label{eq:fibar}
 \bar{f_i}&=\psi_{i+1}-(X_{2i+1; i+1}^{(0,1;-1,0)}+X_{2i; i+1})+\psi''_{i+1}-(X^{(1,0;-1,0)}_{2i+12i+2}+X_{2i+1 ; i})
 \end{align}
  \begin{align}\label{eq:gibar}
 \bar{g_i}&=\psi_{i}-(X_{2i; i}^{(1,0;0,1)}+X_{2i+1; i})+\psi''_{i}-(X^{(0,1;0,1)}_{2i-12i}+X_{2i; i+1})
 \end{align}

\noindent Hence, computing $f_i, \bar{f_i}, g_i, \bar{g_i}$ in terms of generalized  Dynnikov coordinates will require separate consideration of the arcs depicted in Figure \ref{fig:alone} and given in Lemma \ref{lem:except3} and Lemma \ref{lem:htexep}. We first state scale equalities in Lemma \ref{lem:scales1}.

 \begin{lem}\label{lem:scales1} 
Given a scale of type I, type II, type III and type IV as shown in Figure \ref{fig:scaleI&II} we have 
 \begin{align}\label{eq:scale1aa}
\gamma'_{2i+1}+\gamma_{2i+2}=\max(\gamma_{2i+1}+\gamma_{2i}, 	C_i+2\beta_{n+i+1})+2f_i \end{align}
 \begin{align}\label{eq:scale1bb}
\gamma'_{2i}+\gamma_{2i-1}=\max(\gamma_{2i+1}+\gamma_{2i}, 	C_i+2\beta_{n+i-1})+2g_i \end{align}
\begin{align}\label{eq:scale1cc}
\gamma''_{2i+2}+\gamma_{2i+1}=\max(\gamma_{2i+2}+\gamma_{2i-1}, 	C_i+2\beta_{n+i+1})+2\bar{f_i} \end{align}
\begin{align}\label{eq:scale1dd}
\gamma''_{2i-1}+\gamma_{2i}=\max(\gamma_{2i-1}+\gamma_{2i+2}, 	C_i+2\beta_{n+i-1})+2\bar{g_i} \end{align}
 \end{lem}

\begin{lem}\label{lem:except3} Consider the arcs $X_{k; j}$ ($2i-1\leq k\leq 2i+2, j=i,i+1$) in Figure \ref{fig:alone}, and let $\chi(i)<1$. Then,
\begin{align}
X_{2i-1; i+1}&=\min\big(\psi_{i+1}, [\chi_i(A)]^+ \big)\quad X_{2i; i+1}=\min\big(\psi_{i+1}, [\chi_i(B)]^+ \big)\\\
X_{2i+1;i }&=\min\big(\psi_{i}, [\chi_{i+1}(A)]^+ \big)\; \quad X_{2i+2; i}=\min\big(\psi_{i}, [\chi_{i+1}(B)]^+ \big)\
\end{align}

\end{lem}
\begin{proof}

We prove $X_{2i-1; i+1}=\min\big(\psi_{i+1}, [\chi_i(A)]^+ \big)$. The other equalities can be proved in a symmetric way. The proof is similar to the proof of Lemma \ref{lem:clover2coordinates}, and is based on the following facts:
\begin{itemize}
\item[(1)] $X_{2i-1; i+1}$ increases $\psi_{i+1}$ and $\chi_i(A)$ by $1$. 
\item[(2)] If $X$ is compatible with $X_{2i-1; i+1}$ then $\chi_i(A)<1$.
\end{itemize}
Therefore, by fact $(1)$ if $\psi_{i+1}=0$ or $\chi_i(A)=0$, then $X_{2i-1; i+1}=0$. Similarly, by fact $(2)$ if $\chi(i)>1$ then $X_{2i-1; i+1}=0$. So suppose that  $\psi_{i+1}\neq 0$ and $\chi_i(A)\neq 0$, and that $\chi(i)<1$ which is guaranteed by the assumption of the lemma.  Let us say that an arc $X$ has property $Q$ if it satisfies $\psi_{i+1}\neq 0$,  $\chi_i(A)\neq 0$ and $\chi(X)<1$; and that an arc is compatible with property $Q$ if it is compatible with an arc that satisfies property $Q$. Figure \ref{fig:proofarc1} illustrates all arcs apart from $X_{2i-1; i+1}$ which are compatible with property $Q$. Since $[\chi_i(A)]^+=0$ for each of these arcs $X_{2i-1; i+1}=\min\big(\psi_{i+1}, [\chi_i(A)]^+ \big)$ by fact $(1)$.\end{proof}

 \begin{figure}[h!]
\begin{center}

  \includegraphics[scale=0.23]{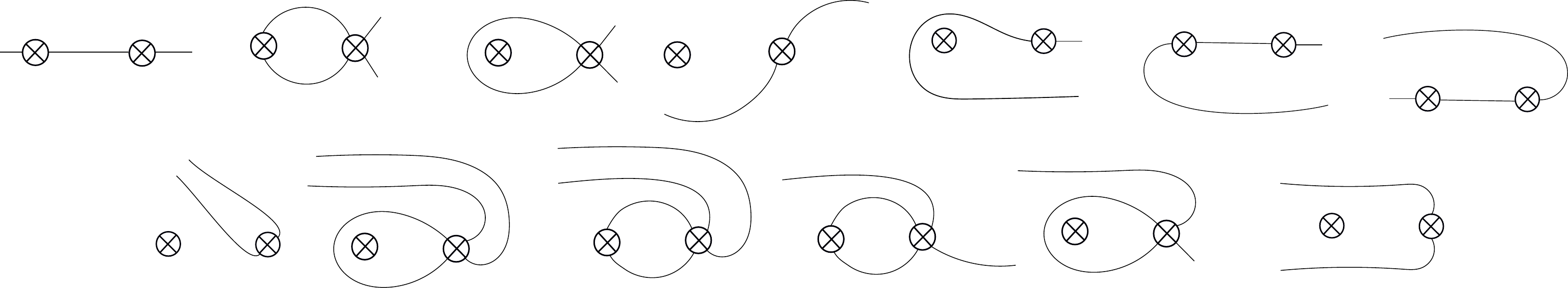}
\caption{Arcs  compatible with property $Q$ }\label{fig:proofarc1}

\end{center}
\end{figure}

\begin{lem}\label{lem:htexep} 
Consider the arcs  $X^{v_i}_{k; j}$ ($2i-1\leq k\leq 2i+2, j=i,i+1$)  in Figure \ref{fig:alone}.  Let
\begin{align*}
M_1&=\min(\Lambda_i, \psi_{i+1}-X_{2i-1; i+1}, \bar{\Delta}_{i+1}(B));~M_2=\min(\Lambda_{i+1}, \psi_{i}-X_{2i+2; i}, \bar{\Delta}_{i}(A))\\
M_3&=\min(\Lambda_i, \psi_{i+1}-X_{2i-1; i+1}, \bar{\Delta}_{i+1}(B));~M_4=\min(\Lambda_{i+1}, \psi_{i}-X_{2i+2; i}, \bar{\Delta}_{i}(A))\end{align*}
\noindent Then

\begin{align}\label{eqs:ex}
 X_{2i+2; i+1}^{(0,1;-1,0)}&=  \min\bigg(M_1,\, \big(c_{i+1}-(c_i-|\chi_{i+1}(B)|)\big)^+,\,\big(c_{i+1}-2[\chi_{i+1}(B)]^+\big)^+\bigg)\\
 X_{2i+1; i+1}^{(0,1;-1,0)}&=\min\bigg(M_4, \big(c_{i+1}-\big(c_i-|(\chi_{i+1}(A)|)\big)^+, \big(c_{i+1}-2[\chi_{i+1}(A)]^+\big)^+\bigg)\\
X_{2i-1; i}^{(1,0;0,1)} &= \min\bigg(M_2, \big(c_{i}-\big(c_{i+1}-|\chi_{i}(A)|)\big)^+, \big(c_{i}-2[\chi_{i}(A)]^+\big)^+\bigg) \\
X_{2i; i}^{(1,0;0,1)}&=\min\bigg(M_3, \big(c_{i}-\big(c_{i+1}-|\chi_{i}(B)|)\big)^+, \big(c_{i}-2[\chi_{i}(B)]^+\big)^+\bigg)
\end{align}
\end{lem}

\begin{proof}
We compute $X_{2i+2; i+1}^{(0,1;-1,0)}$ which is a standard exceptional arc with respect to a clover of type II. Again, the other equalities can be proved in a symmetric way.  To compute this arc separately  we need modification on the formulae given in Lemma \ref{lem:clover2coordinates} to eliminate the values $s'_i$, $r'_i$, $z'_i$ which are parameters related with exceptional arc systems of type II, and the number of highly twisted exceptional arcs in the set $[X_{2i+2; i+1}]$. Using  the value $X_{2i-1; i+1}- \psi_{i+1}$ in $M_1$ rules out the possibility of scissors and hence guarantees that $s'_i=0$. Similarly, since $\big(c_{i+1}-(c_i-|\chi_{i+1}(B)|)\big)^+=0$ for  anchors and ribbons we get $r'_i=z'_i=0$. Finally, for each highly twisted exceptional arc in $[X_{2i+2; i+1}^{(0,1;-1,0)}]$ we have  $\min\big( \big(c_{i+1}-(c_i-|\chi_{i+1}(B)|)\big)^+,\,\big(c_{i+1}-2[\chi_{i+1}(B)]^+\big)^+\big)=0$ (see for instance $l_{19}$ and $l_{20}$).  Since $X_{2i+2; i+1}^{(0,1;-1,0)}$ increases $M_1, \big(c_{i+1}-(c_i-|\chi_{i+1}(B)|)\big)^+$ and $\big(c_{i+1}-2[\chi_{i+1}(B)]^+\big)^+$ by one we conclude that $X_{2i+2; i+1}^{(0,1;-1,0)}$ is as given in equation (\ref{eqs:ex}).
\end{proof}

% \noindent If $L_i$ is mixed  $X_{2i+2; i+1}^{(0,1;-1,0)}=X_{2i; i}^{(1,0;0,1)}=X_{2i; i}^{(1,0;0,1)}=X_{2i+1; i+1}^{(0,1;-1,0)}=0$ by Remark \ref{rem:poss}. Otherwise they are as given in Lemma \ref{lem:htexep}.

\begin{lem}\label{lem:psi'}
Let $\psi_{k}(u_i(L))=\psi'_{k}$ and  $\psi_{k}(u^{-1}_i(L))=\psi''_{k}$ ($k=i, i+1$) denote the number of straight components of $u_i(L)\cap S'_k$ and $u^{-1}_i(L)\cap S'_k$ respectively. 

\noindent Let $L_i$ be negative. Then
\begin{equation*}
\psi'_{i+1}= \begin{cases}
             \lambda^{-}_{c_i}+\chi_{i+1}(B)  & \text{if } \chi_{i+1}(B))\leq 0\\
             (\Delta_{i+1}(B)-\lambda_i^-)^+& \text{if } \chi_{i+1}(B))>0
       \end{cases} \quad \text{and} \quad
\psi'_{i}= \begin{cases}
             \lambda^{+}_{c_{i+1}}+\chi_{i}(A)& \text{if } \chi_{i}(A)\leq 0\\
             (\Delta_{i}(A)-\lambda^+_{i+1})^+& \text{if } \chi_{i}(A)>0
       \end{cases}
\end{equation*}
Let $L_i$ be positive. Then,
\begin{equation*}
\psi''_{i+1}= \begin{cases}
             \lambda^{-}_{c_i}+\chi_{i+1}(A) & \text{if } \chi_{i+1}(A))\leq 0\\\
             (\Delta_{i+1}(A)-\lambda_i^-)^+& \text{if } \chi_{i+1}(A))>0
       \end{cases} \quad \text{and} \quad
\psi''_{i}= \begin{cases}
          \lambda^{+}_{c_{i+1}}+\chi_{i}(B)& \text{if } \chi_{i}(B)\leq 0\\
             (\Delta_{i}(B)-\lambda^+_{i+1})^+& \text{if } \chi_{i}(B)>0
       \end{cases}
\end{equation*}

\end{lem}
\begin{proof}
To compute the number $\psi'_{i+1}$ of straight components  of $u_i(L)\cap S'_{i+1}$ we need to determine which arcs are $u_i$-straight in $S'_{i+1}$. In order to do this, we first list all standard arcs which are straight  in $S'_{i+1}$ (there are finitely many of those) and take their inverse images under $u_i$ from which we obtain the arcs depicted in Figure \ref{fig:straightproof}. Using Notation \ref{not:twist}  we write the following facts:

 \begin{figure}
\begin{center}
\includegraphics[scale=0.23]{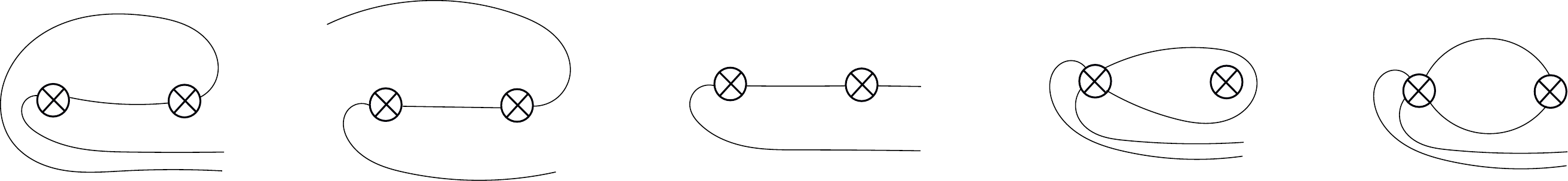}
\caption{ $u_i$-straight arcs in $S'_{i+1}$}\label{fig:straightproof}

\end{center}
\end{figure}
\begin{itemize}

\item[(1)] Each $u_i$-straight arc $X$ is negative with $\widehat{\chi}_{i+1}(B)=0$ and $\lambda^-_{c_i}(X)\neq 0$ (i.e. $X$ has a left core-loop in $S'_i$). The converse is also true.
\item[(2)] $\psi_{i+1}(u_i(X))$ equals the number of left core loops 
of $X\cap S'_i$ which are entirely contained in below components of $X\cap S'_{i+1}$.
\end{itemize}
We have the following cases:

\begin{itemize}

\item If  $\chi_{i+1}(B)<0$, $L_i$ contains $X^{(1, 0; -1, 0)}_{2i+22i+2}$  which satisfies $\lambda^-_{c_i}\neq 0$, is not $u_i$-straight and not compatible with any highly twisted component. Therefore, by (1) and (2) we obtain 
$\psi'_{i+1}=\lambda^-_{c_i}-X^{(1, 0; -1, 0)}_{2i+22i+2}$. It is easy to show that 
$X^{(1, 0; -1, 0)}_{2i+22i+2}=-\chi_{i+1}(B)$ from which we get  
$\psi'_{i+1}=\lambda^{-}_{c_i}+\chi_{i+1}(B)$. Clearly, if $\chi_{i+1}(B)=0$, $L_i$ is some collection of arcs depicted in Figure  \ref{fig:straightproof} each of which satisfies
 $\psi'_{i+1}=\lambda^{-}_{c_i}=(\Delta_{i+1}(B)-\lambda_i^-)^+$ by (2).
 
 \item  If $\chi_{i+1}(B)>0$, $L_i$ contains a highly twisted component $X$, and only left core loops of $X\cap S_{i}$ that do not join right loop components of $X\cap S_{i+1}$ (i.e. those that are contained in  $\Delta_{i+1}(B)$) can be mapped to a straight component of $X\cap S_{i+1}$. That is for each such arc we have $\psi'_{i+1}=(\Delta_{i+1}(B)-\lambda_i^-)^+$. 
 \end{itemize}

 % See e-mail 7 February
%\begin{itemize}
%\item[(1)]  Only a left core loop in $S'_i$ can be mapped to a straight component in $S'_{i+1}$. Therefore, if $\psi_{i+1}(X')\neq 0$, then  $\lambda^-_{c_i}(X)\neq 0$.
%\item[(2)] If $\chi_{i+1}(B)(X)<0$, $\psi'_{i+1}(X)=0$.
% \end{itemize}
%For each arc $X$ in Figure \ref{},  $\psi'_{i+1}(X)$ equals the number or core loop components of $X\cap S'_i$ which are entirely contained in below components of $S'_{i+1}$.  Assuming $\chi_{i+1}(B)\geq 0$ and $\lambda^-_{c_i}\neq 0$
%Since $L$ is a disjoint union of such arcs we get 
\end{proof}

\begin{proof}[Proof of Theorem \ref{thm:update}]
Let $n>0$ and $\mathcal{A}_{g,n}$ denote the set of arcs in Figure \ref{arcsproof}.   $\MCG(N_{g,n})$ acts on both $\mathcal{A}_{g,n}$ and~$\mathfrak{L}_{g,n}$,  and hence $i(\cL, \xi)=i(\delta(\cL),\delta(\xi))$ for any $\delta \in \MCG(N_{g,n})$ and $\xi\in \cA_{g,n}$.  We also recall that  the arcs $\alpha_i$ ($1\leq i\leq 2n-2$) and $\beta_i$ ($1\leq i\leq n$) are not affected by crosscap transpositions. For the crosscap transposition $u_i$, our approach is to compute the number of intersections of $\gamma'_j=u^{-1}_i(\gamma_j)$ $(1\leq j\leq 2g-2)$ and $\beta'_j=u^{-1}_i(\beta_j)$ $(1\leq j\leq n+g-1)$  with $\cL$ instead of computing the number of intersections of $u_i(\cL)$ with $\gamma_j$ and $\beta_j$.  We have,
$$t'_j=\frac{\gamma'_{2j}-\gamma'_{2j-1}}{2}  \quad\text{and}  \quad b'_j=\frac{\beta'_{j}-\beta'_{j+1}}{2}.$$
We shall make use of clover and scale equalities given in Lemma \ref{lem:clover1}, Lemma \ref{lem:clover2}, Lemma \ref{lem:scales1}.   For computational convenience we set $T_{j}=2t_{j}$ ($1\leq j\leq g+n-2$), $B_{j}=2b_{j}$ ($n\leq j\leq g+n-2$), $2\beta_{n+j}=\mathfrak{B}_{n+j}~(1\leq j\leq g-1)$ and $D_j=2d_j$, $E_j=2e_j$, $F_j=2f_j$, $G_j=2g_j$  ($1\leq j\leq g-1$), and work in the max-plus semiring as indicated in Remark \ref{not:semiring}.  Therefore, 

\begin{align}\label{eq:tb}T_j=\left[\frac{\gamma_{2j}}{\gamma_{2j-1}}\right]~\text{and}~B_j=\left[\frac{\beta_{j}}{\beta_{j+1}}\right].\end{align}
and from the clover of type I equality  (\ref{eq:clover1}), 

\begin{align}\label{eq:C}\cC_i=\left[\frac{D_i(\gamma_{2i-1}\gamma_{2i+2}+\gamma_{2i+1}\gamma_{2i})}{\mathfrak{B}_{n+i}} \right]\end{align}

We now consider the two separate cases of the statement.

\begin{itemize}
\item  Suppose that  $1\leq i < g+n-2$/ Observe that $\beta'_j=\beta_j$ for $j\neq n+i$ and $\gamma'_j=\gamma_j$ for $j<2i-1$ and $j>2i+2$. Therefore, $T'_j=T_j$ for except $j=i$ and $j=i+1$; and $B'_j=B_j$ for except $j=n+i-1$ and $j=n+i$. Next we compute $T'_i$, $T'_{i+1}$, $B'_{n+i-1}$ and  $B'_{n+i}$.

\begin{enumerate}
\item\label{birinciupdate}
We shall first compute $T'_{i+1}=\left[\frac{\gamma'_{2i+2}}{\gamma'_{2i+1}}\right]$. We have $\gamma'_{2i+2}=\gamma_{2i}$. To compute $\gamma'_{2i+1}$ we use the scale of type I equality (\ref{eq:scale1aa}) and obtain
\begin{align}\label{eq:gammap2i+1}
\gamma'_{2i+1}= \left[\frac{F_i(\cC_i\mathfrak{B}_{n+i+1}+\gamma_{2i}\gamma_{2i+1})}{\gamma_{2i+2}}\right].
\end{align}

\noindent Then from (\ref{eq:tb}), (\ref{eq:C}) and (\ref{eq:gammap2i+1}) we compute that

\begin{eqnarray*}
\left[\frac{1}{T'_{i+1}}\right]&=&\left[F_i\Bigg(\frac{D_i(\gamma_{2i-1}\gamma_{2i+2}+\gamma_{2i+1}\gamma_{2i})}{\mathfrak{B}_{n+i}\gamma_{2i}\gamma_{2i+2}}\mathfrak{B}_{n+i+1}+\frac{\gamma_{2i}\gamma_{2i+1}}{\gamma_{2i}\gamma_{2i+2}}\Bigg)\right]
%&=&\left[\frac{\Pi_iX_iD_iE_i}{T_{i}B^2_{n+i}}+\frac{\Pi_iX_iD_iE_i}{T_{i+1}B^2_{n+i}}+\frac{\Pi_i}{T_{i+1}}\right]\\
\text{and hence,}\\
T'_{i+1}&=&\left[\frac{T_{i}T_{i+1}B^2_{n+i}}{F_i\big(T_i(D_{i}+B^2_{n+i})+D_iT_{i+1}\big)}\right]. \text{That is,}\\
T'_{i+1}&=&T_{i}+T_{i+1}+2B_{n+i}-\big(F_i+\max\big(T_i+\max(D_i, 2B_{n+i}), D_i+T_{i+1})\big)\end{eqnarray*}
Dividing both sides of the equation by $2 $ we get  $$t'_{i+1}=\left[\frac{t_{i}t_{i+1}B_{n+i}}{f_i\big(t_i(d_{i}+B_{n+i})+ d_it_{i+1}\big)}\right]=\left[\frac{t_{i}t_{i+1}b^2_{n+i}}{f_i\big(t_i(d_{i}+b^2_{n+i})+ d_it_{i+1})}\right].$$

\item\label{ikinciupdate}
We shall now compute $T'_{i}=\left[\frac{\gamma'_{2i}}{\gamma'_{2i-1}}\right]$. We have $\gamma'_{2i-1}=\gamma_{2i+1}$. To compute $\gamma'_{2i}$  we use the scale of type II equality (\ref{eq:scale1bb}) and obtain
$$
\gamma'_{2i}=\left[\frac{G_i(\gamma_{2i}\gamma_{2i+1}+C_i\mathfrak{B_{n+i-1}})}{\gamma_{2i-1}}\right].
$$

Hence, from (\ref{eq:C}) we get $$t'_{i}=\left[g_i\big(t_{i}(1+d_ib^2_{n+i-1})+d_ib^2_{n+i-1}t_{i+1}\big)\right].$$

\item \label{3i} We proceed with $B'_{n+i}=\left[\frac{\beta'_{n+i}}{\beta'_{n+i+1}}\right]$. We have $\beta'_{n+i+1}=\beta_{n+i+1}$ and from the clover of type II equality (\ref{eq:clover2}), 

%$$ 2\beta'_{n+i}+C_i=\max(\gamma'_{2i}+\gamma'_{2i+1}, \gamma_{2i}+\gamma_{2i+1})+X+D_i+E_i$$
%by Lemma \ref{lem:quadrilateral1}. That is,

$$\mathfrak{B}'_{n+i}=\left[\frac{E_i(\gamma'_{2i}\gamma'_{2i+1}+\gamma_{2i}\gamma_{2i+1})}{C_i}\right].$$
 Since $\gamma'_{2i-1}=\gamma_{2i+1}$ and $\gamma'_{2i+2}=\gamma_{2i}$  \begin{align*}
\gamma'_{2i}=\left[T'_{i}\gamma_{2i+1}\right]~~\text{and}~~\gamma'_{2i+1}=\left[\frac{\gamma_{2i}}{T'_{i+1}}\right].
\end{align*}

%\begin{align*}
%B'_{i}=\left[\frac{\alpha'_{2i-2}\alpha'_{2i-1}+\alpha_{2i-2}\alpha_{2i-1}}{u_i\beta_{i+1}}\right].
%\end{align*}

%Recalling that $$C_i=\left[\frac{XD_iE_i(\gamma_{2i-1}\gamma_{2i+2}+\gamma_{2i+1}\gamma_{2i})}{\beta^2_{n+i}} \right]$$

%we obtain,
\begin{align*}
B'_{n+i}&=\left[\frac{E_i\big(\frac{T'_{i}}{T'_{i+1}}\gamma_{2i}\gamma_{2i+1}+\gamma_{2i}\gamma_{2i+1}\big)}{C_i\mathfrak{B}_{n+i+1}}\right]=\left[\frac{E_i\gamma_{2i}\gamma_{2i+1}(\frac{T'_{i}}{T'_{i+1}}+1)}{\frac{D_i\mathfrak{B}_{n+i+1}}{\mathfrak{B}_{n+i}}(\gamma_{2i-1}\gamma_{2i+2}+\gamma_{2i}\gamma_{2i+1})}\right]\\
&=\left[\frac{E_i}{D_i}B^2_{n+i}(\frac{T'_{i}+T'_{i+1}}{T'_{i+1}})\frac{1}{\frac{\gamma_{2i-1}\gamma_{2i+2}+\gamma_{2i}\gamma_{2i+1}}{\gamma_{2i}\gamma_{2i+1}}}\right]=\left[\frac{E_i}{D_i}B^2_{n+i}(\frac{T'_{i}+T'_{i+1}}{T'_{i+1}})\frac{T_{i}}{T_{i}+T_{i+1}}\right]
%&=\left[\frac{A_{i-1}(1+B_{i-1})(1+B_i) + A_iB_{i-1}}{A_i}\right]
\end{align*}
\noindent from which we get $$b'_{n+i}=\left[\frac{e_i}{d_i}b^2_{n+i}(\frac{t'_{i}+t'_{i+1}}{t'_{i+1}})\frac{t_{i}}{t_{i}+t_{i+1}}\right].$$

 \item Now we shall  compute $B'_{n+i-1}=\left[\frac{\beta'_{n+i-1}}{\beta'_{n+i}}\right]$. We have
 $$\beta'_{n+i-1}=\beta_{n+i-1},~\beta'_{n+i+1}=\beta_{n+i+1}~\text{and}~B'_{n+i}=\left[\frac{\beta'_{n+i}}{\beta'_{n+i+1}}\right].$$

 Therefore, $B'_{n+i-1}=\left[\frac{\beta_{n+i-1}}{B'_{n+i}\beta_{n+i+1}}\right]$. Multiplying the numerator and denominator by $\beta_{n+i}$ gives $B'_{n+i-1}=\left[\frac{B_{n+i}B_{n+i-1}}{B'_{n+i}}\right]$. That is,  $$b'_{n+i-1}=\left[\frac{b_{n+i}b_{n+i-1}}{b'_{n+i}}\right].$$

 \end{enumerate}

\item  Now, suppose that $ i=g-1$ (Figure \ref{arcsproof2}). Observe as before that $t'_j=t_j$ for all $j<g-1$ and $b'_j=b_j$ for all $j<n+g-2$. Since there are no teardrops encircling the last crosscap, our approach to compute $t'_{g-1}$ and $b'_{n+g-2}$ is to add  dummy teardrops $\gamma_{2g-1}$,~$\gamma_{2g}$ and $\beta_{n+g}$ as depicted in Figure~\ref{fig:action}, which enables us to make similar calculations as in the previous statement. We first note that $\gamma_{2g-1}=\gamma_{2g}=\frac{\beta_{n+g-1}}{2}$ and $\beta_{n+g}=0$ hence we have $T_g=0$ and $B_{n+g-1}=\beta_{n+g-1}$. Similar calculations give

$$t_{g-1}' = \tropical{g_{g-1}\big(t_{g-1}+d_{g-1}B_{n+g-2}(1+t_{g-1})\big)} \quad \text{and}\quad B_{n+g-2}' = \tropical{ \frac{d_{g-1}}{e_{g-1}}B_{n+g-2}\frac{1+t_{g-1}}{t_{g-1}(1+t'_{g-1})} }$$

%  \begin{aligned}
%  B_{n+g-2}'' &= \tropical{ \frac{d_{g-2}}{\bar{e}_{g-2}}B_{n+g-2}\frac{t''_{g-1}(1+t_{g-1})}{1+t''_{g-1}}}.
%t_{g-1}'' &= \tropical{\frac{t_{g-1}}{\bar{g}_{g-1}\big(d_{g-1}B_{n+g-2}(1+t_{g-1})\big)} } \end{aligned} 

  \begin{figure}[h!]
\begin{center}
\labellist
\small\hair 2pt

 \pinlabel {\begin{turn}{-90}$\scriptstyle{\gamma_{2g-1}}$\end{turn}} [ ] at 355 220
 \pinlabel {\begin{turn}{-90}$\scriptstyle{\gamma_{2g}}$\end{turn}} [ ] at 355 60

 \pinlabel {\begin{turn}{-90}$\scriptstyle{\gamma_{2g-3}}$\end{turn}} [ ] at 275 230
  \pinlabel {\begin{turn}{-90}$\scriptstyle{\gamma_{2g-2}}$\end{turn}} [ ] at 275 60

 \pinlabel {\begin{turn}{-90}$\scriptstyle{\beta_{n+g-1}}$\end{turn}} [ ] at 310 200
 \pinlabel {\begin{turn}{-90}$\scriptstyle{\beta_{n+g-2}}$\end{turn}} [ ] at 205 200

    \pinlabel {\begin{turn}{-90}$\scriptstyle{\beta_{n+g}}$\end{turn}} [ ] at 395 200

 \endlabellist  
 \includegraphics[scale=0.4]{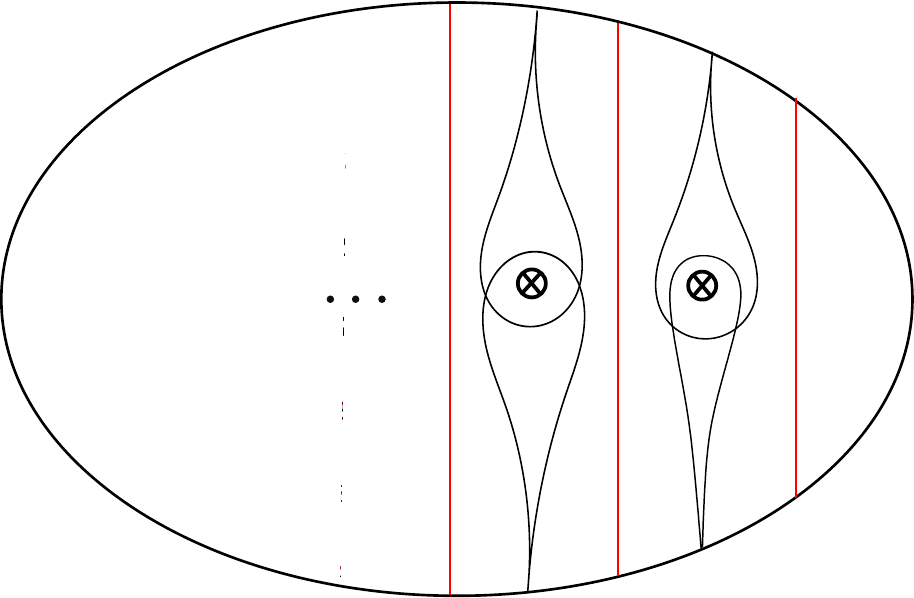}
\caption{Dummy teardrops are used to compute $t'_{g-1}$ and $b'_{n+g-1}$}\label{fig:action}
\end{center}
\end{figure}

Now let $n=0$. The formulae for $1<i\leq g-2$ are obtained similarly replacing $i$ with $i-1$. For $i=1$, we add two dummy punctures around the first crosscap. Similar arguments  give that

  $$t_{1}' =\left[\frac{t_{1}B_{1}}{f_0(1+B_1)+ t_1}\right] \quad \text{and} \quad B_{1}' = \left[B_{1}e_0\frac{(1+t'_1)}{t_1'}\frac{1}{(1+t_1)}\right]$$
  
For $i=g-1$ we note that rotation through~$\pi$ about the center of the surface conjugates each
crosscap generator~$u_i$ to $u_{g-i}$ and the corresponding
  transformation of generalized Dynnikov coordinates in max-plus notation, is given by
\[(t_1,\ldots, t_{g-2},b_1,\ldots,b_{g-2}) \mapsto
\tropical{(1/t_{g-2},\ldots, 1/t_1,1/b_{g-2},\ldots, 1/b_1)}\]
\noindent hence we get
\begin{align*}t_{g-2}' = \tropical{g_{g-2}\big(t_{g-2}+d_{g-2}B_{g-2}(1+t_{g-2})\big)} \quad \text{and} \quad B_{g-2}' = \tropical{ \frac{d_{g-2}}{e_{g-2}}B_{g-2}\frac{1+t_{g-2}}{t_{g-2}(1+t'_{g-2})} }\end{align*}

\end{itemize}

By Remark \ref{reflectionhori} we obtain the rules for $t''_{i}$~and~$b''_{i}$ for each case by symmetry, conjugating the rules for $u_i$ by the involution (\ref{reflectionhori}).\end{proof}

%
%\begin{remark} 
%We note that 
%\end{remark} 
%\begin{example}\label{ex:example}
%Consider the curve $\cL$ in Figure \ref{fig:example1} with $\rho(\cL)=(0,1;-1, -1; 1,1,1)$.  We compute $\rho(u^{-1}_2(\cL))$ in terms of $\rho(\cL)$ applying the formulae in Theorem \ref{thm:update}. We have $c''_1=c_1$, $c''_2=c_3$ and $c''_3=c_2$. Also $t''_1=t_1$ and $b''_1=b_1$. We show that  $t''_2=1$ and $b''_2=0$ using equalities (\ref{eq:updatespecial}). We need the exceptional parameters $d_1$, $\bar{e}_1$ and $\bar{g}_1$. By equalities (\ref{eq:di}),  (\ref{eq:eibar}) and (\ref{eq:gibar}) we compute $d_1=0$, $\bar{e}_1=0$ and $\bar{g_1}=1$. So $\rho(u^{-1}_2(\cL))=()$
%
%\begin{align*}
%t''_2&=\tropical{\frac{t_1}{\bar{g}_1d_1B_1(1+t_1)}}=t_1-(\bar{g_1}+d_1+B_1+\max(0,t_1))=1\\
%B_{2}'' &= \tropical{ \frac{d_{1}}{\bar{e}_{1}}B_{1}\frac{t''_{1}(1+t_{1})}{1+t''_{1}}}=d_1-\bar{e}_1+B_1+t''_1+\max(0, t_1)-\max(0,t''_1)=0.
%\end{align*}
%
%  \begin{figure}[h!]
%\begin{center}
%\labellist
%\small\hair 2pt
%
%
%
% \pinlabel {$\scriptstyle{u^{-1}_{2}}$} [ ] at 270 60
%
%\endlabellist  
% \includegraphics[scale=0.4]{exampleaction}
%\caption{Action of $u^{-1}_2$ on $C$ }\label{fig:action}
%\end{center}
%\end{figure}
%
%\end{example}

\begin{remark}
 Note that the method introduced in this paper can be used to provide an efficient way to solve on non--orientable surfaces many combinatorial and dynamical problems \cite{yurttasyeni2}  that were previously solved only on orientable surfaces before  \cite{bell, paper1, paper3, paper6, paper5, paper2}. 
However, to solve such problems not only for sequences of crosscap transpositions but any element of the mapping class group we need to  describe the action of the mapping class group $\MCG(N_{g,n})$  on $\mathfrak{L}_{g,n}$   in terms of generalized Dynnikov coordinates \cite{yurttasyeni1}. That is we need to  compute the action of the other generators of $\MCG(N_{g,n})$ which are crosscap slides, puncture slides and Dehn twists about certain $2$-sided curves \cite{korkmaz} in terms of generalized Dynnikov coordinates \cite{yurttasyeni1}    which require similar techniques introduced in this paper.
\end{remark}

\begin{ack}
\textnormal{This work was completed during a visit of the author at  Columbia University as a Fulbright scholar. The author would like to thank the Fulbright Scholar Program for their support and  Columbia University for their warm hospitality. }

\end{ack}
%Then, it was studied in [11, 12] as an efficient method for a solution of the word problem of Bn, that is, the problem of determining whether a given braid β ∈ Bn, given in terms of the Artin gen- erators, is the identity. (Here, and throughout the thesis, we use efficient in an informal sense and do not carry out any formal efficiency analysis.)


\begin{thebibliography}{99}
\bibitem{dynn1} Dynnikov, I. \textit {On a {Y}ang-{B}axter mapping and the {D}ehornoy ordering.} Uspekhi Mat. Nauk,  \textbf{57(3(345))}, 151-152,  2002.
\bibitem{bell} Bell, M. \textit{Simplifying triangulations}, arXiv: https://arxiv.org/abs/1604.04314, 2016.
\bibitem{fathi} Fathi, A. and Laudenbach. F, and Poenaru. V. \textit{Travaux de Thurston sur les surfaces}, volume 66 of Ast\'erisque. Soci\'et\'e Math\'ematique de France, Paris, 1979. Seminair\'e Orsay.
\bibitem{paper1} Hall, T. and Yurtta{\c{s}}, S.~{\"O}. \textit{On the topological entropy of families of braids.} Topology Appl. \textbf{156(8)}, 1554-1564, 2009.
\bibitem{paper3} Hall, T. and Yurtta{\c{s}}, S.~{\"O}. \textit{Counting components of an integral lamination.} manuscripta mathematica \textbf{153(1)}, 263-278, 2017.
\bibitem{paper6} Hall, T. and Yurtta{\c{s}}, S.~{\"O}. \textit{Intersections of multicurves from Dynnikov coordinates.} Bull. Aust. Math. Soc., \textbf{98}, 149-158,  2018.

\bibitem{korkmaz} Korkmaz, M. \textit{Mapping Class Groups of Nonorientable Surfaces.} Geometriae Dedicata, \textbf{89(1)}, 107--131, 2002.
\bibitem{pamukyurttas} Pamuk, M. and  Yurtta{\c{s}} S.~{\"O}. \textit{Integral laminations on non--orientable surfaces.} Turkish J. Math., \textbf{42}, 69--82,  2018.
\bibitem{papapenner} Papadopoulos, A. and Penner, R. C. \textit{Hyperbolic metrics, measured foliations and pants decompositions for non-orientable surfaces.} Asian J. Math, \textbf{20}, 157-182,  2016.
\bibitem{penner} Penner, R. C. and  Harer, J. L. \textit{Combinatorics of train tracks}, volume 125 of Annals of Mathematics Studies. Princeton University Press, Princeton, NJ, 1992.
\bibitem{parlak} Parlak, Anna and Stukow, Micha\l. \textit{Roots of crosscap slides and crosscap transpositions.} Periodica Mathematica Hungarica. \textbf{75(2)}, 413--419,  2017.

		
 \bibitem{thurston} Thurston, D. \textit{Geometric intersection of curves on surfaces.}Preprint available from  https://dpthurst.pages.iu.edu/DehnCoordinates.pdf

\bibitem{thurstonw} Thurston, W.P. \textit{On the geometry and dynamics of diffeomorphisms of surfaces}. Bull. Amer. Math. Soc. (N.S.), \textbf{19(2)}:417–431, 1988.
\bibitem{paper2} Yurtta{\c{s}}, S.~{\"O}.  \textit{Geometric intersection of curves on punctured disks.}  Journal of the Mathematical Society of Japan, \textbf{65(4)}, 1554-1564, 2013.
\bibitem{paper5} Yurtta{\c{s}}, S.~{\"O}. \textit{Dynnikov and train track transition matrices of pseudo-Anosov braids}, Discrete and Continuous Dynamical Systems,  \textbf{36(1)}, 541-570, 2016.
\bibitem{yurttasyeni1} Yurtta{\c{s}}, S.~{\"O}. \textit{Action of $y$-homeomorphisms and Dehn twists on non--orientable surfaces}, in preparation.
\bibitem{yurttasyeni2} Yurtta{\c{s}}, S.~{\"O}. \textit{Algorithms for curves on non--orientable surfaces}, in preparation.

\end{thebibliography}
\end{document}